\documentclass[10pt]{IEEEtran}
\usepackage{graphicx}
\usepackage{amsmath,amssymb}
\usepackage{mathrsfs}
\usepackage{hyperref}                    
\usepackage[dvips]{epsfig}
\usepackage{multirow}
\usepackage[noadjust]{cite}
\usepackage{epstopdf}
\usepackage{color}
\usepackage[color=lightgray]{todonotes}
\usepackage{url}
\usepackage{booktabs}
\usepackage{tabularx}
\usepackage{abraces}
\usepackage{hyperref}
\hypersetup{
    colorlinks=true,
    linkcolor=black,
    filecolor=black,      
    urlcolor=blue,
}
\usepackage{enumitem}

\newcommand{\real}{\ensuremath{\mathbb{R}}}

\newcommand{\A}{\ensuremath{\mathcal{A}}}
\newcommand{\smat}[1]{\ensuremath{\left[ 
\begin{smallmatrix} #1 \end{smallmatrix}
\right]
}}
\newcommand{\C}{\ensuremath{\mathsf{c}}}

\DeclareMathOperator{\dom}{dom}
\DeclareMathOperator{\rge}{rge}

\makeatletter
\newcommand{\customlabel}[2]{
\protected@write \@auxout {}{\string \newlabel {#1}{{#2}{}}}}
\makeatother

\newtheorem{thm}{\textbf{Theorem}}
\newenvironment{theorem}{\begin{thm}\itshape}{\end{thm}}
\newtheorem{defi}{\textbf{Definition}}
\newenvironment{definition}{\begin{defi}\itshape}{\end{defi}}
\newtheorem{propo}{\textbf{Proposition}}
\newenvironment{proposition}{\begin{propo}\itshape}{\end{propo}}
\newtheorem{assu}{\textbf{Assumption}}
\newenvironment{assumption}{\begin{assu}\itshape}{\end{assu}}
\newtheorem{lem}{\textbf{Lemma}}
\newenvironment{lemma}{\begin{lem}\itshape}{\end{lem}}
\newtheorem{fact}{\textbf{Fact}}
\newtheorem{rem}{\textbf{Remark}}

 \newenvironment{proof}{\noindent\textit{Proof.}}{\hfill $\square$\\}

\title{Obstacle Avoidance via Hybrid Feedback}
\author{S. Berkane,~\IEEEmembership{Member,~IEEE,} A. Bisoffi,~\IEEEmembership{Member,~IEEE}, D. V. Dimarogonas,~\IEEEmembership{Senior Member,~IEEE}
\thanks{This research work is partially supported by NSERC-DG RGPIN-2020-04759, the European Research Council (ERC), the EU H2020 Co4Robots, the SSF COIN project, the Swedish Research Council (VR) and the Knut och Alice Wallenberg Foundation. S. Berkane (\texttt{soulaimane.berkane@uqo.ca}) is with the Department of Computer Science and Engineering, University of Quebec in Outaouais, Gatineau, Canada. A. Bisoffi (\texttt{a.bisoffi@rug.nl}) is with ENTEG and the J.C. Willems Center
for Systems and Control, University of Groningen, Groningen,
The Netherlands. D. V. Dimarogonas (\texttt{dimos@kth.se}) is with the Division of Decision and Control Systems, KTH Royal Institute of Technology, Stockholm, Sweden. }}
\allowdisplaybreaks

\begin{document}
\maketitle
\begin{abstract}
In this paper we present a hybrid feedback approach to solve the navigation problem of a point mass in the $n-$dimensional space containing an arbitrary number of ellipsoidal shape obstacles. The proposed hybrid control algorithm guarantees both global asymptotic stabilization to a reference and avoidance of the obstacles. The intuitive idea of the proposed hybrid feedback is to switch between two modes of control: stabilization and avoidance. The geometric construction of the flow and jump sets of the proposed hybrid controller, exploiting hysteresis regions, guarantees Zeno-free switching between the stabilization and the avoidance modes. Simulation results illustrate the performance of the proposed hybrid control approach for 2-dimensional and 3-dimensional scenarios.
\end{abstract}
\section{Introduction}
For decades, the obstacle avoidance problem has been an active area of research in the robotics and control communities \cite{hoy2015algorithms}. In a typical robot navigation scenario, the robot is required to reach a given goal (destination) while not colliding with a set of obstacle regions in the workspace. 
Since the pioneering work by Khatib \cite{khatib1986real}, artificial potential fields have been widely used in the obstacle avoidance problem  since  they  offer  the  possibility  to  combine the solution  to  the  global find-path problem  with  a feedback controller for the robot, thus, allowing the high-level planner to address more abstract tasks. The idea is to generate an artificial potential field that renders the goal attractive and the obstacles repulsive. Then, by considering trajectories that navigate along the negative gradient of the artificial potential field, one can ensure that the robot will reach the desired target while avoiding to collide with the obstacles. However, artificial potential field-based algorithms suffer from 1) the presence of local minima preventing the successful navigation to the target point and 2) arbitrarily large repulsive potential near the obstacles which are in conflict with the inevitable actuator saturations. 

The navigation function-based approach, which was initiated by Koditscheck and Rimon \cite{koditschek1990robot} for sphere worlds \cite[p. 414]{koditschek1990robot}, solves both problems. It allows to obtain artificial potential fields with the nice property that all but one of the critical points are saddles with the remaining critical point being the desired reference. Since then, the navigation function-based approach has been extended in many different directions; {\it e.g.,} for multi-agent systems \cite{tanner2005formation,dimarogonas2006feedback,roussos2013decentralized}, for unknown sphere words \cite{lionis2007locally}, and for focally admissible obstacles \cite{filippidis2013navigation}. The major drawback of navigation functions is that they are not correct by construction. In fact, navigation functions are theoretically guaranteed to exist, but their explicit computation is not straightforward since they require an unknown tuning of a given parameter to eliminate local minima.

Recently, Loizou \cite{loizou2017navigation} introduced the navigation transform that diffeomorphically maps the workspace to a trivial domain called the {\it point world} consisting of a closed ball with a finite number of points removed. Once this transformation is found, the navigation problem is solved from almost all initial conditions without requiring any tuning. In addition, the trajectory duration is explicitly available, which provides a timed-abstraction solution to the motion-planning problem. Similarly, the recent work in~\cite{vrohidis2018prescribed} uses the so-called prescribed performance control to design a time-varying control law that drives the robot, in finite time, from all initial conditions to some neighborhood  of the target while avoiding the obstacles. Another approach to the navigation problem is through barrier functions (see \cite{ames2017control} and references therein), which are developed for nonlinear systems with state-space constraints and ensure safety. Model predictive control approaches have been also used for reactive robot navigation, e.g., \cite{falcone2006real,defoort2009motion}.

However, by using any of the approaches described above, it is not possible to ensure safety from all initial conditions in the obstacle-free state space. As pointed out in \cite{koditschek1990robot}, the appearance of additional undesired equilibria is unavoidable when considering continuous time-invariant vector fields. This is a well known topological obstruction to global asymptotic stabilization by continuous feedback when the free state space is not diffeomorphic to an Euclidean space (see, e.g., \cite[Thm.~2.2]{wilson1967structure}).  Furthermore, this problem is more far-reaching since, by using a continuous feedback law, it is always possible to find arbitrarily small adversarial (noise) signals acting on the vector field, such that a set of initial conditions different from the target, possibly of measure zero, can be rendered stable \cite[Thm.~6.5]{sanfelice2007robust}. 
To deal with such limitations, the authors in \cite{sanfelice2006robust} proposed a hybrid state feedback controller, using Lyapunov-based hysteresis switching, to achieve robust global asymptotic regulation in $\mathbb{R}^2$ to a target while avoiding a single obstacle. 
This approach has been exploited in \cite{poveda2018hybrid} to steer a planar vehicle to the source of an unknown but measurable signal while avoiding an obstacle. 
In \cite{braun2018unsafe} and \cite{braun2018explicit}, a hybrid control law was proposed to globally asymptotically stabilize a class of linear systems while avoiding neighbourhoods of unsafe isolated points in $\mathbb{R}^n$. Although such hybrid approaches are promising, they are still challenged by constructing the suitable hybrid feedback for higher dimensions and with more complex obstacles shapes.

In this work, we propose a hybrid control algorithm for the global asymptotic stabilization of a point mass moving in an arbitrary $n-$dimensional space while safely avoiding obstacles that have generic ellipsoidal shapes, based on the preliminary treatment of this problem for a single spherical obstacle in \cite{berkane2019hybrid}. The ellipsoids provide a tighter bounding volume than spheres, and in our scheme this volume can be arbitrarily flat and close to the target, which leads to a significant reduction in the level of conservatism compared, for instance, to \cite[Thm.~3]{paternain2018navigation} (as shown in Section~\ref{section:example}). 

Our proposed hybrid algorithm employs a hysteresis-based switching between the avoidance controller and the stabilizing controller in order to guarantee forward invariance of the obstacle-free region (corresponding to safety) and global asymptotic stability of the reference position. We consider trajectories in an $n-$dimensional Euclidean space and we resort to tools from higher-dimensional geometry \cite{meyer2000matrix} to provide a construction of the flow and jump sets where the different modes of operation of the hybrid controller are activated. Furthermore, 
the hybrid control law guarantees a bounded control input, it matches the stabilizing controller in arbitrarily large subsets of the obstacle-free region by a suitable tuning of its parameters (hence qualifying as minimally invasive), it can be readily extended to a non-point mass vehicle and enjoys some level of inherent robustness to perturbations.
\emph{Structure.} Preliminaries are in Section~\ref{section:preliminaries}. The navigation problem is formulated in Section~\ref{section:problem}. Our proposed hybrid control scheme is discussed in Section~\ref{section:controller}. Section \ref{section:main} presents the main result of forward invariance of the obstacle-free space and global asymptotic stability of the target, together with other desirable complementary properties. Numerical examples are in Section~\ref{section:example}. All the proofs are in the Appendix.%

%
%
\section{Preliminaries}
\label{section:preliminaries}
$\mathbb{N}$, $\mathbb{R}$ and $\mathbb{R}_\geq$ denote, respectively, the set of nonnegative integers, reals and nonnegative reals. $\mathbb{R}^n$ is the $n$-dimensional Euclidean space and $\mathbb{S}^n$ is the $n$-dimensional unit sphere embedded in $\mathbb{R}^{n+1}$. Given the column vectors $v_1 \in \real^{n_1}$ and $v_2 \in \real^{n_2}$, $(v_1,v_2)$ denotes the stack vector $\begin{bmatrix}v_1^\top & v_2^\top\end{bmatrix}^\top$. The Euclidean norm of $x\in\mathbb{R}^n$ is defined as $\|x\|:=\sqrt{x^\top x}$ and  the geodesic distance between two points $x$ and $y$ on the sphere $\mathbb{S}^n$ is defined by $\mathbf{d}_{\mathbb{S}^n}(x,y):=\arccos(x^\top y)$ for all $x,y\in\mathbb{S}^n$. For an arbitrary matrix $A\in\mathbb{R}^{n\times n}$, $\lambda_i(A)$ denotes the $i$-th eigenvalue of $A$. If $A$ is a symmetric matrix, then $\lambda_{\min}(A)$ and $\lambda_{\max}(A)$ denote, respectively, the smallest and largest eigenvalues of $A$. Given a closed set $\mathcal{A}\subset\mathbb{R}^n$, we define the distance to the set $\mathcal{A}$ by $|x|_{\mathcal{A}}:=\inf_{y\in\mathcal{A}}\|x-y\|$.
Given two sets $\mathcal{A}$ and $\mathcal{B}$, we define the distance from $\mathcal{A}$ to $\mathcal{B}$ by $\mathbf{dist}(\mathcal{A},\mathcal{B}):=\inf\{|a-b| \colon a \in \mathcal{A}, b \in \mathcal{B} \}$. The closure, interior and boundary of a set $\mathcal{A}\subset\mathbb{R}^n$ are denoted as  $\overline{\mathcal{A}}, \mathcal{A}^\circ$ and $\partial\mathcal{A}$, respectively. The relative complement of a set $\mathcal{B}\subset\mathbb{R}^n$ with respect to a set $\mathcal{A}$ is denoted by $\mathcal{A}\backslash\mathcal{B}$ and contains the elements of $\mathcal{A}$ which are not in $\mathcal{B}$. In particular, we use $\mathcal{A}^\C$ to denote the complement of $\mathcal{A}$ in $\mathbb{R}^n$, {\it i.e.,} $\mathcal{A}^\C=\mathbb{R}^n\backslash\mathcal{A}$. Given the sets $\mathcal{A}$, $\mathcal{B}$ and $\mathcal{C}$, the following set identities \cite{engelking1989general} will be used
\begin{subequations}
\label{eq:set identities}
\begin{align}
&\!\!\! \mathcal{A} \cup (\mathcal{B} \cap \mathcal{C})=(\mathcal{A} \cup \mathcal{B}) \cap (\mathcal{A} \cup \mathcal{C})\\
&\!\!\! \mathcal{A} \cap (\mathcal{B} \cup \mathcal{C})=(\mathcal{A} \cap \mathcal{B}) \cup (\mathcal{A} \cap \mathcal{C})\\
    &\!\!\! (\mathcal{A}\cap\mathcal{C})\backslash(\mathcal{B}\cap\mathcal{C})=(\mathcal{A}\cap\mathcal{C})\backslash\mathcal{B}=\mathcal{A}\cap(\mathcal{C}\backslash\mathcal{B})\label{eq:set identities:c}\\
    &\!\!\! \mathcal{C}\backslash(\mathcal{A}\!\cup\!\mathcal{B})\!=\!(\mathcal{C}\backslash\mathcal{A})\!\cap\!(\mathcal{C}\backslash\mathcal{B}), \mathcal{C}\backslash(\mathcal{A}\!\cap\!\mathcal{B})\!=\!(\mathcal{C}\backslash\mathcal{A})\!\cup\!(\mathcal{C}\backslash\mathcal{B})\\
    &\!\!\! (\mathcal{A}\cup\mathcal{B})^\C\!=\!\mathcal{A}^\C\!\cap\!\mathcal{B}^\C, (\mathcal{A}\cap\mathcal{B})^\C\!=\!\mathcal{A}^\C\!\cup\!\mathcal{B}^\C,\mathcal{A}\backslash\mathcal{B}\!=\!\mathcal{A}\!\cap\!\mathcal{B}^\C\\
    &\!\!\! \overline{\mathcal{A}\cup\mathcal{B}}=\overline{\mathcal{A}}\cup\overline{\mathcal{B}},\, \overline{\mathcal{A}\cap\mathcal{B}}\subset\overline{\mathcal{A}}\cap\overline{\mathcal{B}},\, \mathcal{A}\backslash\overline{\mathcal{A}^\C}=\mathcal{A}^\circ \label{eq:set identities:h}\\
    &\!\!\! \mathcal{A}\subset\mathcal{B}\implies \overline{\mathcal{A}}\subset\overline{\mathcal{B}},\quad  \partial\mathcal{A}^\C=\partial\mathcal{A}\\
    &\!\!\! \partial(\mathcal{A}\cup\mathcal{B})\subset(\partial\mathcal{A}\backslash\mathcal{B}^\circ)\cup(\partial\mathcal{B}\backslash\mathcal{A}^\circ)\\
    &\!\!\! \partial(\mathcal{A}\cap\mathcal{B})\subset(\partial\mathcal{A}\cap\overline{\mathcal{B}})\cup(\partial\mathcal{B}\cap\overline{\mathcal{A}}). \label{eq:set identities:l}
\end{align}    
\end{subequations}

Two sets $\mathcal{A}$ and $\mathcal{B}$ are said to be disjoint if $\mathcal{A}\cap\mathcal{B}=\emptyset$. They are said to be separated if
$\mathcal{A}\cap\overline{\mathcal{B}}=\emptyset=\overline{\mathcal{A}}\cap\mathcal{B}$. The notion of separated sets is stronger than mere disjointness. If two sets $\mathcal{A}$ and $\mathcal{B}$ are separated then we have \cite[Exercise 1.3.A]{engelking1989general}
\begin{align}
\label{eq:boundary of separated sets}
    \partial(\mathcal{A}\cup\mathcal{B})=\partial\mathcal{A}\cup\partial\mathcal{B}.
\end{align}
The tangent cone to a set $\mathcal{K} \subset \real^n$ at a point $x \in \real^n$, denoted $\mathbf{T}_{\mathcal{K}}(x)$, is defined as in~\cite[Def.~5.12 and Fig.~5.4]{goebel2012hybrid}.
\subsection{Projections Maps}
For $z\in\mathbb{R}^n\backslash\{0\}$, we define the following projection maps:
\begin{equation}
\label{eq:proj-refl-maps}
\pi^\parallel(z):=\tfrac{zz^\top}{\|z\|^2},\, \pi^\perp(z):=\!I_n\!-\tfrac{zz^\top}{\|z\|^2},\, \rho(z):=\!I_n\!-2\tfrac{zz^\top}{\|z\|^2}
\end{equation}
where $I_n$ is the $n\times n$ identity matrix. The map $\pi^\parallel(\cdot)$ is the parallel projection map, $\pi^\perp(\cdot)$ is the orthogonal projection map \cite{meyer2000matrix}, and $\rho(\cdot)$ is the reflector map (also called Householder transformation). Consequently, for any $x\in\mathbb{R}^n$, the vector $\pi^\parallel(z)x$ corresponds to the projection of $x$ onto the line generated by $z$, $\pi^\perp(z)x$ corresponds to the projection of $x$ onto the hyperplane orthogonal to $z$ and $\rho(z)x$ corresponds to the reflection of $x$ about the hyperplane orthogonal to $z$. For $z\in\real^n \backslash\{ 0\}$, some useful properties of these maps follow:
\begin{subequations}
\begin{align}
    \label{eq:propLine1}
    \pi^\parallel(z)z&=z,&\pi^\perp(z)\pi^\perp(z)&=\pi^\perp(z),\\
    \label{eq:propLine2}
    \pi^\perp(z)z&=0,&\pi^\parallel(z)\pi^\parallel(z)&=\pi^\parallel(z), \\
    \label{eq:propLine3}
    \rho(z)z&=-z,&\rho(z)\rho(z)&=I_n.
\end{align}    
\end{subequations}
We also define for $z\in\real^n\backslash\{ 0\}$
and $\theta\in\real$ the parametric map 
\begin{align}
\label{eq:def:piTheta}
    \pi^\theta(z):=\cos^2(\theta)\pi^\perp(z)-\sin^2(\theta)\pi^\parallel(z),
\end{align}
which can also be written (thanks to $2 \pi^\perp(z) - \rho(z) = 2 \pi^\parallel(z)+\rho(z)=I_n$) as
\begin{equation}
\label{eq:piTheta:3}
\pi^\theta(z) =\tfrac{1}{2}\rho(z)+\tfrac{1}{2}\cos(2\theta)I_n.
\end{equation}

\subsection{Geometric Subsets of \texorpdfstring{$\mathbb{R}^n$}{Rn}}
\subsubsection{Line}
A line is the one-dimensional subset of $\mathbb{R}^n$ described by the set 
\begin{align}
        \label{eq:def:line}
     \mathcal{L}(c,v)&:=\{x\in\mathbb{R}^n: x=c+\lambda v, \lambda\in\mathbb{R}\},
\end{align}
which corresponds to the line passing by the point $c\in\mathbb{R}^n$ and with direction parallel to $v\in\mathbb{R}^{n}\backslash\{0\}$. If in \eqref{eq:def:line} $\lambda\ge 0$ (respectively $\lambda\le 0$), then we obtain the half-line denoted by $\mathcal{L}_\geq(c,v)$ (respectively $\mathcal{L}_\leq(c,v)$).
\subsubsection{Hyperplane}
A hyperplane is the $(n-1)$-dimensional subset of $\mathbb{R}^n$ described by the set
\begin{align}
    \label{eq:def:plane}
\mathcal{P}(c,v)&:=\{x\in\mathbb{R}^n: v^\top(x-c)=0\},
\end{align}
which corresponds to the hyperplane that passes through a point $c\in\mathbb{R}^n$ and has normal vector $v\in\mathbb{R}^{n}\backslash\{0\}$. 
The hyperplane $\mathcal{P}(c,v)$ divides the Euclidean space $\mathbb{R}^n$ into two closed subsets $\mathcal{P}_{\geq}(c,v)$ and $\mathcal{P}_\leq(c,v)$, which are obtained by substituting the $=$ in \eqref{eq:def:plane} with $\geq$ and $\leq$, respectively.
\subsubsection{Sphere}
A sphere is the $(n-1)$-dimensional subset of $\mathbb{R}^n$ described by the set
\begin{align}
\label{eq:def:sphere}
    \mathcal{S}(c,r)&:=\{x\in\mathbb{R}^n: \|x-c\|=r\}
\end{align}
where $c$ is the center of the sphere and $r\in\mathbb{R}_{\geq}$ is its radius. The closed interior (respectively exterior) of the sphere, also called a hyperball and denoted by $\mathcal{S}_\leq(c,r)$ (respectively $\mathcal{S}_\geq(c,r)$), is obtained from \eqref{eq:def:sphere} by substituting the $=$ with $\leq$ (respectively $\geq$).
\subsubsection{Ellipsoid}
For a positive definite matrix $E\in\mathbb{R}^{n\times n}$, a ellipsoid is the $(n-1)$-dimensional subset of $\mathbb{R}^n$ described by the set
\begin{align}
\label{eq:def:ellipsoid}
    \mathcal{E}(c,E)&:=\{x\in\mathbb{R}^n: \|E(x-c)\|^2=1\}
\end{align}
where $c$ is the center of the ellipsoid and its $i$-th principal semi-axis is the vector $\lambda_i^{-1}(E) v_i$, with $v_i$ the unit eigenvector corresponding to the eigenvalue $\lambda_i(E)$. The closed interior (respectively exterior) of the ellipsoid, denoted by $\mathcal{E}_\leq(c,E)$ (respectively $\mathcal{E}_\geq(c,E)$), is obtained from \eqref{eq:def:ellipsoid} by substituting the $=$ with $\leq$ (respectively $\geq$). 
\begin{definition}\label{definition:weak:disjoint}
Two ellipsoids $\mathcal{E}_\leq(c_1,E_1)$ and $\mathcal{E}_\leq(c_2,E_2)$ are weakly disjoint if $\mathcal{E}_\leq(c_1,E_1)\cap\mathcal{E}_\leq(c_2,E_2)=\emptyset$.
\end{definition}
Explicit algebraic conditions to test weak disjointness of two ellipsoids can be found in \cite[Thm.~6]{choi2006continuous} for $n=2$ and in \cite[Thm.~8]{wang2001algebraic} for  $n=3$.
\begin{definition}
Two ellipsoids $\mathcal{E}_\le(c_1,E_1)$ and $\mathcal{E}_\le(c_2,E_2)$ are strongly disjoint if $(\lambda_{\min}(E_1))^{-1}+(\lambda_{\min}(E_2))^{-1}<\|c_1-c_2\|$.
\end{definition}
Strong disjointness means that the two smallest spherical balls containing the ellipsoids are disjoint. Strong disjointness is more conservative than weak disjointness. 
\subsubsection{Cone}
For a positive definite matrix $E\in\mathbb{R}^{n\times n}$, a cone is the $(n-1)$-dimensional subset of $\mathbb{R}^n$ described by the set
\begin{equation}
\label{eq:def:cone}
\mathcal{C}(c,\!v,\!\theta,\!E)\!:=\!\{x\in\mathbb{R}^n\!\colon\!\cos(\theta)\|Ev\|\|E(x-c)\|\!=\!v^\top\! E^2(x-c)\}
\end{equation}
where $c\in\mathbb{R}^n$ is its vertex, $v\in\mathbb{R}^n\backslash\{0\}$ is its axis and $2\theta\in[0,\pi]$ is its aperture. The cone defined here is sometimes referred to as nappe or half-cone, as opposed to the double cone. The closed interior (respectively, exterior) of the cone, denoted by $\mathcal{C}_\leq(c,v,\theta,E)$ (respectively $\mathcal{C}_\geq(c,v,\theta,E)$), is obtained from \eqref{eq:def:cone} by substituting the $=$ with $\leq$ (respectively $\geq$). A normal vector to the cone surface $\mathcal{C}(c,v,\theta,E)$ at $x$ is
\begin{equation}
\label{eq:normal vector to cone surface}
n(x) := E \pi^\theta(Ev) E(x-c),
\end{equation}
and can be obtained after squaring in~\eqref{eq:def:cone} and taking the gradient. The next fact will be used. 
\begin{lemma}\label{lemma:cones}
Let $v_1,v_2\in\mathbb{S}^{n-1}$ such that $v_1^\top v_2=\cos \theta$ for some $\theta\in(0,\pi]$. Let $\psi_1,\psi_2\in[0,\pi/2]$ with $\psi_1+\psi_2<\theta$. Then for each $c\in\mathbb{R}^n$ and $E\in\mathbb{R}^{n\times n}$ positive definite,
\begin{align*}
    \mathcal{C}_\leq(c,E^{-1}v_1,\psi_1,E)\cap\mathcal{C}_\leq(c,E^{-1}v_2,\psi_2,E)=\{c\}.
\end{align*}
\end{lemma}
\subsection{Hybrid Systems Framework}
We consider hybrid dynamical systems of the class \cite{goebel2012hybrid}, described through constrained differential and difference inclusions for state $X \in \real^n$:
\begin{equation}
\label{Hybrid:general}
\begin{cases}
\dot X\in\mathbf{F}(X), &X\in\mathcal{F},\\
X^+\in \mathbf{J}(X),  & X\in\mathcal{J},
\end{cases}
\end{equation}
where the \textit{flow map} $\mathbf{F}:\mathbb{R}^n\rightrightarrows\mathbb{R}^n$ governs the continuous evolution, the \textit{flow set} $\mathcal{F}\subseteq\mathbb{R}^n$ dictates where continuous evolution can occur. The \textit{jump map} $\mathbf{J}:\mathbb{R}^n\rightrightarrows\mathbb{R}^n$ governs the discrete evolution, and the \textit{jump set} $\mathcal{J}\subseteq\mathbb{R}^n$ defines where discrete evolution can occur. The hybrid system \eqref{Hybrid:general} is defined by its data and denoted $\mathscr{H}=(\mathbf{\mathcal{F}},\mathbf{F},\mathcal{J},\mathbf{J})$.

A subset $\mathbb{T}\subset\mathbb{R}_{\geq}\times\mathbb{N}$ is a \textit{hybrid time domain} if it is a union of a finite or infinite sequence of intervals $[t_j ,t_{j+1}]\times\{j\}$, with the last interval (if existent) possibly of the form $[t_j,T)$ with $T$ finite or $T= + \infty$. The ordering of points on each hybrid time domain is such that $(t,j)\preceq(t^\prime,j^\prime)$ if $t < t^\prime$, or $t = t^\prime$ and $j\leq j^\prime$. A \textit{hybrid solution} is defined in \cite[Def.~2.6]{goebel2012hybrid}. A hybrid solution $\phi$ is maximal if it cannot be extended and complete if its domain $\dom\phi$ (which is a hybrid time domain) is unbounded.
\section{Problem Formulation}\label{section:problem}
We consider a point mass vehicle moving in the $n$-dimensional Euclidean space containing $I\in\mathbb{N}$ obstacles denoted by $\mathcal{O}_1,\cdots,\mathcal{O}_I$. For each 
\begin{equation*}
i\in\{1,\cdots,I\}=: \mathbb{I},
\end{equation*}
the obstacle $\mathcal{O}_i$ has an ellipsoidal shape such that $\mathcal{O}_i:=\mathcal{E}_\leq(c_i,E_i)$, for some center $c_i\in\mathbb{R}^n$ and some positive definite matrix $E_i\in\mathbb{R}^{n\times n}$ defining the orientation and the shape of the obstacle. The free workspace (obstacle-free region) is then defined by the closed set
\begin{align}\label{eq:W}
  \mathcal{W}:=\bigcap_{i^\prime \in \mathbb{I}}\mathcal{E}_\geq(c_{i^\prime},E_{i^\prime}).
\end{align}
The vehicle is moving according to the dynamics
\begin{equation}
\label{eq:integrator}
\dot x=u,
\end{equation}
where $x\in\mathbb{R}^n$ is the state and $u\in\mathbb{R}^n$ is the control input. The vehicle is required to stabilize its position to a target position while avoiding the obstacles. Without loss of generality we consider the target position to be $x=0$ (the origin).
\begin{assumption}\label{assumption:n}
$n\geq 2$.
\end{assumption}
We consider $n\geq 2$ since for $n=1$ (i.e., the state space is a line), global asymptotic stabilization with obstacle avoidance is infeasible.
\begin{assumption}\label{assumption:obstacle}
For all $i\in\mathbb{I}$, $\|E_ic_i\|>1$.
\end{assumption}
Assumption \ref{assumption:obstacle} requires that the target position $x=0$ is not inside any of the obstacle regions $\mathcal{O}_i$, otherwise the considered navigation problem would be infeasible.
\begin{assumption}\label{assumption:obstacle2}
$\{\mathcal{O}_i\}_{i\in\mathbb{I}}$ are weakly pairwise disjoint. 
\end{assumption}
In Assumption \ref{assumption:obstacle2} we impose that there is no intersection region between the obstacles. Otherwise, the union of the two intersecting obstacles forms another region which might have a different shape than an ellipsoid. Our \emph{objectives} in designing a control strategy are:
\begin{itemize}
    \item[i)] the obstacle-free region $\mathcal{W}$ in~\eqref{eq:W} is forward invariant,
    \item[ii)] the target $x=0$ is globally asymptotically stable.
    \end{itemize}
Objective~i) guarantees that all solutions of the closed-loop system are safely avoiding the obstacles by remaining in the obstacle-free region $\mathcal{W}$  for all times while objective~ii) corresponds to global stabilization of the target. 
\section{Hybrid Control for Obstacle Avoidance}\label{section:controller}
In this section, we propose a hybrid controller that switches suitably between a {\it stabilizing} and an {\it avoidance} controller. 
Let us define a discrete variable 
\begin{equation*}
m\in \{-1,0,1\}=:\mathbb{M}.
\end{equation*}
The value $m=0$ corresponds to the activation of the stabilizing controller and the values $m=-1$, $m=1$ correspond to the activation of one of the two configurations of the avoidance controller. The avoidance controller depends also on the current obstacle $\mathcal{O}_i$, as  detailed in the next sections.
\vskip -0.5cm
\subsection{Control Input}
\begin{figure}
    \centering
    \includegraphics[scale=.8]{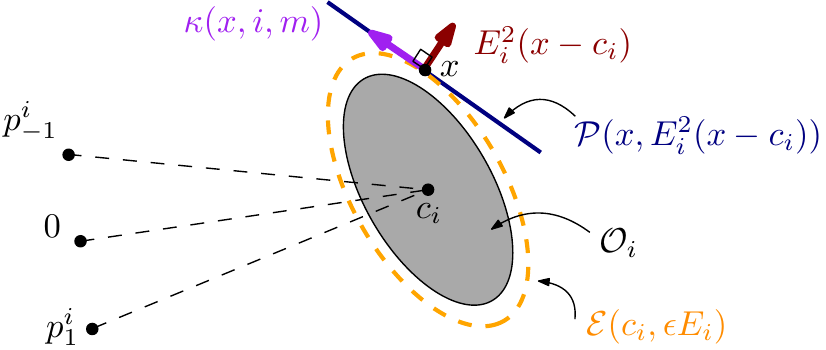}
    \caption{Illustration of the projection-based avoidance controller. The vehicle is attracted to the auxiliary point $p_m^i$ while sliding on a neighbouring ellipsoid.}
    \label{fig:projection}
\end{figure}

In this section we propose the feedback law for the control input $u$ in~\eqref{eq:integrator}. $u$ depends on the state $x\in\mathbb{R}^n$, the obstacle $i\in\mathbb{I}$ and the control mode $m\in \mathbb{M}$ as 
\begin{align}
    &u =\kappa(x,i,m) \label{eq:u}\\
    &:=\begin{cases}
    -k_0 x, & m=0,\\
    - k_mE_i^{-1}\pi^\perp(E_i(x-c_i))E_i(x-p_m^i),& m \in \{-1,1\},
    \end{cases} \nonumber
\end{align}
where $k_{-1},\,k_0,\,k_1>0$ are the control gains for each control mode $m\in \mathbb{M}$ and the points $p_m^i\in\mathbb{R}^n$, $m\in\{-1,1\}$ and $i\in\mathbb{I}$, are design parameters defined below. In the stabilization mode ($m=0$), the control input in~\eqref{eq:u} steers $x$ towards the origin under a state feedback. In the avoidance mode depicted in Fig.~\ref{fig:projection}, the control input minimizes the distance to the {\it auxiliary} attractive point $p_m^i$ {\it while} maintaining a constant distance to the obstacle $\mathcal{O}_i$. 
Indeed, the time derivative of $\|E_i(x-c_i)\|^2$ along solutions of $\dot x=\kappa(x,i,m)$ for $m\in\{-1,1\}$ and $i\in\mathbb{I}$, reads
\begin{equation}\label{eq:avoidance}
\begin{aligned}
    &\tfrac{1}{2}\tfrac{d}{dt}\|E_i(x-c_i)\|^2=(x-c_i)^\top E_i^2\kappa(x,i,m)\\
    &=- k_m(x-c_i)^\top E_i\pi^\perp(E_i(x-c_i))E_i(x-p_m^i)=0
\end{aligned}
\end{equation}
by~\eqref{eq:propLine2}. Then, if we activate the avoidance mode sufficiently away from the obstacle, the avoidance feedback $u=\kappa(x,i,m)$ guarantees that the vehicle does not hit the obstacle. Whereas the logic variable $i$ corresponds to obstacle $\mathcal{O}_i$, the logic variable $m$ is selected according to a hybrid mechanism that exploits a suitable construction of the flow and jump sets as detailed in Section~\ref{section:construction}. 

In order to clear the obstacle while approaching the desired target position at the origin, we select the points $p_1^i$ and $p_{-1}^i$ in the region between the obstacle and the origin, see Fig.~\ref{fig:projection}. The motivation is that the avoidance task is equivalent (up to a linear transformation) to a stabilization problem on the unit sphere $\mathbb{S}^{n-1}$. Therefore, as pointed out for instance in \cite{mayhew2013global}, global asymptotic stabilization cannot be accomplished by only one continuous time-invariant controller, but it can be by a hybrid feedback with at least two configurations. For this reason, we consider two avoidance modes with $m=-1$ and $m=1$ and, hence, the points $p_1^i$ and $p_{-1}^i$ must be distinct. More precisely, for $\theta_i>0$ (which will be further bounded in Lemma \ref{lemma:empty1}), the points $p_1^i$ and $p_{-1}^i$ are selected as
\begin{subequations}
\label{eq:p-1}
\begin{align}
    \label{eq:p_1}p_1^i& \in\mathcal{C}(c_i,-c_i,\theta_i,E_i)\backslash\{c_i\},\\
    \label{eq:p_-1}p_{-1}^i&:=-E_i^{-1}\rho(E_i c_i)E_i p_1^i.
\end{align}
\end{subequations}
By~\eqref{eq:p-1}, $p_{-1}^i$ opposes $p_1^i$ diametrically with respect to the cone axis (for $E_i=I_n$, $p_{-1}^i$ is obtained by an orthogonal reflection) and also belongs to $\mathcal{C}(c_i,-c_i,\theta_i,E_i)\backslash\{c_i\}$ as shown in the next lemma.
\begin{lemma}\label{lemma:p-1}
$p_{-1}^i\in\mathcal{C}(c_i,-c_i,\theta_i,E_i)\backslash\{c_i\}.$
\end{lemma}
Note that the results of the paper hold for any selection of the point $p_1^i$ as long as it lies on the surface of the cone as in \eqref{eq:p_1}. An explicit guided choice for those points is given in Section~\ref{section:example} for the 2D and 3D cases. Finally, further motivation about the choice of the avoidance controller mode in \eqref{eq:u} is detailed in Section~\ref{section:construction} and, in particular, in Lemma~\ref{lemma:equilibria}, which is important for the construction of flow and jump sets.

\subsection{Geometric Construction of the Flow and Jump sets}\label{section:construction}
In this section we construct explicitly the flow and jump sets where the stabilization and avoidance controllers are activated. 
\subsubsection{Safety Helmets} 
\begin{figure}
    \centering
    \includegraphics[width=.75\columnwidth]{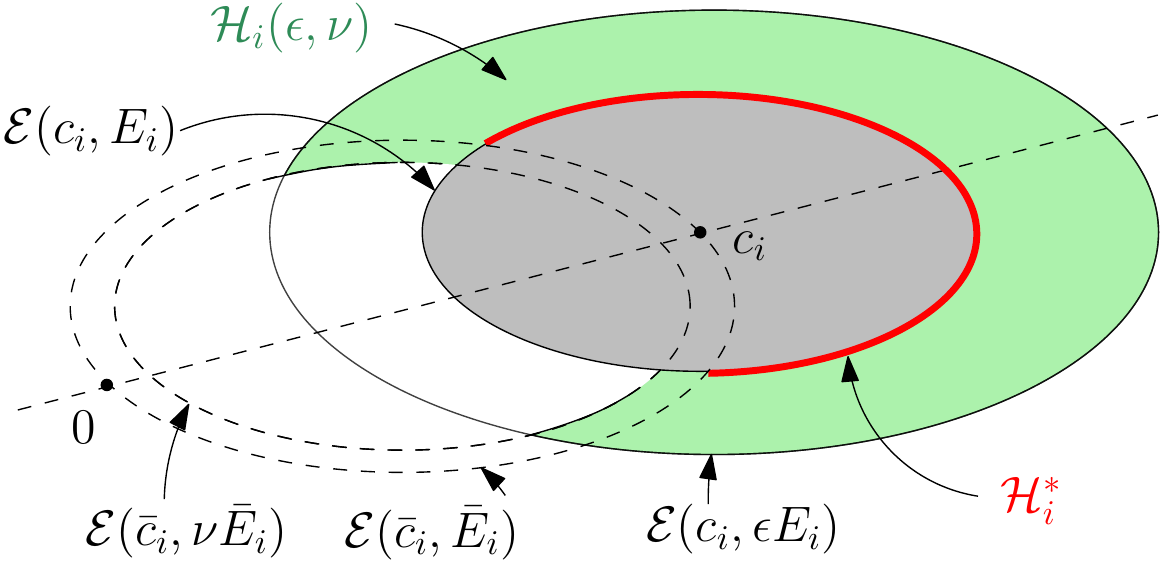}
    \caption{The helmet $\mathcal{H}_i^*$ in \eqref{eq:helmet} (red) corresponds to all boundary points where the stabilization vector field is pointing inside the obstacle (grey). The safety helmet $\mathcal{H}_i(\epsilon,\nu)$ in \eqref{eq:safety:helmet} (green) corresponds to an dilated version of $\mathcal{H}_i^*$.}
    \label{fig:helmet}
\end{figure}

Our proposed construction of flow and jump sets is based on regions that have the shape of a {\it helmet}, whose construction is now motivated. In the stabilization mode $m=0$, the closed-loop system should \emph{not} flow when: 1) $x$ is close enough to any of the obstacle regions  $\mathcal{E}_\leq(c_i,E_i)$ and 2) the vector field $-k_0x$ points inside $\mathcal{E}_\leq(c_i,E_i)$. Otherwise, the vehicle ends up hitting the obstacle $i$. Indeed, by computing the time derivative of $\|E_i(x-c_i)\|^2$ along solutions of the vector field $-k_0x$, we obtain
\begin{equation}
    \begin{aligned}
        &\tfrac{1}{2}\tfrac{d}{dt}\|E_i(x-c_i)\|^2 =-k_0 x^\top E_i^2(x-c_i)\\
        &=k_0c_i^\top E_i^2c_i/4-k_0(x-c_i/2)^\top E_i^2(x-c_i/2)\\
        &=k_0\|E_i \bar c_i\|^2\left(1-\|\bar E_i(x-\bar c_i)\|^2\right)
        \label{eq:dist_to_c_decreases}
  \end{aligned}
\end{equation}
where $\bar c_i$ and $\bar E_i$ are defined as
\begin{align}
\label{eq:ci,Ei}
    \bar c_i:={c_i}/{2},\quad\bar E_i:={2E_i}/(\|E_i c_i\|).
\end{align}
\eqref{eq:dist_to_c_decreases} implies that the distance function $\|E_i(x-c_i)\|^2$ decreases for all $x$ in the closed set $\mathcal{E}_\geq(\bar c_i,\bar E_i)$. Consider now Fig.~\ref{fig:helmet} for a sketch of the next sets and for obstacle $i$, define the {\it helmet}-shaped set 
\begin{align}\label{eq:helmet}
    \mathcal{H}^*_i:=\mathcal{E}(c_i,E_i)\cap\mathcal{E}_\geq(\bar c_i,\bar E_i).
\end{align}
$\mathcal{H}^*_i$ is the set of all points that lie on the boundary of the obstacle $\mathcal{O}_i$ and generate a vector field pointing towards the obstacle. Then, for obstacle $i$, we define the {\it safety helmet} as:
\begin{align}\label{eq:safety:helmet}
    \mathcal{H}_i(\epsilon,\nu):=\mathcal{E}_\leq(c_i,\epsilon E_i)\cap\mathcal{E}_\geq(c_i,E_i)\cap\mathcal{E}_\geq(\bar c_i,\nu\bar E_i)
\end{align}
for some parameters $\epsilon,\,\nu>0$. $\epsilon$ and $\nu$ determine the thickness of the safety helmet by tuning the dilation/shrinking of the ellipsoids $\mathcal{E}(c_i,E_i)$ and $\mathcal{E}(\bar c_i,\bar E_i)$, thereby generating a dilated version of $\mathcal{H}^*_i$. The safety helmet $\mathcal{H}_i(\epsilon,\nu)$ constitutes the main ingredient of our following constructions. 

\subsubsection{Stabilization Mode \texorpdfstring{$m=0$}{m=0}}
\begin{figure}
    \centering
    \includegraphics[width=.9\columnwidth]{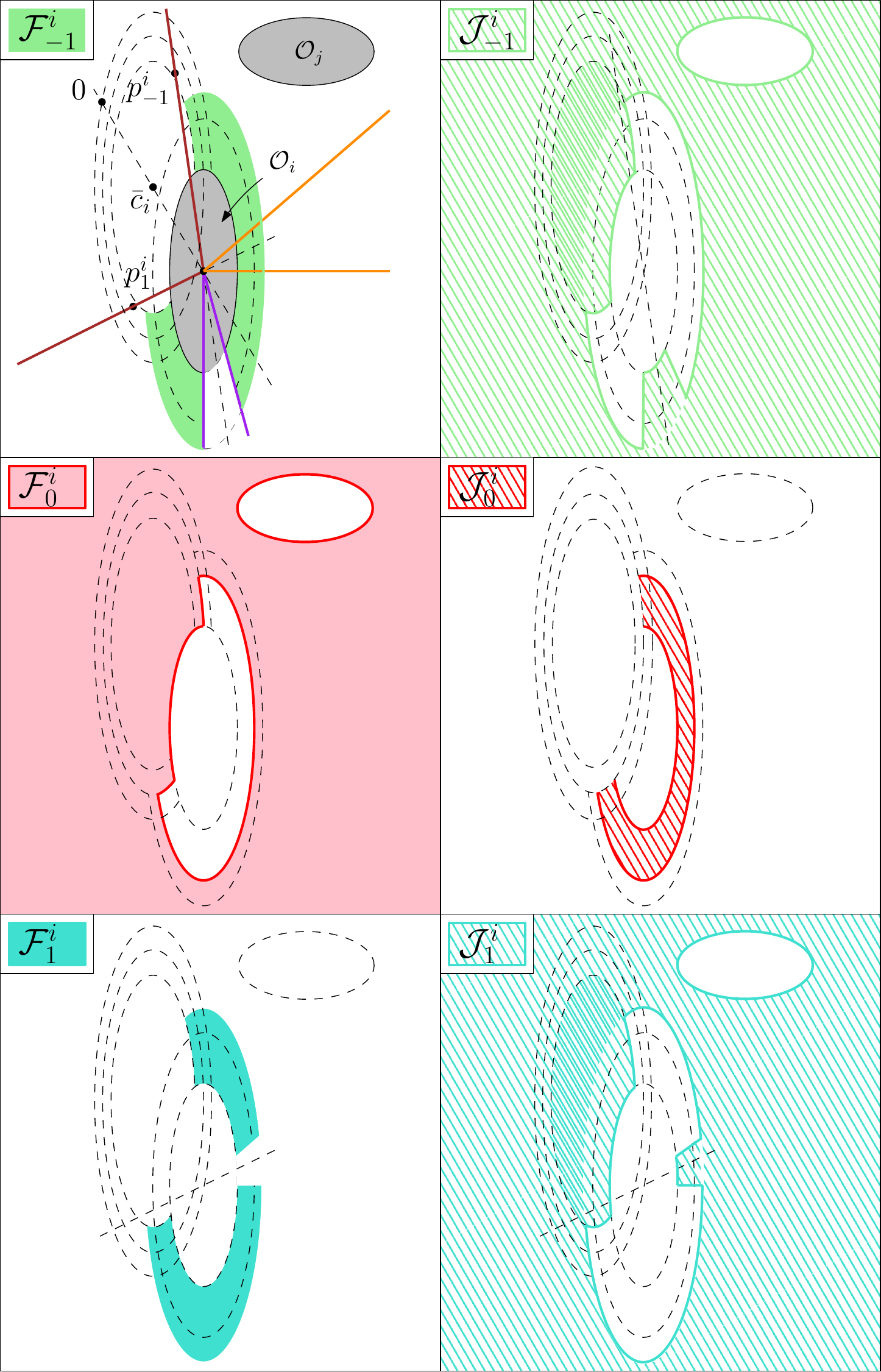}
    \caption{2D illustration of flow and jump sets considered in Sections~\ref{section:controller}-\ref{section:main} corresponding to obstacle $\mathcal{O}_i$ (in the presence of a second obstacle $\mathcal{O}_j$). The stabilization-mode jump set $\mathcal{J}_0^i$ (hatched red) is constructed  by using the helmet $\mathcal{H}_i(\epsilon_i,\nu_i)$, while the corresponding flow set $\mathcal{F}_0^i$ is the complement of $\mathcal{J}_0^i$ in the free workspace. For the avoidance mode we select $p_1^i$ and $p_{-1}^i$ to lie on the cone $\mathcal{C}(c_i,-c_i,\theta_i,E_i)$ (solid brown line). The avoidance flow set $\mathcal{F}_m^i$, with $m\in\{-1,1\}$, corresponds to the helmet $\mathcal{H}_i(\delta_i,\mu_i)$ deprived of the interior of the the cone region defined by $\mathcal{C}(c_i,c_i-p_m^i,\psi_i,E_i)$ (solid purple line for $m=-1$ and solid orange line for $m=1$). The corresponding jump set $\mathcal{J}_m^i$ is the complement of $\mathcal{F}_m^i$ in the free workspace.} 
    \label{fig:flowAndJumpSets}
\end{figure}

Consider from now on Fig.~\ref{fig:flowAndJumpSets} for a visualization of the sets we are introducing in our construction. In stabilization mode ($m=0$), we create around each obstacle $\mathcal{O}_i$ a safety helmet $\mathcal{H}_i(\epsilon_i,\nu_i)$ which adds a safety layer to the given obstacle. The controller mode must be switched to the avoidance mode whenever the vehicle reaches this safety helmet. Specifically, we define for each $i\in\mathbb{I}$, a jump set 
\begin{align}
\label{eq:J0i}
&\mathcal{J}_0^i:=\mathcal{H}_i(\epsilon_i,\nu_i)\cap\mathcal{W},
\end{align}
where $\epsilon_i\in(0,1)$ (dilating $\mathcal{E}_\leq(c_i, E_i)$ to $\mathcal{E}_\leq(c_i,\epsilon_i E_i)$), and $\nu_i\in(1,\infty)$ (shrinking $\mathcal{E}_\geq(\bar c_i,\bar E_i)$ to $\mathcal{E}_\geq(\bar c_i,\nu_i\bar E_i)$) and $\mathcal{W}$ is the free workspace defined in~\eqref{eq:W}. 
We emphasize that we consider the intersection with $\mathcal{W}$ in~\eqref{eq:J0i} for convenience, but later we tune the parameters such that $\mathcal{H}_i(\epsilon_i,\nu_i)\subset\mathcal{W}$, which implies $\mathcal{J}_0^i$ will equal to $\mathcal{H}_i(\epsilon_i,\nu_i)$. The selection of $\mathcal{J}_0^i$ in~\eqref{eq:J0i} leads naturally to the following flow set of the stabilization mode (corresponding to the closed complement of $\mathcal{J}_0^i$ in the free workspace)
\begin{align}
\label{eq:F0i}
&\mathcal{F}_0^i:=\Big(\mathcal{E}_\geq(c_i,\epsilon_i E_i)\cup\mathcal{E}_\leq(\bar c_i,\nu_i\bar E_i)\Big)\cap\mathcal{W}.
\end{align}
Finally, from~\eqref{eq:J0i} and \eqref{eq:F0i}, we take all the obstacles into account and define the flow and jump sets for the stabilization mode $m=0$ as 
\begin{equation}
\label{eq:FJ0}
\mathcal{F}_0:=\Big(\bigcap\limits_{i\in\mathbb{I}}\mathcal{F}_0^i\Big)\times \mathbb{I},\qquad
\mathcal{J}_0:=\Big(\bigcup\limits_{i\in\mathbb{I}}\mathcal{J}_0^i\Big) \times \mathbb{I}.
\end{equation}
Indeed, the stabilization mode can be selected when the state $x$ belongs to the intersection of the sets $\mathcal{F}_0^i$ (and for \emph{any} obstacle index $i\in\mathbb{I}$), and a jump to the avoidance mode can occur when the state $x$ belongs to the union of the sets $\mathcal{J}^i_0$ (and for \emph{any} obstacle index $i\in\mathbb{I}$). In other words, if during the stabilization mode the vehicle reaches any one of the safety helmets, then the controller jumps to one of the avoidance modes with $m$ equal to $-1$ or $1$.
\subsubsection{Avoidance Mode \texorpdfstring{$m\in\{-1,1\}$}{m=-1,1}}
We consider now the construction of flow and jump sets for the avoidance modes $m \in \{-1,1\}$ and the specific obstacle $i \in \mathbb{I}$ with the aid of Fig.~\ref{fig:flowAndJumpSets}. To highlight their motivation, we first define such flow sets and state later the corresponding jump sets (see~\eqref{eq:J^i_m}). For each $i \in \mathbb{I}$ and $m \in\{-1,1\}$, the avoidance flow set is
\begin{equation}
\label{eq:F-11i}
\begin{aligned}
\mathcal{F}_m^i:=\mathcal{H}_i(\delta_i,\mu_i)\cap\mathcal{C}_\geq(c_i,c_i-p_m^i,\psi_i,E_i)\cap\mathcal{W},
\end{aligned}    
\end{equation}
where $\delta_i\in(0,\epsilon_i)$ (dilating $\mathcal{E}_\leq(c_i, \epsilon_iE_i)$ to $\mathcal{E}_\leq(c_i,\delta_i E_i)$), $\mu_i\in(\nu_i,\infty)$ (shrinking $\mathcal{E}_\geq(\bar c_i,\nu_i \bar E_i)$ to $\mathcal{E}_\geq(\bar c_i,\mu_i\bar E_i)$), and $\psi_i\in(0,\pi/2]$. In the two configurations $m\in\{-1,1\}$ of the avoidance of obstacle $i \in \mathbb{I}$, we want the vehicle to slide on the safety helmet $\mathcal{H}_i(\delta_i,\mu_i)$ while maintaining a constant distance to the obstacle. 
By selecting $\delta_i\in(0,\epsilon_i)$ and $\mu_i\in(\nu_i,\infty)$, one obtains a dilated version of $\mathcal{H}_i(\epsilon_i,\nu_i)$ used in $\mathcal{J}_0^{i}$ and, thus, creates a hysteresis region useful to prevent infinitely many consecutive jumps (Zeno behavior). 
However, the avoidance vector field $\kappa(x,i,m)$ in \eqref{eq:u} has some undesirable equilibria, which we need to rule out from the flow sets $\mathcal{F}_1^i$ and $\mathcal{F}_{-1}^i$. These are characterized in the next lemma.
\begin{lemma}\label{lemma:equilibria}
Let $c\in \real^n$, $p\in \real^n\backslash\{c\}$ and $E\in \real^{n\times n}$ positive definite. 
For each $x\in\real^n \backslash \{c\}$, $\pi^\perp(E(x-c))E(x-p)=0$ if and only if $x\in\mathcal{L}(c,p-c)$.
\end{lemma}
For each $m\in\{-1, 1\}$, $i \in \mathbb{I}$, we want solutions to eventually leave the set $\mathcal{F}_m^i$ of the avoidance mode, so it is necessary to select point $p_m^i$ and flow set $\mathcal{F}_m^i$ such that $\mathcal{L}(c_i,p_m^i-c_i)\cap\mathcal{F}_m^i = \emptyset$ based on Lemma~\ref{lemma:equilibria}, otherwise solutions could stay in avoidance mode indefinitely. This motivates the intersection with the cone in~\eqref{eq:F-11i}, and the next lemma.
\begin{lemma}\label{lemma:empty1}
For each $i \in \mathbb{I}$, define the quantities
\begin{subequations}
\label{eq:bar delta,bar mu,bar theta}
\begin{align}
   &\underline{\delta}_i :=\|E_i c_i\|^{-\frac{1}{2}} \label{eq:bar delta}\\
   &\bar\mu_i(\delta_i):=\big(1-4\underline{\delta}_i^2(1-\underline{\delta}_i^2/\delta_i^2)\big)^{-\frac{1}{2}} \label{eq:bar mu}\\
  & \bar\theta_i(\delta_i,\mu_i) :=\arccos\left(\frac{\underline{\delta}_i^2}{\delta_i^2}+\frac{1}{4\underline{\delta}_i^2}\left(1-\frac{1}{\mu_i^2}\right)\right) \label{eq:bar theta}
\end{align}
\end{subequations}
and select the parameters $\delta_i, \mu_i,\theta_i, \psi_i$ as in Table \ref{table:parameters} so that $\bar\mu_i(\delta_i)$ and $\bar\theta_i(\delta_i,\mu_i)$ are well-defined. Then, for each $m\in\{-1,1\}$, $\mathcal{L}(c_i,p_m^i-c_i)\cap\mathcal{F}_m^i = \emptyset$.
\end{lemma}

From the flow set in~\eqref{eq:F-11i}, we suitably define the jump set for the avoidance mode, of an obstacle $i \in \mathbb{I}$ with configuration $m \in \{-1,1\}$, to be the closed complement of $\mathcal{F}_m^i$ in the free workspace. For $i \in \mathbb{I}$ and $m \in\{-1,1\}$,
\begin{align}
\label{eq:J^i_m}
\mathcal{J}_m^i:=\Big(\mathcal{E}_\geq(c_i,\delta_i E_i)&\cup\mathcal{E}_\leq(\bar c_i,\mu_i\bar E_i)\\
\nonumber
&\cup\mathcal{C}_\leq(c_i,c_i\!-\!p_m^i,\psi_i,E_i)\Big)\cap\mathcal{W}.
\end{align}
Finally, the avoidance mode has overall flow and jump sets
\begin{subequations}
\begin{align}
\label{eq:FJ1}
\mathcal{F}_1\!&:=\bigcup\limits_{i\in\mathbb{I}}\left(\mathcal{F}_1^i\!\times\! \{i\}\right),
& 
\mathcal{J}_1\! &:=\bigcup\limits_{i\in\mathbb{I}}\left(\mathcal{J}_1^i\!\times\! \{i\}\right),\\
\label{eq:FJ-1}
\mathcal{F}_{-1}\!&:=\bigcup\limits_{i\in\mathbb{I}}\left(\mathcal{F}_{-1}^i\!\times\! \{i\}\right),
&
\mathcal{J}_{-1}\!&:=\bigcup\limits_{i\in\mathbb{I}}\left(\mathcal{J}_{-1}^i\!\times\! \{i\}\right),
\end{align}
\end{subequations}
where $\mathcal{F}_m^i$ and $\mathcal{J}_m^i$ ($m\in\{-1,1\}$) are defined in~\eqref{eq:F-11i} and \eqref{eq:J^i_m}. Indeed, \emph{each} obstacle $i$ gives rise, for the avoidance mode, to a specific flow (jump) set with two configurations $\mathcal{F}_1^i$ and $\mathcal{F}_{-1}^i$ ($\mathcal{J}_1^i$ and $\mathcal{J}_{-1}^i$), as we motivated in this section. 
\begin{table}[t]
\caption{Selection of the design parameters of~\eqref{eq:hs}, with $i \in \mathbb{I}$.}
\centering
\begin{tabularx}{0.88\columnwidth}{XX|XX}
\toprule
Parameter     				&  Selection						&Parameter				&  Selection\\
$\delta_i$                	&$(\underline{\delta}_i,1)$			&$\epsilon_i$			&$(\delta_i,1)$\\
$\mu_i$                   	&$(1,\bar\mu_i(\delta_i))$			&$\nu_i$				&$(1,\mu_i)$\\
$\theta_i $               	&$(0,\bar\theta_i(\delta_i,\mu_i))$	&$\bar \psi_i$			&$(0,\theta_i)$\\
$k_0$, $k_1$, $k_{-1}$      &$(0,+\infty)$	 					&$\psi_i$				&$(0,\bar\psi_i)$\\
  							& 									&$p^i_1$, $p^i_{-1}$	& as in \eqref{eq:p-1}\\
\bottomrule
\end{tabularx}
\label{table:parameters}
\end{table}

\subsection{Hybrid Mode Selection}
\label{section:selection}
In this section we define the hybrid switching strategy that permits a Zeno-free transition between the different control modes. The hybrid selection of the logical variables $i\in \mathbb{I}$ and $m\in \mathbb{M}$ is implemented in the hybrid system
\begin{subequations}
\label{eq:hs}
\begin{align}
&\left\{
\begin{aligned}
\dot x&=\kappa(x,i,m)\\
\dot{\aoverbrace[L1R]{i}} & =0\\
\dot m&=0
\end{aligned}
\right.&&(x,i,m)\in \mathcal{F}
\label{eq:flow}\\
&\left\{
\begin{aligned}
 x^+ &=x\\
\smat{
i^+\\
m^+
}
&\in\mathbf{L}(x,i,m)
\end{aligned}
\right.&&(x,i,m)\in \mathcal{J}
\label{eq:jump}
\end{align}
where $\kappa(x,i,m)$ is the control input as defined in~\eqref{eq:u} and the flow and jump sets are given by 
\begin{align}
\label{eq:FJ}
\mathcal{F}& :=\bigcup\limits_{m\in\mathbb{M}}\left(\mathcal{F}_m\!\times\! \{m\}\right), & \mathcal{J}&:=\bigcup\limits_{m\in\mathbb{M}}\left(\mathcal{J}_m\!\times\! \{m\}\right).
\end{align}
with $\mathcal{F}_m$ and $\mathcal{J}_m$ being defined in \eqref{eq:FJ0} for $m=0$ and in \eqref{eq:FJ1}-\eqref{eq:FJ-1} for $m\in\{-1,1\}$. We define now the (set-valued) jump map $\mathbf{L}$ in~\eqref{eq:jump}. To this end, for $i\in \mathbb{I}$ and $m \in \{-1,1\}$, define the sets $\mathcal{C}^i_m$ as
\begin{align}
\label{eq:Cim}
    \mathcal{C}^i_m:=\mathcal{C}_\geq(c_i,c_i-p^i_m,\bar\psi_i,E_i)
\end{align}
which corresponds to the region outside the cone with vertex at $c_i$, axis $c_i-p_m^i$ and aperture $2\bar\psi_i$, where $\bar \psi_i$ is a design parameter selected below. The jump map $\mathbf{L}$ for $m\in\{-1,1\}$ is then defined as
\begin{equation}
\label{eq:L_m1,L_m-1}
\mathbf{L}(x,i,-1) := 
\mathbf{L}(x,i,1) := 
\{
\smat{
i\\0
}
\},
\end{equation}
i.e., when jumping to stabilization mode, the obstacle index $i$ is not used in the control law $\kappa$ in~\eqref{eq:u} and consequently is not updated. The jump map $\mathbf{L}$ for $m=0$ is
\begin{equation}
\label{eq:L_m0}
\mathbf{L}(x,i,0) :=
\left\{
\smat{
i^\prime\\
m^\prime
} \colon
x \in \mathcal{J}^{i^\prime}_0, m^\prime \in \mathbf{M}(x,i^\prime)
\right\}
\end{equation}
where $\mathbf{M}$ is defined, based on~\eqref{eq:Cim}, as
\begin{align}\label{eq:Mi}
    \mathbf{M}(x,i):=\begin{cases}
    \{-1\}&x\in\mathcal{C}^i_{-1}\backslash\mathcal{C}^{i}_1\\
    \{1\}&x\in\mathcal{C}^i_{1}\backslash\mathcal{C}^i_{-1}\\
    \{-1,1\}&x\in\mathcal{C}^i_{-1}\cap\mathcal{C}^{i}_1.\\
    \end{cases}
\end{align}
$\mathbf{L}(\cdot,\cdot,0)$ captures that when jumping from the stabilization mode $m=0$, the suitable avoidance mode of obstacle $i^\prime \in \mathbb{I}$ with configuration $m^\prime \in \{-1,1\}$ is selected based on the position $x$ of the vehicle ($m^\prime$, in particular, is selected based on whether $x$ is within the cone region $\mathcal{C}^{i^\prime}_{-1}$ or $\mathcal{C}^{i^\prime}_1$). A necessary condition to implement our hybrid controller is that the jump map is nonempty, for which we have the next lemma.
\begin{lemma}\label{lemma:jump_map}
Select the parameters $\bar\psi_i$ and $\psi_i$ as in Table \ref{table:parameters}. Then, the set $\mathbf{L}(x,i,m)$ is nonempty for all $(x,i,m)\in\mathcal{J}$.
\end{lemma}
For compact notation, we write flow and jump maps as
\begin{align}
    (x,i,m) & \mapsto \mathbf{F}(x,i,m) := (\kappa(x,i,m),0,0)\\
    (x,i,m) & \mapsto \mathbf{J}(x,i,m) := (x,\mathbf{L}(x,i,m)),
\end{align}
and the overall state of the hybrid system as
\begin{equation}
\label{eq:overall state}
    \xi:=(x,i,m)\in\mathbb{R}^n\times\mathbb{I}\times\mathbb{M}.
\end{equation}
\end{subequations}
This completes the description of the hybrid controller in~\eqref{eq:hs}.
The selections we made in this section for the parameters of~\eqref{eq:hs} are summarized in Table~\ref{table:parameters}.
\section{Main Result}\label{section:main}
In this section, we show that the hybrid controller achieves forward invariance (Section~\ref{section:main:fwdInv}) and global asymptotic stability (Section~\ref{section:main:GAS}) (related to the objectives in Section \ref{section:problem}), as well as some complementary properties (Section~\ref{section:main:complProp}). 

The mild regularity conditions satisfied by the hybrid system~\eqref{eq:hs}, as in the next lemma, allows us to invoke useful results on hybrid systems in the proof of our results.
\begin{lemma}
\label{lemma:hbc}
The hybrid system with data $(\mathcal{F},\mathbf{F},\mathcal{J},\mathbf{J})$ satisfies the hybrid basic conditions in~\cite[Assumption~6.5]{goebel2012hybrid}.
\end{lemma}

\subsection{Forward Invariance}
\label{section:main:fwdInv}

In this section, we show that all generated solutions are complete and safe. Since the state $x$ must evolve always within the free workspace $\mathcal{W}$ in~\eqref{eq:W} regardless of the logic variables $i$ and $m$, we seek forward invariance of the set $\mathcal{K}$ defined as:
\begin{equation}
\label{eq:K}
\mathcal{K}:=\bigcap_{i^\prime \in \mathbb{I}}\mathcal{E}_\geq(c_{i^\prime},E_{i^\prime})\times\mathbb{I}\times\mathbb{M}=\mathcal{W}\times\mathbb{I}\times\mathbb{M}.
\end{equation}
The next lemma shows that the union of flow and jump sets covers exactly the obstacle-free state space $\mathcal{K}$ and that solutions cannot leave $\mathcal{K}$ through jumps.
\begin{lemma}\label{lemma:union}
$\mathcal{F}\cup\mathcal{J}=\mathcal{K}$ and $\mathbf{J}(\mathcal{J})\subset\mathcal{K}$.
\end{lemma}
Forward invariance of $\mathcal{K}$ is proven in the next theorem.
\begin{theorem}\label{theorem:invariance}
Under Assumptions~\ref{assumption:n}-\ref{assumption:obstacle2}, consider the hybrid system \eqref{eq:hs} with parameters selected as in Table~\ref{table:parameters}. Assume also that the controller parameters $\delta_i$ are tuned so that the ellipsoids $\{\mathcal{E}_{\le}(c_i,\delta_i E_i)\}_{i\in\mathbb{I}}$ are weakly pairwise disjoint. Then, the obstacle-free set $\mathcal{K}$ in~\eqref{eq:K} is forward invariant.
\end{theorem}
The existence of tuning parameters $\delta_1,\dots, \delta_I$ satisfying the weak pairwise disjointness of the sets $\{\mathcal{E}_{\le}(c_i,\delta_iE_i)\}_{i\in\mathbb{I}}$ is guaranteed by Assumption \ref{assumption:obstacle2}, which implies that weak pairwise disjointness holds when $\delta_i=1$ for all $i\in\mathbb{I}$. Hence, by a continuity argument, we can always tune each $\delta_i$ sufficiently close to $1$ in order to guarantee the weak pairwise disjointness of the dilated obstacles $\{\mathcal{E}_{\le}(c_i,\delta_iE_i)\}_{i\in\mathbb{I}}$. 
Note that algebraic tests of weak pairwise disjointness (provided in \cite[Thm.~6]{choi2006continuous} for $n=2$ and in \cite[Thm.~8]{wang2001algebraic} for $n=3$) can be used for this tuning purpose. 
\subsection{Global Asymptotic Stability}
\label{section:main:GAS}
\begin{figure}
    \centering
    \includegraphics[width=\columnwidth]{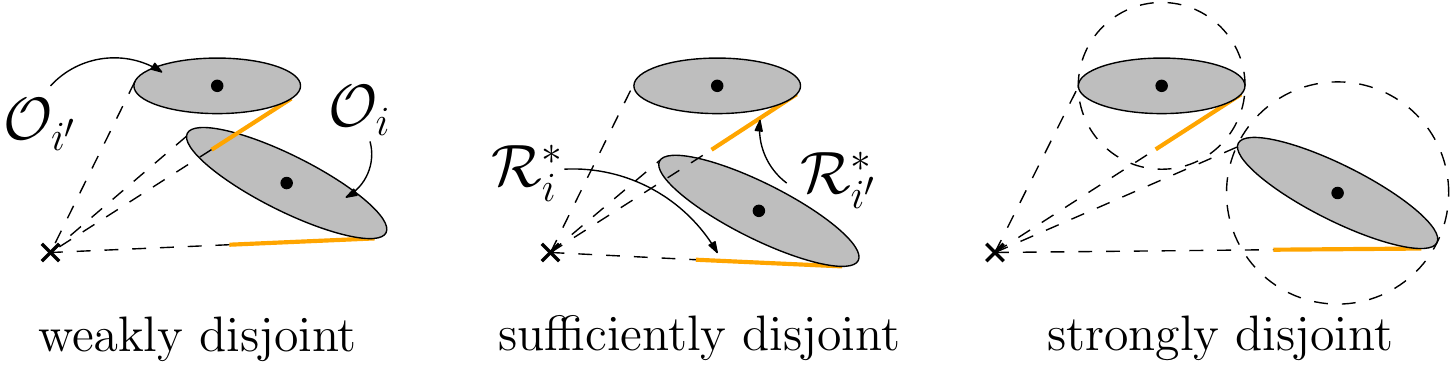}
    \caption{Different types of disjointness introduced in the paper with set $\mathcal{R}_i^*$ (orange, see~\eqref{eq:Ri}). For global attractivity, sufficient disjointness is asked.}
    \label{fig:different_disjoint}
\end{figure}

In this section we show that from all initial conditions in
the free workspace, all solutions converge asymptotically to the origin. 
To this end, we define the notion of \emph{sufficient disjointness} of a set of ellipsoids, 
which is slightly stronger than weak disjointness but less conservative than strong disjointness, and guarantees that each obstacle is avoided at most one time. The motivation behind the assumption of sufficient disjointness is that the ellipsoids considered here can
be arbitrarily large and flat, which might lead to long detours during the avoidance mode that take
the vehicle far away from the origin. In this case, specific configurations of the obstacles exist
such that from a set of initial conditions, the vehicle does not converge to the origin although it remains safe. 
Similarly, in the Bug~$0$ planning algorithm \cite{lumelsky1986dynamic}, termination 
(i.e., convergence to the target) is not \emph{always} guaranteed since the algorithm is
designed to ``walk toward the target whenever you
can'' \cite{lumelsky1986dynamic}. Our hybrid feedback shares a similar philosophy
since the vehicle jumps from avoidance to stabilization mode whenever the stabilization controller generates a vector field not
pointing towards the obstacle (see~\eqref{eq:dist_to_c_decreases}). To proceed, the next lemma characterizes the intersection of two ellipsoids of interest.  
\begin{lemma}\label{lemma:intersection}
Consider an arbitrary $i \in \mathbb{I}$. For $\underline{\delta}_i$, $\delta \mapsto \bar\mu_i(\delta)$ and $(\delta,\mu) \mapsto \bar\theta_i(\delta,\mu)$ defined in~\eqref{eq:bar delta,bar mu,bar theta}, let $\delta\in[\underline{\delta}_i,1]$, $\mu\in[1,\bar\mu_i(\delta)]$ and $\vartheta_i(\delta,\mu)$ be such that
\begin{equation}
\label{eq:cosine expr}
    \cos(\vartheta_i(\delta,\mu)):=\frac{1-\cos(\bar\theta_i(\delta,\mu))\underline{\delta}_i^2}{\sqrt{(1+\mu^{-2})/2-\delta^{-2}\underline{\delta}_i^4}}.
\end{equation}
The expression in~\eqref{eq:cosine expr} is well-defined and positive, and
\begin{equation}
\label{eq:ellips inters contained in cone}
    \mathcal{E}(c_i,\delta E_i)\cap\mathcal{E}(\bar c_i,\mu\bar E_i)\subset\mathcal{C}(0,c_i,\vartheta_i(\delta,\mu),E_i).
\end{equation}
\end{lemma}
Let us consider for each obstacle $i\in\mathbb{I}$ the sphere $\mathcal{S}(0,\bar r_i)$ with center at the origin and radius $\bar r_i$ defined by the next quadratic optimization problem
\begin{equation}\label{eq:rBari}
    \begin{aligned}
       \bar r_i^2:=&\min \|x\|^2\quad \textrm{subject to }&&x\in\mathcal{H}^*_i
    \end{aligned}
\end{equation}
where $\mathcal{H}^*_i$ is the helmet defined in \eqref{eq:helmet}. The radius $\bar r_i$ defines the minimum distance from the helmet $\mathcal{H}^*_i$ to the origin. Let $x$ be a point belonging to the intersection of the two ellipsoids $\mathcal{E}(c_i,E_i)$ and $\mathcal{E}(\bar c_i,\bar E_i)$. 
Taking $\delta$ and $\mu$ equal to $1$ in Lemma~\ref{lemma:intersection}, one obtains $x\in\mathcal{C}(0,c_i,\bar\vartheta_i,E_i)$ with 
\begin{equation}
\label{eq:cos bar vartheta}
\cos(\bar\vartheta_i) := \cos(\vartheta_i(1,1)) = \sqrt{1-\|E_ic_i\|^{-2}},
\end{equation}
from~\eqref{eq:cosine expr}, \eqref{eq:bar theta} and \eqref{eq:bar delta}. Now, let us define the set
\begin{equation}
\label{eq:Ri}
    \mathcal{R}^*_i:=\mathcal{C}(0,c_i,\bar\vartheta_i,E_i) \cap \mathcal{S}_\geq(0,\bar r_i)\cap\mathcal{E}_{\geq}(c_i,E_i)\cap\mathcal{E}_{\leq}(\bar c_i,\bar E_i),
\end{equation}
whose geometry is sketched in Fig.~\ref{fig:different_disjoint}. Intuitively speaking, it is a subset of all points on the cone $\mathcal{C}(0,c_i,\bar\vartheta_i,E_i)$ that have a distance to the origin greater than the distance $\bar r_i$ of the helmet $\mathcal{H}^*_i$ to the origin. The idea is that the vehicle should not to start avoiding another obstacle while it is still in $\mathcal{R}^*_i$, otherwise there is no guarantee that the number of times the vehicle avoids the obstacles is bounded and that global attractivity holds. This motivates the next definition. 
\begin{definition}
\label{def:suff pairwise disjoint}
The ellipsoids $\{\mathcal{E}(c_i,E_i)\}_{i\in\mathbb{I}}$ are sufficiently pairwise disjoint if they are weakly pairwise disjoint and 
\begin{equation}\label{eq:sufficient_disjoint}
     \forall i,i^\prime \in\mathbb{I} \text{ with } i\neq i^\prime,\quad \mathcal{R}^*_i\cap\mathcal{E}_\leq(c_{i^\prime},E_{i^\prime})=\emptyset.
\end{equation}
\end{definition}

\begin{figure}
    \centering
    \includegraphics[width=0.65\columnwidth]{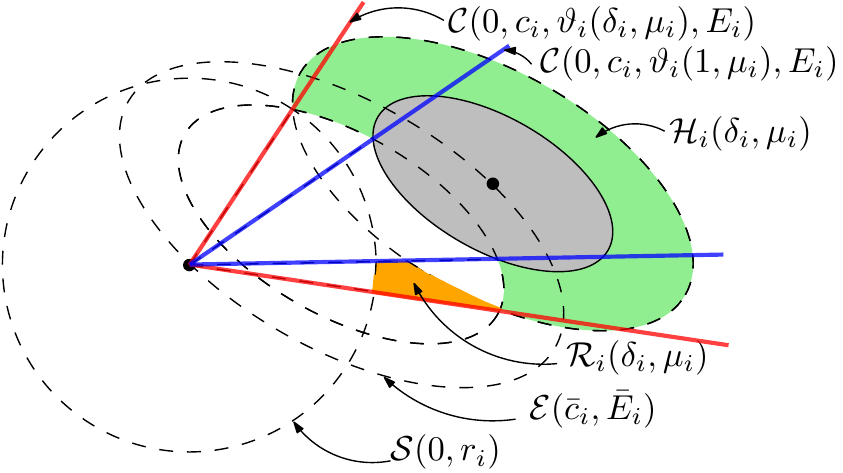}
    \caption{Safety helmet $\mathcal{H}_i(\delta_i,\mu_i)$ (green) and the corresponding escape region $\mathcal{R}_i(\delta_i,\mu_i)$ (orange). The region $\mathcal{R}_i(\delta_i,\mu_i)$ must not intersect with any other jump set $\mathcal{J}_0^{i^\prime}$, $i^\prime\neq i$, to avoid starting another avoidance while the distance to target has not yet decreased.}
    \label{fig:Ri}
\end{figure}

Now, let us introduce the ingredients for a dilated version of $\mathcal{R}^*_i$ as in~\eqref{eq:Ri_deltai_mui} below and refer to Fig.~\ref{fig:Ri}. First, consider the \textit{escape annulus} cone where solutions escape from the avoidance mode by applying the stabilization vector field. This region lies between the two cones $\mathcal{C}(0,c_i,\vartheta_i(1,\mu_i),E_i)$ and $\mathcal{C}(0,c_i,\vartheta_i(\delta_i,\mu_i),E_i)$ which are related, according to Lemma \ref{lemma:intersection}, to the intersections $\mathcal{E}(c_i,E_i)\cap\mathcal{E}(\bar c_i,\mu_i\bar E_i)$ and $\mathcal{E}(c_i,\delta_iE_i)\cap\mathcal{E}(\bar c_i,\mu_i\bar E_i)$, respectively. Second, consider for each obstacle $i\in\mathbb{I}$ the ball $\mathcal{S}_\geq(0,r_i)$ where the radius $r_i$ is defined by the quadratic optimization problem
\begin{equation}
\label{eq:ri}
    \begin{aligned}
       r_i^2:=&\min \|x\|^2\quad\textrm{subject to }&&x\in\mathcal{H}_i(\delta_i,\mu_i).
    \end{aligned}
\end{equation}
Note the following on~\eqref{eq:ri}. 1) The safety helmet $\mathcal{H}_i(\delta_i,\mu_i)$ is compact and, hence, the solution to~\eqref{eq:ri} exists. 2) For each $i\in\mathbb{I}$, $r_i>0$. Indeed, for each $i\in\mathbb{I}$, $\|\delta_iE_ic_i\|=\delta_i\underline{\delta}_i^{-2}>\delta_i\underline{\delta}_i^{-1}>1$ by Assumption~\ref{assumption:obstacle} and the selection of $\delta_i$ in Table~\ref{table:parameters}, so that $0\notin\mathcal{E}_\leq(c_i,\delta_iE_i)$ and in turn $0\notin\mathcal{H}_i(\delta_i,\mu_i)$ ($\mathcal{H}_i(\delta_i,\mu_i)\subset\mathcal{E}_\leq(c_i,\delta_iE_i)$). Hence, since $\mathcal{H}_i(\delta_i,\mu_i)$ is compact there exists $r_i>0$ such that $\|x\|\geq r_i$ for all $x\in\mathcal{H}_i(\delta_i,\mu_i)$. Finally, we can define the considered dilated version of $\mathcal{R}^*_i$ as
\begin{multline}
\label{eq:Ri_deltai_mui}
    \mathcal{R}_i(\delta_i,\mu_i):=\mathcal{S}_\geq(0,r_i)\cap\mathcal{E}_{\geq}(c_i,\delta_iE_i)\cap\mathcal{E}_{\leq}(\bar c_i,\bar E_i)\\\cap\mathcal{C}_{\geq}(0,c_i,\vartheta_i(1,\mu_i),E_i)\cap\mathcal{C}_{\leq}(0,c_i,\vartheta_i(\delta_i,\mu_i),E_i).
\end{multline}
\begin{lemma}\label{lemma:disjoint}
Assume that the obstacles $\{\mathcal{O}_i\}_{i\in\mathbb{I}}$ are sufficiently pairwise disjoint. Then, for each $i\in\mathbb{I}$, there exist $\delta_i^*, \mu_i^*$ such that for all $\delta_i\in(\delta_i^*,1)$ and $\mu_i\in(1,\mu_i^*)$, we have
\begin{equation}\label{eq:Ri:disjoint}
     \forall i',i'' \in\mathbb{I}, i'\neq i'',\,\, \mathcal{R}_{i'}(\delta_{i'},\mu_{i'})\cap\mathcal{E}_\leq(c_{i''},\delta_{i''} E_{i''})=\emptyset.
\end{equation}
\end{lemma}
Property \eqref{eq:Ri:disjoint} of Lemma \ref{lemma:disjoint} will be used to show global attractivity. Intuitively, we require that after avoiding an obstacle, the distance $\|x\|$ to the target decreases before the vehicle reaches the proximity of another obstacle. Although the bounds $\delta_i^*$ and $\mu_i^*$ are not defined explicitly for generic ellipsoids, the parameters $\delta_i$ and $\mu_i$ can be tuned offline. Now, we are ready to state our main result for this section.
\begin{theorem}\label{theorem:main}
Consider the hybrid system \eqref{eq:hs} under the same assumptions as Theorem~\ref{theorem:invariance}. Assume also that the obstacles $\{\mathcal{O}_i\}_{i\in\mathbb{I}}$ are sufficiently pairwise disjoint, and $\delta_i$ and $\mu_i$ are tuned such that \eqref{eq:Ri:disjoint} holds. Then, the set $\A:=\{0\}\times\mathbb{I}\times\mathbb{M}$ is globally asymptotically stable for~\eqref{eq:hs} and the number of jumps is bounded.
\end{theorem}
%


For spherical obstacles, we show next that the extra tuning of the parameters to satisfy \eqref{eq:Ri:disjoint} is not needed.
\begin{theorem}(Spherical obstacles)\label{theorem:spherical}
Let $E_i=\lambda_iI_n$ for all $i\in\mathbb{I}$. Under the same assumptions as Theorem~\ref{theorem:invariance}, the set $\A:=\{0\}\times\mathbb{I}\times\mathbb{M}$ is globally asymptotically stable for~\eqref{eq:hs} and the number of jumps is bounded.
\end{theorem}

\subsection{Complementary Properties}
\label{section:main:complProp}

In this section we present four relevant complementary properties of the proposed hybrid law for obstacle avoidance. 
\subsubsection{Bounded Control}
\label{boundedControl}
First, we can show that $x$ remains always in a given ball. Indeed, let $\mathcal{S}_\leq(0,r_b)$, with $r_b>0$, be the smallest ball containing all the dilated ellipsoids $\mathcal{E}(c_i,\delta_iE_i)$ (which must exist since these ellipsoids are compact). During stabilization mode the distance $\|x\|$ is decreasing and during avoidance mode the vehicle stays within the dilated ellipsoids $\mathcal{E}(c_i,\delta_iE_i)$. Then, it is guaranteed that from all $x(0,0)\in\mathcal{S}_\leq(0,r_b)$, $x(t,j)\in\mathcal{S}_\leq(0,r_b)$ for all $(t,j)\in\dom x$.  Moreover, since the projection matrix $\pi^\perp(E_i(x-c_i))$ has eigenvalues in $0$ and $1$, it follows that we can upper bound the control input in \eqref{eq:u} by
$
\|u\|\leq k\alpha(r_b+p)
$
where $k=\max\{k_1,k_0,k_{-1}\}$, $\alpha=\max_{i\in\mathbb{I}} (\lambda_{\max}(E_i)/\lambda_{\min}(E_i))$ and $p=\max_{i\in\mathbb{I}}\|p_1^i\|$. The control gains can then be tuned to satisfy the inherent practical saturation of the actuators.
\subsubsection{Semiglobal Preservation}
The second property is the so-called {\it semiglobal preservation} property  \cite[\S II]{braun2018unsafe}.  This  property is desirable when the original controller parameters are optimally tuned and the controller modifications imposed by the presence of  the  obstacles  should  be  as  minimal  as  possible. Such a property is also accounted for in the quadratic programming formulation of~\cite[III.A.]{wang2017safety}. We summarize this property for our case in the next proposition.
\begin{proposition}\label{proposition:semiglobal_preservation}
Be $\epsilon\in(0,1)$ and $\mathcal{W}_\epsilon:=\bigcap_{i^\prime \in \mathbb{I}}\mathcal{E}_\geq(c_{i^\prime},\epsilon E_{i^\prime})$. There exist controller parameters such that the control law matches, in  $\mathcal{W}_\epsilon$, the stabilization feedback $u=-k_0x$  ($k_0>0$) used in the absence of obstacles.
\end{proposition}
\subsubsection{Non-point Mass Vehicles} There is no loss of generality in considering a point-mass vehicle in this work. Let us rather consider that the vehicle has some volume, {\it e.g.,} bounded by $\mathcal{S}_\leq(x,r_v)$. 
Then, for the navigation scenario to be feasible, the radius $r_v$ of the vehicle needs to be smaller than the smallest distance between the obstacles, {\it i.e.,} for all $i,i^\prime \in \mathbb{I}$ with $i \neq i^\prime$, $r_v<\mathbf{dist}(\mathcal{E}_\leq(c_i,E_i),\mathcal{E}_\leq(c_{i^\prime},E_{i^\prime}))$. For the safety of the vehicle during the stabilization mode, selecting the parameter $\epsilon_i$ as $\epsilon_i<(1+\lambda_{\max}(E_i)r_v)^{-1}$ is sufficient (in addition to Table~\ref{table:parameters}) to guarantee that the vehicle starts the avoidance mode away from the obstacle. 
Indeed, under this condition, it is easy to show that for all $x\in\mathcal{E}_{\geq}(c_i,\epsilon_i E_i)$ ({\it i.e.,} the vehicle center is outside the dilated ellipsoid $\mathcal{E}(c_i,\epsilon_i E_i)$) and for all $x^\prime\in\mathcal{S}_\leq(x,r_v)$, one has $x^\prime\in\mathcal{E}_\geq(c_i,E_i)$, which guarantees safety of the whole volume of the vehicle.

\subsubsection{Robustness}
The constructed hybrid controller guarantees some level of robustness to perturbations (e.g., in the form of measurement noise).
Hysteresis switching is one of the typical ways to ensure robustness to measurement noise, and hysteresis switching is indeed behind the designed hybrid feedback, in particular the hysteresis regions of flow and jump sets in Section~\ref{section:construction} and the logical selections of the jump sets in Section~\ref{section:selection}. 
More generally, fundamental results in \cite[Chap.~7]{goebel2012hybrid} guarantee structurally that  global asymptotic stability of $\A$ in Theorem~\ref{theorem:main} is also uniform (by \cite[Thm.~7.12]{goebel2012hybrid}) and robust (by \cite[Thm.~7.21]{goebel2012hybrid}) with respect to perturbations since $\A$ is a compact set and the hybrid basic conditions are satisfied as in Lemma~\ref{lemma:hbc}.

\section{Simulations}
\label{section:example}
\begin{figure}
    \centering
    \includegraphics[width=0.49\columnwidth]{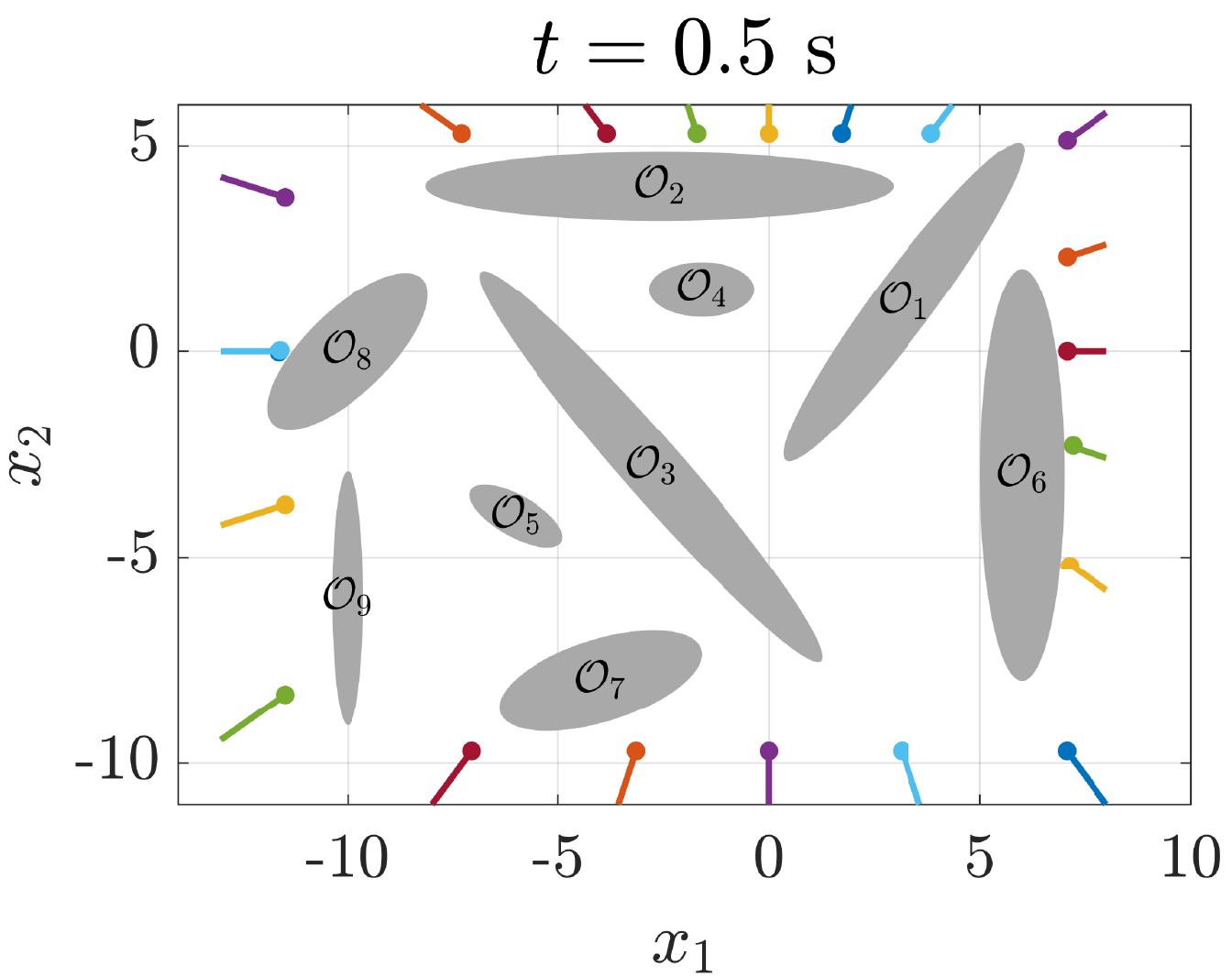}
     \includegraphics[width=0.49\columnwidth]{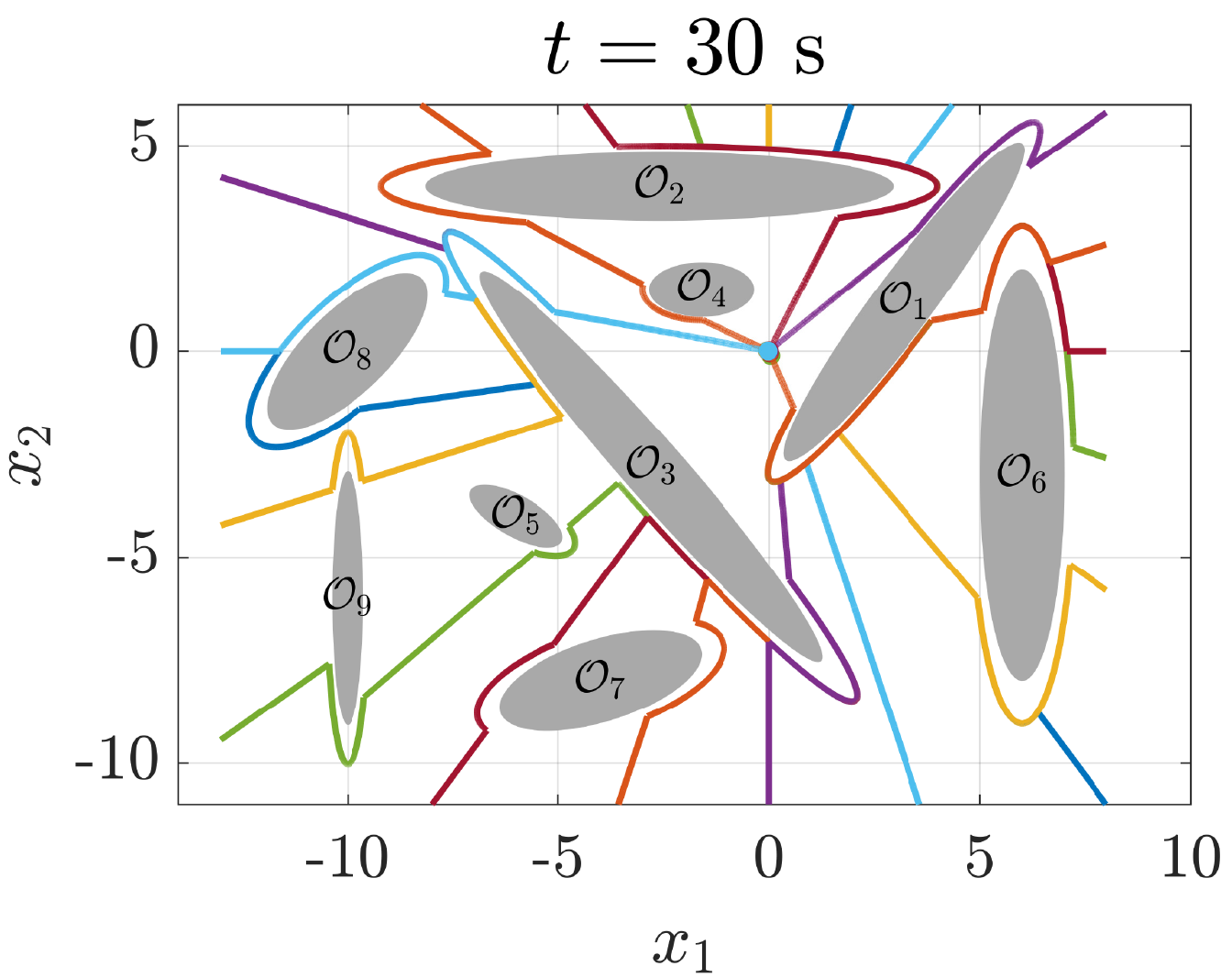}
    \caption{Plot (at time $t=0.5$ and $t=30$ seconds) of the $2$-dimensional trajectory of the vehicle starting at different initial conditions. }
    \label{fig:2Dobstacles}
\end{figure}

\begin{figure}
    \centering
    \includegraphics[width=0.7\columnwidth]{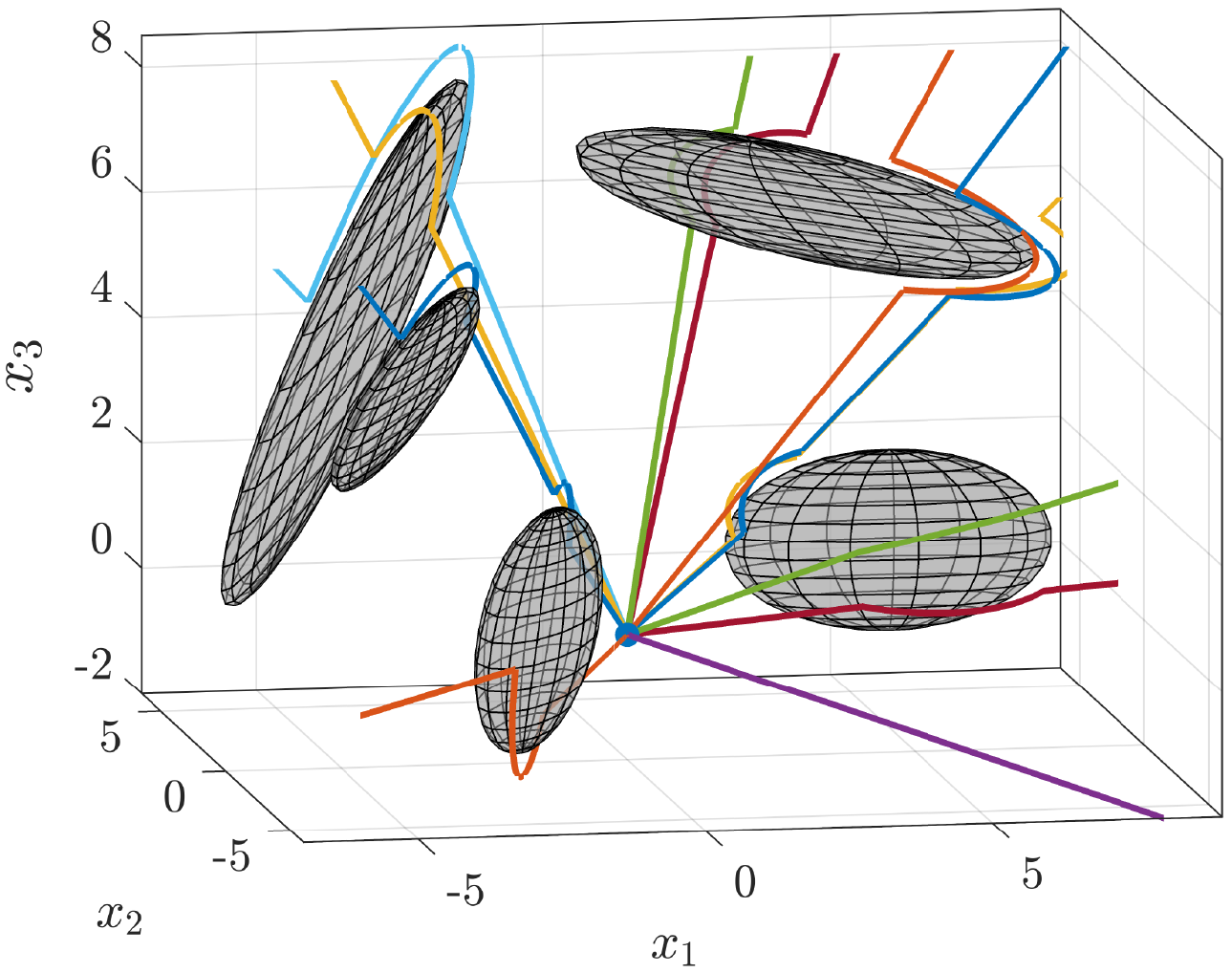}
    \caption{Plot (at time $t=30$ seconds) of the $3$-dimensional trajectory of the vehicle starting at different initial conditions.}
    \label{fig:3Dobstacles}
\end{figure}

We illustrate the effectiveness of the proposed hybrid control strategy through two simulation scenarios. The first scenario considers $9$ obstacles in $2$D (see Fig.~\ref{fig:2Dobstacles}) while the second one considers $5$ obstacles in $3$D (see Fig.~\ref{fig:3Dobstacles}). 
For both cases, Table~\ref{table:parameters} provides a suitable order to choose the parameters for each $i \in \mathbb{I}$, as follows.
\begin{enumerate}[leftmargin=12pt,label=\arabic*)]
\item \label{guideline step 1} For $\underline{\delta}_i$ in~\eqref{eq:bar delta}, select $\delta_i$ and $\epsilon_i$ so that $\underline{\delta}_i < \delta_i < \epsilon_i < 1$; 
\item \label{guideline step 2} For $\delta_i$ and $\bar\mu_i(\delta_i)$ in~\eqref{eq:bar mu}, select $\nu_i$ and $\mu_i$ so that $1<\nu_i<\mu_i<\bar\mu_i(\delta_i)$ (possibly iterating steps \ref{guideline step 1} and \ref{guideline step 2} so that $\delta_i$ and $\mu_i$ satisfy \eqref{eq:Ri:disjoint});
\item For $\delta_i$, $\mu_i$ and $\bar\theta_i(\delta_i,\mu_i)$ in~\eqref{eq:bar theta}, select $\psi_i$, $\bar \psi_i$ and $\theta_i$ so that $0<\psi_i<\bar \psi_i<\theta_i< \bar\theta_i(\delta_i,\mu_i)$.
\end{enumerate}
Any parameter selection according to this guideline guarantees our results, and can be carried out keeping in mind the physical interpretation illustrated in Section~\ref{section:construction} for these parameters.
The gains are $k_0=k_1=k_{-1}=1/4$ and determine the speed of convergence of the scheme. By~\eqref{eq:p_1}, the point $p_1^i$ can be selected arbitrarily as long as it is on $\mathcal{C}(c_i,-c_i,\theta_i,E_i)\backslash\{c_i\}$. A suitable choice is given by
\begin{equation}
\label{eq:p1i construction}
p_1^i=\pi^\perp(E_i^{-1}\mathbf{R}(\theta_i)E_ic_i)c_i
\end{equation}
where $\mathbf{R}(\theta_i)$ is the standard $2\times 2$ rotation matrix with angle $\theta_i$ or the standard $3\times 3$ axis-angle rotation matrix with angle $\theta_i$ and an arbitrary vector of $\mathbb{S}^2$ as axis.
The idea behind~\eqref{eq:p1i construction} is to project $c_i$ on the plane orthogonal to a rotated version of $c_i$, in order to obtain the point lying on the cone and closest to the origin. Having all points $p_m^i$ close enough to the origin is an effective way so that $k_0$, $k_1$, $k_{-1}$ can take the same values and yield comparable speeds for avoidance and stabilization, independently of the obstacles.

\begin{figure}
    \centering
    \includegraphics[width=0.49\columnwidth]{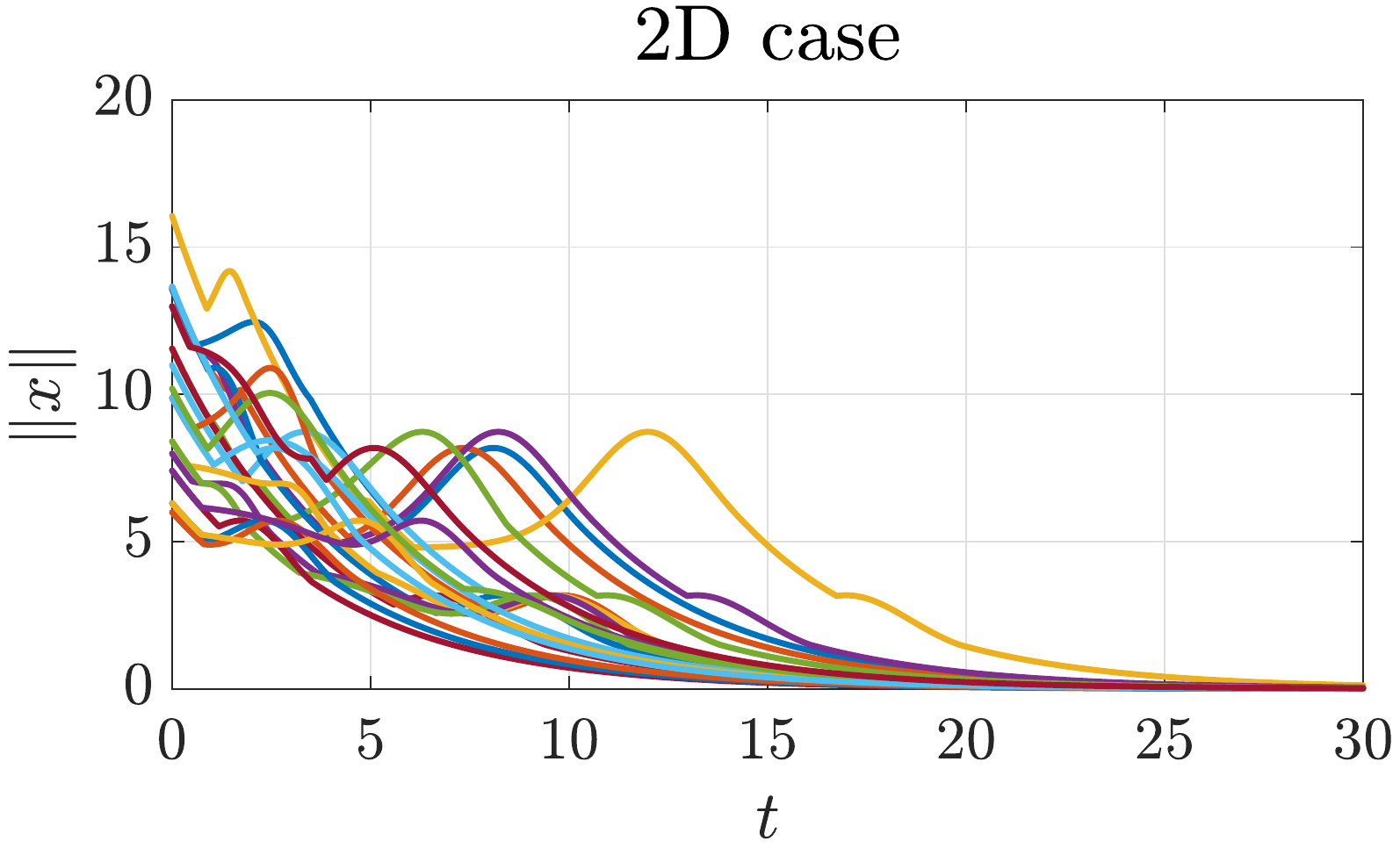}
    \includegraphics[width=0.49\columnwidth]{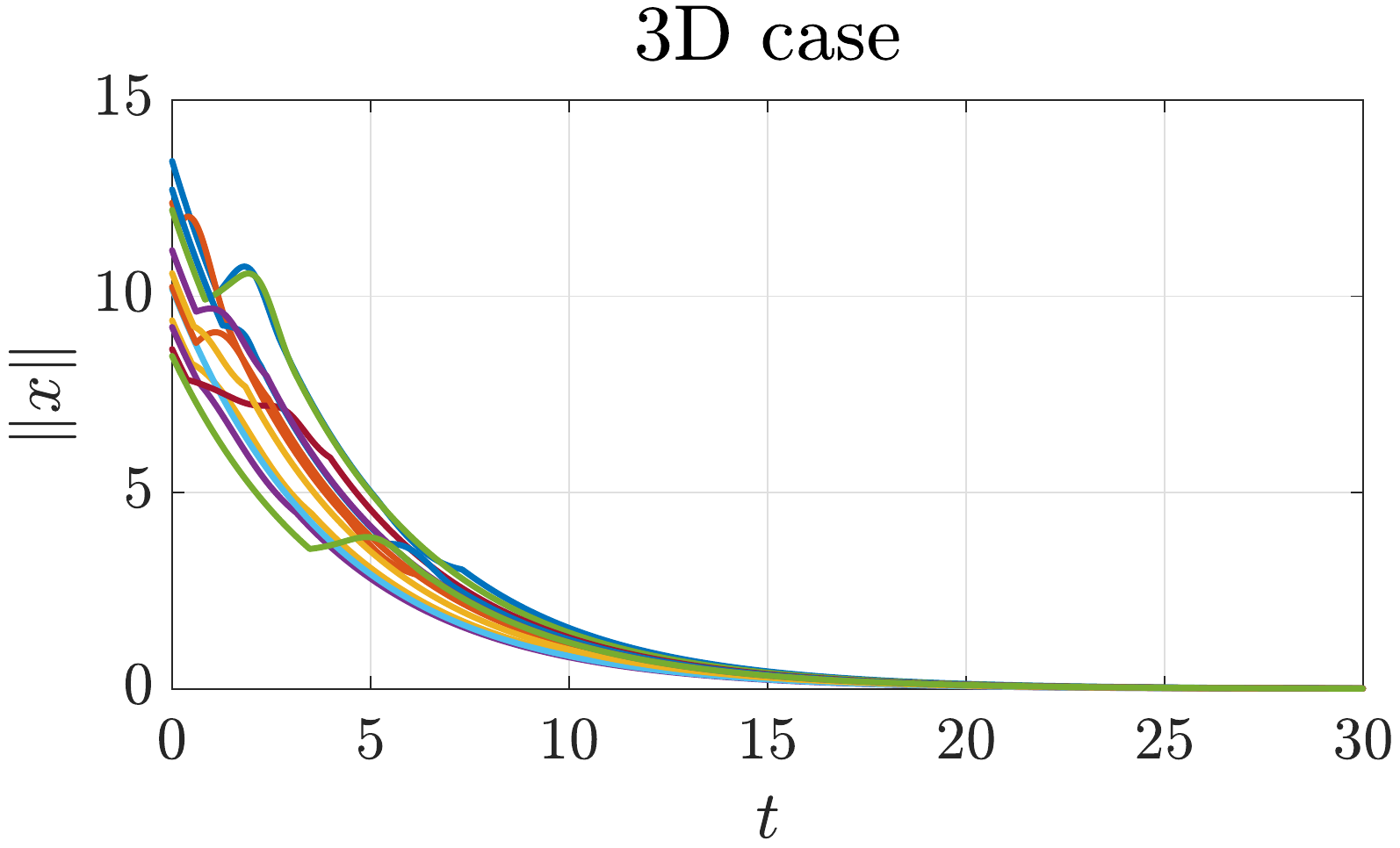}
    \caption{Plot of the position norm $\|x\|$ versus time showing the convergence of the solutions to the origin from the considered initial conditions.}
    \label{fig:2Dconvergence}
\end{figure}

Fig.~\ref{fig:2Dobstacles} (Fig.~\ref{fig:3Dobstacles}, respectively) shows that the solution generated by the closed loop hybrid system avoids the $2$D obstacles ($3$D obstacles, respectively) and Fig.~\ref{fig:2Dconvergence} shows the convergence of solutions to the origin. Complete simulation videos for the $2$D and $3$D cases can be found at \url{https://youtu.be/CnXJlhzlzd8}, \url{https://youtu.be/4mzTXPR6D9Y}.

Finally, we note that for the very obstacle configuration of the 2D scenario, the state-of-the-art approach of navigation functions  \cite{koditschek1990robot,paternain2018navigation} cannot be applied since the condition \cite[Thm.~3, Eq.~(23)]{paternain2018navigation} is violated for all obstacles except obstacle $\mathcal{O}_5$, where \cite[Eq.~(23)]{paternain2018navigation} intuitively corresponds to the fact that obstacles are not too flat and not too close to the target position. 
(\cite[Eq.~(23)]{paternain2018navigation} is violated for all obstacles of the 3D scenario.)  
Moreover, navigation function approaches require tuning a parameter sufficiently large ($k$ in \cite[Eq.~(17) and Remark~5]{paternain2018navigation}), which may conflict with actuator limitations. Instead, our approach provides a clear tuning guideline for all parameters (given in this section) and actuator limitations can be taken into account (see Section~\ref{boundedControl}).

\section{Conclusions}
We proposed a novel hybrid feedback on $\mathbb{R}^n$ to solve the obstacle avoidance problem for arbitrarily flat ellipsoidal obstacles. Our control strategy ensures global asymptotic stabilization to the target and safety (thus, successful navigation from all initial conditions) while guaranteeing a Zeno-free switching between the avoidance and stabilization modes. Moreover, the control input remains bounded (in particular, arbitrarily close to any obstacle) and matches semi-globally in the free-state space the nominal feedback used in the absence of obstacles. 

Future work will be devoted to considering more complex vehicle dynamics ({\it e.g.}, under-actuated and second-order dynamics) and more generic obstacle shapes ({\it e.g.}, convex obstacles). Furthermore, although our scheme considers static obstacles to obtain formal guarantees for global asymptotic stability and safety, extending this approach to deal with unknown environments is an interesting research direction that we aim at pursuing in the future.


\appendix
\label{section:appendix}
In the appendix, an equation number over an (in)equality indicates which equation has been used to obtain the (in)equality.
\subsubsection{Proof of Theorem \ref{theorem:invariance}}
Define $\mathbf{S}_\mathscr{H}(\mathcal{K})$ as the set of all maximal solutions $\phi$ to $\mathscr{H}$ with $\phi(0,0) \in \mathcal{K}$. 
Each $\phi\in\mathbf{S}_\mathscr{H}(\mathcal{K})$ has range $\rge\phi\subset\mathcal{K}=\mathcal{F}\cup\mathcal{J}$ by Lemma~\ref{lemma:union} and the definition of hybrid solution~\cite[p.~124]{goebel2012hybrid}, so $\mathcal{K}$ is forward pre-invariant \cite[Def.~3.3]{chai2018forward}. 
The set $\mathcal{K}$ is in fact forward invariant \cite[Def.~3.3]{chai2018forward} if for each $\xi \in \mathcal{K}$ there exists one solution and each $\phi\in\mathbf{S}_\mathscr{H}(\mathcal{K})$ is complete, which we show in the rest of the proof through~\cite[Prop.~6.10]{goebel2012hybrid}. 
In the rest of the proof, let
\begin{equation}
\label{eq:F0J0star}
    \mathcal{F}_0^*:=\bigcap_{i\in\mathbb{I}}\mathcal{F}_0^i, \quad\mathcal{J}_0^*:=\bigcup_{i\in\mathbb{I}}\mathcal{J}_0^i.
\end{equation}
\begin{lemma}\label{lemma:boundary}
Under the assumptions of Theorem \ref{theorem:invariance}, we have for each $i \in \mathbb{I}$ and $m\in\{-1,1\}$
\begin{subequations}
\begin{align}
        \label{eq:J0i_2}
      &\mathcal{J}_0^i=\mathcal{H}_i(\epsilon_i,\nu_i),\\
      \label{eq:Fmi}
       &\mathcal{F}_m^i=\mathcal{H}_i(\delta_i,\mu_i)\cap\mathcal{C}_\geq(c_i,c_i-p_m^i,\psi_i,E_i), \\
      & \partial\mathcal{F}_0^*\backslash\mathcal{J}_0^* \subset \bigcup_{i\in\mathbb{I}}\big(\mathcal{E}(c_i,E_i)\backslash\mathcal{E}_\geq(\bar c_i,\bar E_i)\big) \label{eq:boundary F_0 star},\\
      &\begin{aligned}
      \partial\mathcal{F}_m^i\backslash\mathcal{J}_m^i \subset & \mathcal{E}(c_i,E_i)\backslash(\mathcal{E}_\leq(\bar c_i,\mu_i\bar E_i)\\
    &\cup\mathcal{C}_\leq(c_i,c_i\!-\!p_m^i,\psi_i,E_i)).        
      \end{aligned} \label{eq:boundary F^i_m}
\end{align}
\end{subequations}
\end{lemma}
First, let us show that the viability condition
\begin{align}\label{eq:viability}
    \mathbf{F}(x,i,m)\cap\mathbf{T}_{\mathcal{F}}(x,i,m)\neq\emptyset
\end{align}
holds for all $(x,i,m)\in\mathcal{F}\backslash\mathcal{J}$. Let $(x,i,m)\in\mathcal{F}\backslash\mathcal{J}$, which implies by~\eqref{eq:FJ} that $(x,i)\in\mathcal{F}_m\backslash\mathcal{J}_m$ for some $m\in\mathbb{M}$, and divide into the cases $m=0$ and $m \in \{-1,1\}$. When $m=0$, from~\eqref{eq:FJ0} there exists $i\in\mathbb{I}$ such that $x\in\mathcal{F}_0^*\backslash\mathcal{J}_0^*$. If $x\in(\mathcal{F}_0^*)^\circ\backslash\mathcal{J}_0^*$ (hence, $x$ is in the interior of $\mathcal{F}_0^*$), then $\mathbf{T}_{\mathcal{F}_0^*}(x)=\mathbb{R}^n$, so that $\mathbf{T}_{\mathcal{F}}(\xi)=\mathbb{R}^n\times\{0\}\times\{0\}$ and \eqref{eq:viability} holds.
If $x\in\partial\mathcal{F}_0^*\backslash\mathcal{J}_0^*$, which satisfies the set inclusion \eqref{eq:boundary F_0 star}, the weak pairwise disjointness of $\{\mathcal{E}(c_i,E_i)\}_{i\in\mathbb{I}}$ yields:
\begin{equation}
\label{eq:tangent cone m=0}
\begin{aligned}
& x\in\mathcal{E}(c_i,E_i),\, i\in\mathbb{I}  \\
& \mathbf{T}_{\mathcal{F}}(x,i,0)=\mathcal{P}_\geq(0,E_i^2(x-c_i))\times\{0\}\times\{0\}.
\end{aligned}
\end{equation}
By~\eqref{eq:dist_to_c_decreases} and $x\notin\mathcal{E}_\geq(\bar c_i,\bar E_i)$ by~\eqref{eq:boundary F_0 star}, we obtain
\begin{equation}
-\!k_0x^\top\!E_i^2(x\!-\!c_i)=k_0\|E_i \bar c_i\|^2\!\left(1\!-\!\|\bar E_i(x-\bar c_i)\|^2\right)\!>\!0,
\end{equation}
hence, $\kappa(x,i,0)\in\mathcal{P}_\geq(0,E_i^2(x-c_i))$ in~\eqref{eq:tangent cone m=0}, and \eqref{eq:viability} holds for $m=0$. When $m\in\{-1,1\}$, we have $i\in\mathbb{I}$ and $x\in\partial\mathcal{F}_m^i\backslash\mathcal{J}_m^i$, which satisfies the set inclusion \eqref{eq:boundary F^i_m}, and so 
\begin{equation}
\label{eq:tangent cone |m|=1}
    \mathbf{T}_{\mathcal{F}}(x,i,m)=\mathcal{P}_\geq(0,E_i^2(x-c_i))\times\{0\}\times\{0\}.
\end{equation}
$\kappa(x,i,m)\in\mathcal{P}_\geq(0,E_i^2(x-c_i))$ in~\eqref{eq:tangent cone |m|=1} because
\begin{equation}
\label{eq:vector field tangent cone |m|=1}
-k_m(x-p_m^i)^\top E_i\pi^\perp(E_i(x-c_i))E_i(x-c_i)=0,
\end{equation}
so the viability condition \eqref{eq:viability} holds for $m\in\{-1,1\}$ as well. 

Second, we apply~\cite[Prop.~6.10]{goebel2012hybrid}. By it and \eqref{eq:viability}, there exists a nontrivial solution to $\mathscr{H}$ from each initial condition in $\mathcal{K}$. Finite escape times can only occur through flow. They can neither occur for $x$ in the set $\mathcal{F}_{-1}^i \cup \mathcal{F}_1^i$ ($\mathcal{F}_{-1}^i$ and $\mathcal{F}_1^i$ are bounded by their definitions in~\eqref{eq:F-11i}) nor for $x$ in the set $\mathcal{F}_{0}^*$ because they would make $x^\top x$ grow unbounded, and this would contradict that $\tfrac{d}{dt}(x^\top x) \le 0$ by the definition of $\kappa(x,i,0)$ and by~\eqref{eq:flow}. So, all maximal solutions do not have finite escape times. By Lemma~\ref{lemma:union}, $\mathbf{J}(\mathcal{J})\subset\mathcal{K}= \mathcal{F} \cup \mathcal{J}$. Hence, by \cite[Prop.~6.10]{goebel2012hybrid}, all maximal solutions are complete.
\subsubsection{Proof of Theorem \ref{theorem:main}}
We prove global asymptotic stability of $\A$ by~\cite[Def.~7.1]{goebel2012hybrid}. For each $i\in\mathbb{I}$, $\|\delta_iE_ic_i\|=\delta_i\underline{\delta}_i^{-2}>\delta_i\underline{\delta}_i^{-1}>1$ by Assumption~\ref{assumption:obstacle} and the selection of $\delta_i$ in Table~\ref{table:parameters}, so $0\notin\mathcal{E}_\leq(c_i,\delta_iE_i)$. As a consequence, there exists $\varepsilon^*>0$ such that the ball $\mathcal{S}_\leq(0,\varepsilon^*)$ does not intersect with any of the dilated obstacles $\mathcal{E}_\leq(c_i,\delta_iE_i)$. It can be shown easily that for each $\varepsilon\in [0,\varepsilon^*]$, the set $\mathcal{S}:=\mathcal{S}_\leq(0,\varepsilon)\times\mathbb{I}\times\mathbb{M}$ is forward invariant because $\mathcal{S}_\leq(0,\varepsilon)$ is disjoint from $\mathcal{J}_0^*$ and the component $x$ of solutions evolves, after at most one jump, with the stabilization mode $\dot x = - k_0 x$. 
Thanks to forward invariance of $\mathcal{S}$, stability of $\A$ for~\eqref{eq:hs} is immediate from \cite[Def.~7.1]{goebel2012hybrid}.
Let us prove global attractivity of $\A$.Before that, we need the next intermediate result.
\begin{lemma}\label{lemma:decrease}
There exists $\sigma>0$ such that for all solutions $\xi=(x,i,m)$ with $\xi(t,j)\in\mathcal{F}_l\times\{l\}$ for some $l\in\{-1,1\}$ and $(t,j)\in\dom\xi$, there exists $(s, \ell) \in \dom \xi$ such that $(s,\ell)\succeq(t,j)$ and
\begin{equation}
\label{eq:strict decrease}
\|x(s, \ell)  \|  \leq\|x(t,j)\| - \sigma.
\end{equation}
\end{lemma}
Now, for each solution $\xi$ to~\eqref{eq:hs}, there exists a finite time $(T,J)\succeq (0,0)$ after which the solution does not evolve with the avoidance controller any longer, {\it i.e.}, $m(t,j)=0$ for all $(t,j)\succeq (T,J)$. Otherwise, there would exist a sequence of hybrid times $\{(t_k,j_k)\}_{k=0}^\infty$ such that $\xi(t_k,j_k)\in\mathcal{F}_{l_k}\times\{l_k\}$ with $l_k\in\{-1,1\}$ and this would imply by Lemma~\ref{lemma:decrease} that $\|x(t_{k+1},j_{k+1})\|\leq\|x(t_k,j_k)\|-\sigma$ for all $k\in\mathbb{N}$. This is indeed a contradiction as it would lead to $\|x(\cdot,\cdot)\|$ becoming negative. Then, the solution $\xi$ enters the stabilizing mode $m=0$ after $(T,J)$ and its flow map $\dot x = - k_0 x$ guarantees in turn global attractivity. Moreover, $J$ is the maximum number of jumps of the hybrid system since any extra jump will cause $m$ to take values in $\{-1,1\}$, which is not possible after $(T,J)$.
\subsubsection{Proof of Theorem \ref{theorem:spherical}}
To prove the theorem, it is sufficient to show that for spherical obstacles the result of Lemma \ref{lemma:decrease} holds. The proof of Lemma \ref{lemma:decrease} under the assumptions of Theorem~\ref{theorem:spherical} is the same up to~\eqref{eq:intersection of tilde F^i_l and J^i_l}.
From~\eqref{eq:intersection of tilde F^i_l and J^i_l}  we have $x(t^\prime,j+1)=x(t^\prime,j) \in \tilde{\mathcal{P}}_{m(t^\prime,j),3}^{i(t^\prime,j)}(\delta,\psi)$, $i(t^\prime,j+1)= i(t^\prime,j)=:\iota$ and $m(t^\prime,j)=0$. However, since $x(t,j)$ and $x(t^\prime,j+1)$ both belong to $\mathcal{E}(c_\iota,\delta E_\iota)$, we can write $\|x(t,j)-c_\iota\|^2=\|x(t^\prime,j+1)-c_\iota\|^2$ (since $E_\iota=\lambda_\iota I_n$) and hence
\begin{align}\label{eq:thm3:1}
    \|x(t,j)\|^2-2c_\iota^\top x(t,j)=\|x(t^\prime,j+1)\|^2-2c_\iota^\top x(t^\prime,j+1). 
\end{align}
$x(t^\prime,j+1)\in  \tilde{\mathcal{P}}_{0,3}^{\iota}(\delta,\psi) \subset \mathcal{E}(\bar c_\iota,\mu_\iota \bar E_\iota)$ implies also that
\begin{equation*}
    c_\iota^\top x(t^\prime,j+1)=\|x(t^\prime,j+1)\|^2+(1-\mu_\iota^{-2})\|\bar c_\iota\|^2,
\end{equation*}
thus, with \eqref{eq:thm3:1}, we have
\begin{equation*}
    c_\iota^\top x(t,j)=\frac{ \|x(t,j)\|^2+\|x(t^\prime,j+1)\|^2}{2}+(1-\mu_\iota^{-2})\|\bar c_\iota\|^2.
\end{equation*}
However, since $x(t,j)\in\mathcal{E}_{\geq}(\bar c_\iota,\mu_\iota\bar E_\iota)$, we have
\begin{equation*}
    c_\iota^\top x(t,j)\leq\|x(t,j)\|^2+(1-\mu_\iota^{-2})\|\bar c_\iota\|^2,
\end{equation*}
and, hence, $\|x(t^\prime,j+1)\|^2\leq \|x(t,j)\|^2$ must hold. Also, by~\eqref{eq:intersection of tilde F^i_l and J^i_l}, $x(t^\prime,j+1)\in\mathcal{E}(c_\iota,\delta E_\iota)\cap\mathcal{E}(\bar c_\iota,\mu_\iota\bar E_\iota)$, and, by Lemma~\ref{lemma:intersection}, $x(t^\prime, j+1)\in\mathcal{C}(0,c_\iota,\vartheta_\iota(\delta,\mu_\iota),E_\iota)$. In view of Step~2 of the proof of Lemma \ref{lemma:decrease}, both the sets $\mathcal{E}_{\leq}(\bar c_\iota,\bar E_\iota)$ and $\mathcal{C}(0,c_\iota,\vartheta_\iota(\delta,\mu_\iota),E_\iota)$ are forward invariant under the stabilization flow map for $x$, i.e., $-k_0x$. Since the obstacles are weakly disjoint, the solution then flows in $\mathcal{E}_\le(\bar c_\iota,\bar E_\iota) \cap \mathcal{C}(0,c_\iota,\vartheta_\iota(\delta,\mu_\iota),E_\iota)$ until it reaches the set $\mathcal{E}(c_\iota,\delta_\iota E_\iota)$ at $(t'',j+1)$. Since flow with stabilization mode decreases the distance to the origin we have
\begin{equation*}
    \|x(t'',j+1)\|^2\leq\|x(t^\prime,j+1)\|^2\leq \|x(t,j)\|^2
\end{equation*}
Also, the solution must flow after $(t'',j+1)$ up to some $(s,j+1)$ with the stabilization mode (since obstacles are weakly disjoint) such that at least $\sigma_1$ in \eqref{eq:sigma1} is traversed, i.e.,
\begin{equation*}
   \|x(s,j+1)\|^2+\sigma_1\leq \|x(t'',j+1)\|^2\leq \|x(t,j)\|^2.
\end{equation*}
This proves Lemma \ref{lemma:decrease} and in turn Theorem \ref{theorem:main}.

\subsubsection{Proof of Proposition \ref{proposition:semiglobal_preservation}}
Note preliminarily that thanks to $\epsilon<1$, $\mathcal{W}_\epsilon\subset\mathcal{W}$ in~\eqref{eq:W}. 
It is sufficient to show that the closed loop system under the proposed hybrid feedback cannot flow except with stabilization mode $m=0$ when $x\in\mathcal{W}_\epsilon$. Indeed, if in Table~\ref{table:parameters} we further constrain $\delta_i$ as $\delta_i\in(\max(\underline{\delta}_i,\epsilon),1)$ for all $i\in\mathbb{I}$, then $\mathcal{F}_m^i\subset\mathcal{H}(\delta_i,\mu_i)\subset\mathcal{E}_\leq(c_i,\delta_i E_i)\subset\mathcal{E}_\leq(c_{i},\epsilon E_{i})$ and $\mathcal{E}_\leq(c_i,\delta_i E_i)\neq \mathcal{E}_\leq(c_{i},\epsilon E_{i})$. Therefore, we have $\mathcal{F}_m^i\cap\mathcal{W}_{\epsilon}=\emptyset$ for all $i\in\mathbb{I}$ and $m\in\{-1,1\}$. This implies that solutions cannot flow with the avoidance mode when $x$ belongs to $\mathcal{W}_\epsilon$ and must then flow with the stabilization mode.
\subsubsection{Proof of Lemma \ref{lemma:cones}}
Let $x_l\in\mathcal{C}_{\leq}(c,E^{-1}v_l,\psi_l,E)\backslash\{c\}$, $l=1,2$, and be otherwise arbitrary. Define then $z_l:=E(x_l-c)/\|E(x_l-c)\|\in\mathbb{S}^{n-1}$ for $l=1,2$. Hence, $z_l\in \mathcal{S}_l:=\mathcal{C}_{\leq}(0,v_l,\psi_l,I_n)\cap\mathbb{S}^{n-1}$, $l=1,2$. For $l=1,2$, $z_l$ satisfies, by \eqref{eq:def:cone}, $\cos(\psi_l)\leq z_l^\top v_l$, and consequently $\mathbf{d}_{\mathbb{S}^{n-1}}(v_l, z_l) \le \psi_l$. It follows from the triangle inequality that
$
\theta = \mathbf{d}_{\mathbb{S}^{n-1}} (v_1,v_2) \leq \mathbf{d}_{\mathbb{S}^{n-1}}(v_1,z_1)+\mathbf{d}_{\mathbb{S}^{n-1}}(z_1,z_2)+\mathbf{d}_{\mathbb{S}^{n-1}}(v_2,z_2)\leq\mathbf{d}_{\mathbb{S}^{n-1}}(z_1,z_2)+\psi_1+\psi_2.
$
Hence, in view of the condition $\psi_1+\psi_2<\theta$,  $\mathbf{d}_{\mathbb{S}^{n-1}}(z_1,z_2)>0$. This fact  implies that the compact sets $\mathcal{S}_1$ and $\mathcal{S}_2$ (and in turn $\mathcal{C}_{\leq}(c,E^{-1}v_l,\psi_l,E)\backslash\{c\}$, $l=1,2$) 
are disjoint.
%
\subsubsection{Proof of Lemma \ref{lemma:p-1}}
First, by~\eqref{eq:p_-1} and \eqref{eq:proj-refl-maps} we have
\begin{align}
	&p_{-1}^i-c_i=-E_i^{-1}\rho(E_ic_i)E_i(p_1^i-c_i) \text{, hence } \label{eq:difference}\\
   &\|E_i(p_{-1}^i-c_i)\|^2\overset{\eqref{eq:difference}}{=}(p_1^i-c_i)^\top E_i\rho(E_ic_i)\rho(E_ic_i)E_i(p_1^i-c_i) \nonumber\\
   &\hspace*{10pt}\overset{\eqref{eq:propLine3}}{=}(p_1^i-c_i)^\top E_i^2(p_1^i-c_i) = \|E_i(p_1^i-c_i)\|^2  \label{eq:equalnorm}
\end{align}
(so that $p^i_1 \neq c_i$ implies $p^i_{-1} \neq c_i$).
Based on~\eqref{eq:def:cone} for $\mathcal{C}(c_i,-c_i,\theta_i,E_i)$, one has
\begin{equation*}
\begin{aligned}
    & -c_i^\top E_i^{2}(p_{-1}^i-c_i) \overset{\eqref{eq:difference}}{=} c_i^\top E_i\rho(E_ic_i)E_i(p_1^i-c_i) \\
    & \quad\overset{\eqref{eq:propLine3}}{=}-c_i^{\top}  E_i^{2}(p_1^i-c_i)\! \overset{\eqref{eq:p_1}}{=}\!\cos(\theta_i)\|E_i(-c_i)\|\|E_i(p_1^i-c_i)\|\\
    &\quad \overset{\eqref{eq:equalnorm}}{=}\cos(\theta_i)\|E_i(-c_i)\|\|E_i(p_{-1}^i-c_i)\|.
\end{aligned}
\end{equation*}
This concludes by~\eqref{eq:def:cone} that $p_{-1}^i\in\mathcal{C}(c_i,-c_i,\theta_i,E_i)\backslash\{c_i\}$.

\subsubsection{Proof of Lemma \ref{lemma:equilibria}}
As for the $\Longrightarrow$ implication, let $x\in\mathbb{R}^n\backslash\{c\}$ be such that $\pi^\perp(E(x-c))E(x-p)=0$, which is equivalent to $\pi^\perp(E(x-c))E(p-c)=0$. By substituting the definition of $\pi^\perp(\cdot)$ in~\eqref{eq:proj-refl-maps}, one obtains $\|E(x-c)\|^2(p-c)=\big((p-c)^\top E^2(x-c)\big)(x-c)$. This very equation excludes that $(p-c)^\top E^2(x-c) = 0$ since $E$ is positive definite, $x\neq c$, and $p \neq c$ by assumption. So, letting $\lambda = \|E(x-c)\|^2/\big( (p-c)^\top E^2(x-c)\big)$ in~\eqref{eq:def:line}, one deduces that $x\in\mathcal{L}(c,p-c)$. The $\Longleftarrow$ implication is straightforward. 

\subsubsection{Proof of Lemma \ref{lemma:empty1}}
The quantities in~\eqref{eq:bar mu}-\eqref{eq:bar theta} are well-defined. Indeed, we have for~\eqref{eq:bar mu} that $1-4\underline{\delta}_i^2(1-\underline{\delta}_i^2/\delta_i^2)=(2\underline{\delta}_i^2-1)^2+4\underline{\delta}_i^4(\delta_i^{-2}-1)>0$ thanks to $\delta_i \in (\underline{\delta}_i,1)$. Moreover, by $\mu_i\in(1,\bar\mu_i(\delta_i))$, the argument of the $\arccos$ in~\eqref{eq:bar theta} belongs to $(0,1)$, so $\bar\theta_i(\delta_i,\mu_i)$ is also well-defined. Now, define
\begin{equation}
\label{eq:F-11iHat}
\begin{aligned}
\hat{\mathcal{F}}_m^i:=&\mathcal{E}_\leq(c_i,\delta_i E_i)\cap \mathcal{E}_\geq(\bar c_i,\mu_i\bar E_i)\cap\\ & \mathcal{C}_\geq(c_i,c_i-p_m^i,\psi_i,E_i)\cap\mathcal{E}_\geq(c_{i},E_{i}) \supset \mathcal{F}_m^i.
\end{aligned}    
\end{equation}
By proving that for each $i\in \mathbb{I}$ and $m \in \{ -1,1 \}$
\begin{equation}
\label{eq:the claim we actually prove}
\mathcal{L}(c_i,p_m^i-c_i)\cap\hat{\mathcal{F}}_m^i = \emptyset,
\end{equation}
the claim of the lemma is also proven. We prove then \eqref{eq:the claim we actually prove} for an arbitrary $i\in \mathbb{I}$ and an arbitrary $m \in \{ -1,1 \}$. For this proof, select the following angle $\psi_i^\prime$ as \emph{any} angle $\psi_i^\prime \in (0,\psi_i)$. First, let us show that the following set inclusions hold
\begin{subequations}
\label{eq:setincl}
\begin{align}
\label{eq:setincl1}
    &\mathcal{L}_\leq(c_i,p_m^i-c_i)\subset\mathcal{C}_\leq(c_i,c_i-p_m^i,\psi_i^\prime,E_i),\\
\label{eq:setincl2}
    &\mathcal{L}_\geq(c_i,p_m^i-c_i)\subset\mathcal{C}(c_i,-c_i,\theta_i,E_i).
\end{align}
\end{subequations}
Let $x\in\mathcal{L}_\leq(c_i,p_m^i-c_i)$. Then there exists $\lambda\leq 0$ such that $x-c_i=\lambda(p_m^i-c_i)$. Such $x-c_i$ verifies the condition $\cos(\psi_i^\prime) \|E_i (c_i-p_m)\| \|E_i(x-c_i)\| \le (c_i-p_m^i)^\top E^2 (x-c_i)$ corresponding to $\mathcal{C}_{\le}(c_i,c_i-p_m,\psi_i^\prime,E_i)$ by simple computations for any $0< \psi_i^\prime < \psi_i$ (since $\cos(\psi_i^\prime) \le 1$). This proves \eqref{eq:setincl1}. Now, let $x\in\mathcal{L}_\geq(c_i,p_m^i-c_i)$. Then there exists $\lambda\geq 0$ such that $x-c_i=\lambda(p_m^i-c_i)$. Such $x-c_i$ verifies the condition $\cos(\theta_i) \| E_i(-c_i)\| \|E_i(x-c_i)\| = -c_i^\top E_i^2 (x-c_i)$ corresponding to $\mathcal{C}(c_i,-c_i,\theta_i,E_i)$ by simple computations using that $\cos(\theta_i) \|E_i(-c_i)\| \| E_i(p_m^i-c_i)\| = - c_i^\top E_i^2 (p_m^i-c_i)$ (corresponding to $p_m^i\in\mathcal{C}(c_i,-c_i,\theta_i,E_i)$ from~\eqref{eq:p_1} and Lemma~\ref{lemma:p-1}). This proves \eqref{eq:setincl2}. Second, from \eqref{eq:setincl} one has
    \begin{align}
        \nonumber
        &\!\!\!\!\!\hat{\mathcal{F}}_m^i\cap\mathcal{L}(c_i,p_m^i-c_i)\\
        \nonumber
        &\!\!\!\!\!\!=(\hat{\mathcal{F}}_m^i\cap\mathcal{L}_\leq(c_i,p_m^i-c_i))\cup(\hat{\mathcal{F}}_m^i\cap\mathcal{L}_\geq(c_i,p_m^i-c_i))\\
        \label{eq:two intersections}
        &\!\!\!\!\!\!\!\subset\!\!(\hat{\mathcal{F}}_m^i\!\cap\!\mathcal{C}_\leq(c_i,c_i\!-\!p_m^i,\psi_i^\prime,E_i))\!\cup\!(\hat{\mathcal{F}}_m^i\!\cap\!\mathcal{C}(c_i,\!-c_i,\theta_i,E_i))
    \end{align}
and we prove that the two intersections in~\eqref{eq:two intersections} are empty. Since $0<\psi_i^\prime<\psi_i$, one obtains readily from the definition of the cone in~\eqref{eq:def:cone} that
\begin{equation*}
   \mathcal{C}_\geq(c_i,c_i-p_m^i,\psi_i,E_i)\cap \mathcal{C}_\leq(c_i,c_i-p_m^i,\psi_i^\prime,E_i)=\{c_i\}.    
\end{equation*}
This relationship and the definition of $\hat{\mathcal{F}}_m^i$ in~\eqref{eq:F-11iHat} imply that the first intersection in~\eqref{eq:two intersections} is empty. We show now that the second intersection in~\eqref{eq:two intersections} is also empty. Let $x\in\hat{\mathcal{F}}_m^i\cap\mathcal{C}(c_i,-c_i,\theta_i,E_i)$. So,
\begin{align}
    \nonumber
    &c_i^\top E_i^2(x-c_i)\overset{\eqref{eq:ci,Ei}}{=}-(x - c_i)^\top\! E_i^2(x - c_i)\\
    \nonumber
    &\qquad+\tfrac{1}{4}c_i^\top\! E_i^2 c_i\, (x\!-\!\bar c_i)^\top\! \bar E_i^2 (x\!-\!\bar c_i)-\tfrac{1}{4}c_i^\top\! E_i^2 c_i\\
    \nonumber
    &=-\|E_i(x-c_i)\|^2+\tfrac{1}{4}\|E_i c_i \|^2\|\bar E_i (x-\bar c_i)\|^2-\tfrac{1}{4}\|E_i c_i \|^2\\
    \label{eq:whenInH(epsilon_h,mu)}
    &\geq-\frac{1}{\delta_i^2}-\frac{1}{4}\|E_i c_i\|^2\Big(1-\frac{1}{\mu_i^2}\Big)
\end{align}
where the bound holds since $x\in\hat{\mathcal{F}}_m^i$ implies $x\in\mathcal{E}_\leq(c_i,\delta_i E_i)$ and $x\in\mathcal{E}_\geq(\bar c_i,\mu_i \bar E_i)$. We continue~\eqref{eq:whenInH(epsilon_h,mu)} as
\begin{align}
    &c_i^\top E_i^2(x-c_i)\geq-\frac{1}{\delta_i^2}-\frac{1}{4}\|E_i c_i\|^2\Big(1-\frac{1}{\mu_i^2}\Big) \nonumber\\
    &\overset{\eqref{eq:bar theta}, \eqref{eq:bar delta}}{=}-\cos(\bar\theta_i(\delta_i,\mu_i))/\underline{\delta}_i^2 \overset{\eqref{eq:bar delta}}{=}-\cos(\bar\theta_i(\delta_i,\mu_i))\|E_i c_i\| \nonumber\\
    &\geq - \cos(\bar\theta_i(\delta_i,\mu_i))\|E_i c_i \|\|E_i (x-c_i)\| \label{eq:cE(x-c) expr 1}
\end{align}
since $\cos(\bar \theta_i(\delta_i,\mu_i))\ge 0$ and $\|E_i(x-c_i)\| \ge 1$ ($x\in\hat{\mathcal{F}}_m^i$ implies $x\in\mathcal{E}_\geq(c_i,E_i)$). It is also $x\in\mathcal{C}(c_i,-c_i,\theta_i,E_i)$. So,
\begin{equation}
\label{eq:cE(x-c) expr 2}
\begin{aligned}
  c_i^\top E_i^2(x-c_i)&=-\cos(\theta_i)\|E_i c_i\|\|E_i (x-c_i)\|\\
  &<-\cos(\bar\theta_i(\delta_i,\mu_i))\|E_i c_i \|\|E_i (x-c_i)\|
\end{aligned}
\end{equation}
from the bound on $\theta_i$ in Table~\ref{table:parameters}. \eqref{eq:cE(x-c) expr 1} and \eqref{eq:cE(x-c) expr 2} contradict each other, so the second intersection in~\eqref{eq:two intersections} is \emph{also} empty. Then, \eqref{eq:the claim we actually prove} is proven.


\subsubsection{Proof of Lemma \ref{lemma:jump_map}}
Given~\eqref{eq:L_m1,L_m-1}-\eqref{eq:L_m0}, we just need to show that $\mathbf{L}(x,i,0) \neq \emptyset$ for all $(x,i) \in \mathcal{J}_0$. This holds if we show, as we do in the rest of the proof, that for each $x \in \real^n$ and $i \in \mathbb{I}$, $\mathbf{M}(x,i) \neq \emptyset$.  First, we show that $\cap_{m=-1,1}\mathcal{C}_{\leq}(c_i,c_i-p_m^i,\bar\psi_i,E_i)=\{ c_i \}$ for each $i \in \mathbb{I}$. To this end, note that $p_1^i\in\mathcal{C}(c_i,-c_i,\theta_i,E_i)$ by~\eqref{eq:p_1}, and this implies $\cos^2(\theta_i) \| E_i c_i \|^2 \| E_i(p^i_1-c_i) \|^2 = (c_i^\top E_i^2 (p^i_1-c_i))^2$ or, equivalently,
\begin{equation}
\label{eq:implied by p_1 on cone}
(p^i_1-c_i)^\top E_i \pi^{\theta_i}(E_i c_i) E_i (p^i_1-c_i)=0.
\end{equation}
Introduce then $v^i_m:=E_i(c_i-p_m^i)/\|E_i(c_i-p_m^i)\|$ for $m\in\{-1,1\}$, and compute
\begin{align*}
(v^i_1&)^\top v^i_{-1}
=\frac{(p_1^i-c_i)^\top E_i^2(p_{-1}^i-c_i)}{\|E_i(p_1^i-c_i)\|\|E_i(p_{-1}^i-c_i)\|}\\
&\overset{\eqref{eq:equalnorm}}{=}
-\frac{(p_1^i-c_i)^\top E_i \rho(E_ic_i) E_i (p_1^i-c_i)}{\|E_i(p_1^i-c_i)\|^2}\\
& \overset{\eqref{eq:piTheta:3}}{=} 
-\frac{(p_1^i-c_i)^\top E_i (2\pi^{\theta_i}(E_ic_i)-\cos(2\theta_i)I_n) E_i (p_1^i-c_i)}{\|E_i(p_1^i-c_i)\|^2}\\
& \overset{\eqref{eq:implied by p_1 on cone}}{=} \frac{\cos(2\theta_i) (p_1^i-c_i)^\top E_i^2 (p_1^i-c_i)}{\|E_i(p_1^i-c_i)\|^2}= \cos(2\theta_i)
\end{align*}
Then, by Lemma~\ref{lemma:cones} and $2 \bar\psi_i< 2\theta_i$, $\cap_{m=-1,1}\mathcal{C}_{\leq}(c_i,c_i-p_m^i,\bar\psi_i,E_i)=\{ c_i \}$. Second, note that
\begin{multline}
\label{eq:real^n minus c}
    \cup_{m=-1,1} \mathcal{C}_{>}(c_i,c_i-p_m^i,\bar\psi_i,E_i) =\\ \big( 
\cap_{m=-1,1} \mathcal{C}_{\leq}(c_i,c_i-p_m^i,\bar\psi_i,E_i)
\big)^\C= \real^n \backslash \{c_i\}.
\end{multline}
Therefore, we have 
\begin{equation}
\label{eq:real^n union c}
\cup_{m=-1,1} \mathcal{C}^i_m \overset{\eqref{eq:Cim}}{=}
\cup_{m=-1,1} \mathcal{C}_{\geq}(c_i,c_i-p_m^i,\bar\psi_i,E_i)=\mathbb{R}^n
\end{equation}
since $\cup_{m=-1,1} \mathcal{C}^i_m$ is a superset of the set in~\eqref{eq:real^n minus c} and contains $c_i$. So, for each $x \in \real^n$ and $i \in \mathbb{I}$, $\mathbf{M}(x,i) \neq \emptyset$ in~\eqref{eq:Mi}.  

\subsubsection{Proof of Lemma \ref{lemma:hbc}}

$\mathcal{F}$ and $\mathcal{J}$ are closed subsets of $\real^{n}\times\{-1,0,1\}$. $\mathbf{F}$ is continuous on $\mathcal F$.
$\mathbf{J}(x,i,m) \neq \emptyset$ for each $(x,i,m) \in \mathcal{J}$ thanks to Lemma~\ref{lemma:jump_map} and $\mathbf{J}$ has a closed graph relative to $\mathcal{J}$ because, in particular, the construction in~\eqref{eq:Mi} allows $\mathbf{M}$ to be set-valued whenever $x\in\cap_{m=-1,1} \mathcal{C}^i_m$. Then, $\mathbf{J}$ is outer semicontinuous and locally bounded relative to $\mathcal{J}$.

\subsubsection{Proof of Lemma \ref{lemma:union}}
If we prove that $\forall i \in \mathbb{I}, \,m\in\{-1,1\}$
\begin{equation}
\label{eq:key unions}
\Big( \bigcap\limits_{i^\prime\in\mathbb{I}}\mathcal{F}_0^{i^\prime} \Big) \cup \Big( \bigcup\limits_{i^\prime\in\mathbb{I}}\mathcal{J}_0^{i^\prime}  \Big)\! =\! \mathcal{W} \!= \! \mathcal{F}_m^i \cup \mathcal{J}_m^i
\end{equation}
then \eqref{eq:FJ0}, \eqref{eq:FJ1}, \eqref{eq:FJ-1} imply straightforwardly $\mathcal{F}_0 \cup \mathcal{J}_0 = \mathcal{W} \times \mathbb{I} = \mathcal{F}_1 \cup \mathcal{J}_1 = \mathcal{F}_{-1} \cup \mathcal{J}_{-1} $ and, in turn, \eqref{eq:FJ} implies $\mathcal{F} \cup \mathcal{J} = \mathcal{W} \times \mathbb{I} \times \mathbb{M} =: \mathcal{K}$. Therefore, we just need to prove \eqref{eq:key unions} in the remainder. For each $i \in \mathbb{I}$ and $m\in \{-1,1\}$,
\begin{multline}
    \mathcal{F}_m^i\cup \mathcal{J}_m^i\overset{\eqref{eq:F-11i},\eqref{eq:J^i_m}}{=}
\Big( \mathcal{E}_\leq(c_i,\delta_i E_i)\cap \mathcal{E}_\geq(\bar c_i,\mu_i\bar E_i)\\ 
    \cap \mathcal{C}_\geq(c_i,c_i-p_m^i,\psi_i,E_i)\cap \mathcal{W} \Big)\cup \Big( \big(\mathcal{E}_\geq(c_i,\delta_i E_i)\cup\mathcal{E}_\leq(\bar c_i,\mu_i\bar E_i)\\ 
    \cup\mathcal{C}_\leq(c_i,c_i\!-\!p_m^i,\psi_i,E_i)\big)\cap \mathcal{W} \Big) =\mathcal{W}.  
\end{multline}
We are left with proving $\big( \bigcap\limits_{i^\prime\in\mathbb{I}}\mathcal{F}_0^{i^\prime} \big) \cup \big( \bigcup\limits_{i^\prime\in\mathbb{I}}\mathcal{J}_0^{i^\prime}  \big)\! =\! \mathcal{W}$ in~\eqref{eq:key unions}. First, note that for each $i^\prime$, $\mathcal{H}_{i^\prime}\!(\epsilon_{i^\prime},\nu_{i^\prime})\! \cap\! \mathcal{W}=\mathcal{E}_\le(c_{i^\prime},\epsilon_{i^\prime}E_{i^\prime}) \cap \mathcal{E}_\ge(\bar c_{i^\prime},\nu_{i^\prime}\bar E_{i^\prime}) \cap \mathcal{W}$ and, hence,
\begin{equation*}
\begin{aligned}
& \big( \bigcap\limits_{i^\prime\in\mathbb{I}}\mathcal{F}_0^{i^\prime} \big) \cup \big( \bigcup\limits_{i^\prime\in\mathbb{I}}\mathcal{J}_0^{i^\prime}  \big) = \bigcap_{i^\prime\in \mathbb{I}}\Big(\big(\mathcal{E}_\ge(c_{i^\prime},\epsilon_{i^\prime}E_{i^\prime})  \\
&\quad\cup \mathcal{E}_\le(\bar c_{i^\prime},\nu_{i^\prime}\bar E_{i^\prime})\big)\cap \mathcal{W}\Big)\cup\bigcup_{i^\prime\in \mathbb{I}}\Big(\mathcal{H}_{i^\prime}(\epsilon_{i^\prime},\nu_{i^\prime}) \cap \mathcal{W}\Big) \\
&=
\Big(\!\bigcap_{i^\prime \in \mathbb{I}}\!\!\big(\mathcal{E}_\ge(c_{i^\prime},\epsilon_{i^\prime}E_{i^\prime}) \cup \mathcal{E}_\le(\bar c_{i^\prime},\nu_{i^\prime}\bar E_{i^\prime}) \big)
\\
&\quad\cup\bigcup_{i^\prime \in \mathbb{I}}\!\!\big(\mathcal{E}_\le(c_{i^\prime},\epsilon_{i^\prime}E_{i^\prime}) \!\cap\! \mathcal{E}_\ge(\bar c_{i^\prime},\nu_{i^\prime}\bar E_{i^\prime}) \big)\!\Big)\!\cap\! \mathcal{W}\!=\!\real^n \!\cap\! \mathcal{W}\!=\!\mathcal{W}.
\end{aligned}
\end{equation*}
From $\mathcal{F} \cup \mathcal{J}=\mathcal{K}$, the definition of $\mathcal{K}$ in~\eqref{eq:K} and $x^+=x$ in the jump map $\mathbf{J}$, it follows immediately that $\mathbf{J}(\mathcal{J}) \subset \mathcal{K}$.

\subsubsection{Proof of Lemma \ref{lemma:intersection}}
The intersection of $\mathcal{E}(c_i,\delta E_i)$ and $\mathcal{E}(\bar c_i,\mu\bar E_i)$ corresponds to the two quadratic equations
\begin{equation}
\label{eq:ellips intersection}
    \begin{cases}
    \delta^2\|E_i(x-c_i)\|^2=1\\
    \mu^2\|\bar E_i(x-\bar c_i)\|^2=1
    \end{cases}.
\end{equation}
By expanding squares and using \eqref{eq:ci,Ei}, \eqref{eq:ellips intersection} is equivalent to 
\begin{equation*}
    \begin{cases}
    \|E_ix\|^2-2c_i^\top E_i^2 x+\|E_i c_i\|^2-\delta^{-2}=0\\
    \|E_ix\|^2-c_i^\top E_i^2 x+\|E_i c_i\|^2(1-\mu^{-2})/4=0
    \end{cases}.
\end{equation*}
Solving for $\|E_ix\|^2$ and $c_i^\top E_i^2 x$, we obtain using~\eqref{eq:bar delta}
\begin{equation}
\label{eq:solving for cone quantities}
    \begin{cases}
    \|E_ix\|^2=\|E_ic_i\|^2\big((1+\mu^{-2})/2-\delta^{-2} \underline{\delta}_i^4\big)\\
    c_i^\top E_i^2 x=\|E_ic_i\|^2\big((3+\mu^{-2})/4-\delta^{-2} \underline{\delta}_i^4\big)
    \end{cases}
\end{equation}
and both right-hand sides of~\eqref{eq:solving for cone quantities} are positive because
\begin{subequations}
\begin{align}
& \label{eq:well posedness 1}
\begin{aligned}
     &(1+\mu^{-2})/2-\delta^{-2}\underline{\delta}_i^4\\
      &\quad \geq(1+\bar\mu_i(\delta)^{-2})/2-\delta^{-2}\underline{\delta}_i^4 \overset{\eqref{eq:bar mu}}{=}1-2\underline{\delta}_i^2+\delta^{-2}\underline{\delta}_i^4\\
      &\quad =(\underline{\delta}_i^2-1)^2+\underline{\delta}_i^4(\delta^{-2}-1) \geq(\underline{\delta}_i^2-1)^2>0
    \end{aligned}\\
& \label{eq:well posedness 2}
\begin{aligned}
    &(3+\mu^{-2})/4-\delta^{-2}\underline{\delta}_i^4\geq(3+\bar\mu_i(\delta)^{-2})/4-\delta^{-2}\underline{\delta}_i^4\\
    &\quad =1-\underline{\delta}_i^2>0, \text{ by Assumption~\ref{assumption:obstacle}}.
\end{aligned}
\end{align}
\end{subequations}
From~\eqref{eq:solving for cone quantities}, one obtains with some computations
\begin{equation}
\label{eq:ellips inters interm step}
\begin{aligned}
    & \frac{c_i^\top  E_i^2 x}{\|E_i c_i\|\|E_i x\|}
    \!\overset{\eqref{eq:solving for cone quantities}}{=}\!\frac{\|E_ic_i\|^2\big((3\!+\!\mu^{-2})/4\!-\!\delta^{-2} \underline{\delta}_i^4 \big)}{\|E_i c_i\|\|E_ic_i\|\sqrt{(1\!+\!\mu^{-2})/2\!-\!\delta^{-2}\underline{\delta}_i^4}}\\
    & \overset{\eqref{eq:bar theta}}{=}\!\!\!\frac{1-\cos(\bar\theta_i(\delta,\mu))\underline{\delta}_i^2}{\sqrt{(1+\mu^{-2})/2-\delta^{-2}\underline{\delta}_i^4}}.
\end{aligned}
\end{equation}
For $\mu \in [1, \big(1-4\underline{\delta}_i^2(1-\underline{\delta}_i^2/\delta^2)\big)^{-\frac{1}{2}}]$ and $\delta <1$, we can prove
\begin{equation}
\label{eq:well posedness 3}
\!\!\!\!\frac{1-\cos(\bar\theta_i(\delta,\mu))\underline{\delta}_i^2}{\sqrt{(1+\mu^{-2})/2-\delta^{-2}\underline{\delta}_i^4}}
    \!=\!\frac{(3+\mu^{-2})/4-\delta^{-2}\underline{\delta}_i^4}{\sqrt{(1+\mu^{-2})/2-\delta^{-2}\underline{\delta}_i^4}} \!<\!1
\end{equation}
(\textit{e.g.}, set $\chi:=(1+\mu^{-2})/2$, obtain the bounds of $\chi$ from the bounds of $\mu$, substitute $\chi$ in~\eqref{eq:well posedness 3}, and note that the obtained quadratic inequality holds true for such bounds of $\chi$ due to $\delta<1$).
Because of~\eqref{eq:well posedness 1}, \eqref{eq:well posedness 2} and \eqref{eq:well posedness 3}, the expression in~\eqref{eq:cosine expr} is well-defined and positive. Since \eqref{eq:ellips inters interm step} yields
\begin{equation*}
\begin{aligned}
    & \frac{c_i^\top  E_i^2 x}{\|E_i c_i\|\|E_i x\|}= \frac{1-\cos(\bar\theta_i(\delta,\mu))\underline{\delta}_i^2}{\sqrt{(1+\mu^{-2})/2-\delta^{-2}\underline{\delta}_i^4}} \overset{\eqref{eq:cosine expr}}{=}\cos(\vartheta_i(\delta,\mu)),
\end{aligned}
\end{equation*}
\eqref{eq:ellips inters contained in cone} holds as well by the cone definition in \eqref{eq:def:cone}.
\subsubsection{Proof of Lemma~\ref{lemma:disjoint}}
The obstacles $\{\mathcal{O}_i\}_{i\in\mathbb{I}}$ are sufficiently pairwise disjoint, so \eqref{eq:sufficient_disjoint} holds. Moreover, the sets $\mathcal{R}_{i'}(\delta_{i'},\mu_{i'})$ in~\eqref{eq:Ri_deltai_mui} and $\mathcal{E}_\leq(c_{i''},\delta_{i''} E_{i''})$ are bounded, and for $\delta_{i'},\delta_{i''}, \mu_{i'}\to 1$, $\mathcal{R}_{i'}(\delta_{i'},\mu_{i'})$ and $\mathcal{E}_\leq(c_{i''},\delta_{i''} E_{i''})$ reduce respectively to $\mathcal{R}^*_{i'}$ and $\mathcal{E}_\leq(c_{i''},E_{i''})$. By a continuity argument, there exist
parameters $\delta_i^*$ and $\mu_i^*$ such that the lemma holds.

\subsubsection{Proof of Lemma \ref{lemma:boundary}}
We prove the claim for arbitrary $i\in\mathbb{I}$ and $m \in \{-1,1\}$. 
Let us prove \eqref{eq:J0i_2}-\eqref{eq:Fmi}.
Thanks to the weak pairwise disjointness of $\{\mathcal{E}_\leq(c_i,\delta_iE_i)\}_{i\in\mathbb{I}}$ we have $\mathcal{E}_\leq(c_i,\delta_iE_i)\cap\mathcal{E}_\geq(c_{i^{\prime}},\delta_{i^\prime}E_{i^{\prime}})=\mathcal{E}_\leq(c_i,\delta_iE_i)$ for all $i^\prime\neq i$. Then,
$
    \mathcal{E}_\leq(c_i,\delta_iE_i)\cap\mathcal{W}=\mathcal{E}_\leq(c_i,\delta_iE_i)\cap\mathcal{E}_\geq(c_i,E_i),
$
which implies by~\eqref{eq:F-11i} that $\mathcal{F}_m^i$ satisfies \eqref{eq:Fmi}. By a similar argument, $\mathcal{J}_0^i$ satisfies \eqref{eq:J0i_2} as well.

Let us prove~\eqref{eq:boundary F_0 star}. Write the complement of $\mathcal{F}_0^*$ as
\begin{equation}
\label{eq:complement of F_0star}
\begin{aligned}
&\!\!\!\!\!(\mathcal{F}_0^*)^\C\overset{\eqref{eq:F0J0star}}{=}\big(\bigcap_{i\in\mathbb{I}}\mathcal{F}_0^i\big)^\C \overset{\eqref{eq:set identities}}{=}\bigcup_{i\in\mathbb{I}}\big(\mathcal{F}_0^i\big)^\C\\
&\!\!\!\!\!\overset{\eqref{eq:F0i}}{=}\bigcup_{i\in\mathbb{I}}\Big(\big(\mathcal{E}_\geq(c_i,\epsilon_i E_i)\cup\mathcal{E}_\leq(\bar c_i,\nu_i\bar E_i)\big)\cap\mathcal{W}\Big)^\C\\
&\!\!\!\!\!\overset{\eqref{eq:set identities}}{=}\bigcup_{i\in\mathbb{I}}\big(\mathcal{E}_\geq(c_i,\epsilon_i E_i)\cup\mathcal{E}_\leq(\bar c_i,\nu_i\bar E_i)\big)^\C\cup\mathcal{W}^\C\\
&\!\!\!\!\!\overset{\eqref{eq:W}}{=}\bigcup_{i\in\mathbb{I}}\mathcal{E}_<(c_i,\epsilon_iE_i)\cap\mathcal{E}_>(\bar c_i,\nu_i\bar E_i)\cup\bigcup_{i\in\mathbb{I}}\mathcal{E}_<(c_i,E_i)\\
&\!\!\!\!\!=\!\bigcup_{i\in\mathbb{I}}\big(\mathcal{E}_<(c_i,\epsilon_iE_i)\cap\mathcal{E}_>(\bar c_i,\nu_i\bar E_i)\big)\cup\mathcal{E}_<(c_i,E_i)\\
&\!\!\!\!\!=\!\bigcup_{i\in\mathbb{I}}(\mathcal{E}\!_<\!(c_i,\!\epsilon_iE_i)\!\cup\!\mathcal{E}\!_<\!(c_i\!,E_i))\!\cap\!(\mathcal{E}\!_>\!(\bar c_i,\nu_i\bar E_i)\!\cup\!\mathcal{E}\!_<\!(c_i,\!E_i))\\
&\!\!\!\!\!=\!\bigcup_{i\in\mathbb{I}}\mathcal{E}_<(c_i,\epsilon_iE_i)\cap(\mathcal{E}_>(\bar c_i,\nu_i\bar E_i)\cup\mathcal{E}_<(c_i,E_i))
    \end{aligned}
\end{equation}
because $\epsilon_i<1$. Thanks to the weak pairwise disjointness of $\{\mathcal{E}(c_i,\delta_iE_i)\}_{i\in\mathbb{I}}$ and $\delta_i<\epsilon_i$, the sets $\{\mathcal{E}_<(c_i,\epsilon_iE_i)\cap(\mathcal{E}_>(\bar c_i,\nu_i\bar E_i)\cup\mathcal{E}_<(c_i,E_i))\}_{i\in\mathbb{I}}$ can actually be proven to be pairwise separated. Then, we can use \eqref{eq:boundary of separated sets} to obtain the boundary of the set $\mathcal{F}_0^*$ as
\begin{align*}
& \partial\mathcal{F}_0^*=\partial\big((\mathcal{F}_0^*)^\C\big)\\
&\overset{\eqref{eq:complement of F_0star}}{=}\partial\Big(\bigcup_{i\in\mathbb{I}}\mathcal{E}_<(c_i,\epsilon_iE_i)\cap(\mathcal{E}_>(\bar c_i,\nu_i\bar E_i)\cup\mathcal{E}_<(c_i,E_i))\Big)\\
&\overset{\eqref{eq:boundary of separated sets}}{=}\bigcup_{i\in\mathbb{I}}\partial\left(\mathcal{E}_<(c_i,\epsilon_iE_i)\cap(\mathcal{E}_>(\bar c_i,\nu_i\bar E_i)\cup\mathcal{E}_<(c_i,E_i))\right)\\
&\overset{\eqref{eq:set identities}}{\subset}\bigcup_{i\in\mathbb{I}}\Big(\partial\mathcal{E}_<(c_i,\epsilon_iE_i)\cap\overline{\mathcal{E}_>(\bar c_i,\nu_i\bar E_i)\cup\mathcal{E}_<(c_i,E_i)}\\
&\hspace*{12pt}\cup\partial\big(\mathcal{E}_>(\bar c_i,\nu_i\bar E_i)\cup\mathcal{E}_<(c_i,E_i)\big)\cap\overline{\mathcal{E}_<(c_i,\epsilon_i E_i)}\Big)\\
&\overset{\eqref{eq:set identities}}{\subset}\bigcup_{i\in\mathbb{I}}\bigg(\mathcal{E}(c_i,\epsilon_iE_i)\cap\big(\mathcal{E}_\geq(\bar c_i,\nu_i\bar E_i)\cup\mathcal{E}_\leq(c_i,E_i)\big)\\
&\hspace*{12pt}\!\cup\!\Big(\!\big(\mathcal{E}(\bar c_i,\nu_i\bar E_i)\backslash \mathcal{E}_<(c_i,\!E_i)\big)\!\cup\!\big(\mathcal{E}(c_i,\!E_i)\backslash \mathcal{E}_>(\bar c_i,\nu_i\bar E_i)\big)\!\Big)\\
&\hspace*{24pt}\!\cap\!\mathcal{E}_\leq(c_i,\!\epsilon_i E_i)\!\bigg)\\
&=\bigcup_{i\in\mathbb{I}}\Big(\!\big(\mathcal{E}(c_i,\epsilon_iE_i)\!\cap\!\mathcal{E}_\geq(\bar c_i,\nu_i\bar E_i)\big)\!\cup\!\big((\mathcal{E}(\bar c_i,\nu_i\bar E_i)\backslash\mathcal{E}_<(c_i,E_i))\\
&\hspace*{12pt} \cap\mathcal{E}_\leq(c_i,\epsilon_i E_i)\big) \cup\big(\mathcal{E}(c_i,E_i)\backslash\mathcal{E}_>(\bar c_i,\nu_i\bar E_i)\big)\Big) =: \bigcup_{i\in\mathbb{I}}\mathcal{P}_0^i
\end{align*}
\eqref{eq:boundary F_0 star} is finally proven in~\eqref{eq:boundary1}, for which we note that: \textit{(a)} $\mathcal{J}_0^*$ is simplified into $\cup_{i\in\mathbb{I}}\mathcal{H}_i(\epsilon_i,\nu_i)$ thanks to~\eqref{eq:J0i_2}; \textit{(b)} $\mathcal{P}_0^i\cap\mathcal{H}_{i^\prime}(\epsilon_{i^\prime},1)=\emptyset$ for all $i\neq i^\prime$; \textit{(c)} most of the sets in the next-to-last expression are empty, so the last expression follows. 

Let us prove~\eqref{eq:boundary F^i_m}. First, note that by~\eqref{eq:Fmi} and \eqref{eq:safety:helmet}
\begin{equation*}
\begin{aligned}
&\mathcal{F}_m^i = \mathcal{E}_\leq(c_i,\delta_i E_i)\cap\mathcal{E}_\geq(c_i,E_i)\\
& \hspace*{2cm}\cap\mathcal{E}_\geq(\bar c_i,\mu_i\bar E_i)\cap\mathcal{C}_\geq(c_i,c_i-p_m^i,\psi_i,E_i),
\end{aligned}
\end{equation*}
which is an intersection of four closed sets. By successive applications of~\eqref{eq:set identities:l} and \eqref{eq:set identities:h}, the boundary of $\mathcal{F}_m^i$ satisfies
\begin{align}\label{eq:boundary:Fmi}
    \partial\mathcal{F}_m^i\subset\bigcup_{k\in\{1,2,3,4\}}\mathcal{P}_{m,k}^i
\end{align}
with the following definitions 
\begin{subequations}
\label{eq:definitions of auxiliary sets}
\begin{align}
    \nonumber
    &\mathcal{P}_{m,1}^i:=\mathcal{E}(c_i,\delta_i E_i)\cap\mathcal{E}_\geq(c_i,E_i)\cap\mathcal{E}_\geq(\bar c_i,\mu_i\bar E_i)\\
    &\hspace{4cm}\cap\mathcal{C}_\geq(c_i,c_i-p_m^i,\psi_i,E_i) \label{eq:definitions of auxiliary sets:1}\\
     \nonumber
    &\mathcal{P}_{m,2}^i:=\mathcal{E}_\leq(c_i,\delta_i E_i)\cap\mathcal{E}(c_i,E_i)\cap\mathcal{E}_\geq(\bar c_i,\mu_i\bar E_i)\\
    &\hspace{4cm}\cap\mathcal{C}_\geq(c_i,c_i-p_m^i,\psi_i,E_i) \label{eq:definitions of auxiliary sets:2}\\
     \nonumber
     &\mathcal{P}_{m,3}^i:=\mathcal{E}_\leq(c_i,\delta_i E_i)\cap\mathcal{E}_\geq(c_i,E_i)\cap\mathcal{E}(\bar c_i,\mu_i\bar E_i)\\
    &\hspace{4cm}\cap\mathcal{C}_\geq(c_i,c_i-p_m^i,\psi_i,E_i) \label{eq:definitions of auxiliary sets:3}\\
     \nonumber
     &\mathcal{P}_{m,4}^i:=\mathcal{E}_\leq(c_i,\delta_i E_i)\cap\mathcal{E}_\geq(c_i,E_i)\cap\mathcal{E}_\geq(\bar c_i,\mu_i\bar E_i)\\
    &\hspace{4cm}\cap\mathcal{C}(c_i,c_i-p_m^i,\psi_i,E_i). \label{eq:definitions of auxiliary sets:4}
\end{align}
\end{subequations}
Second, note that since $\mathcal{F}_m^i$ is closed, $\partial\mathcal{F}_m^i\subset\mathcal{F}_m^i\subset\mathcal{W}$, and hence $\partial\mathcal{F}_m^i\backslash\mathcal{W}\subset\mathcal{W}\backslash\mathcal{W}=\emptyset$. By this fact, we can write that 
\begin{align*}
&\partial\mathcal{F}_m^i\backslash\mathcal{J}_m^i\overset{\eqref{eq:J^i_m}}{=}\partial\mathcal{F}_m^i\backslash\bigg(\Big(\mathcal{E}_\geq(c_i,\delta_i E_i)\cup\mathcal{E}_\leq(\bar c_i,\mu_i\bar E_i)\\
&\qquad\cup\mathcal{C}_\leq(c_i,c_i\!-\!p_m^i,\psi_i,E_i)\Big)\cap\mathcal{W}\bigg)\\
&\overset{\eqref{eq:set identities}}{=}\partial\mathcal{F}_m^i\backslash\mathcal{W}\cup\partial\mathcal{F}_m^i\backslash\Big(\mathcal{E}_\geq(c_i,\delta_i E_i)\cup\mathcal{E}_\leq(\bar c_i,\mu_i\bar E_i)\\
&\qquad\cup\mathcal{C}_\leq(c_i,c_i\!-\!p_m^i,\psi_i,E_i)\Big)\\
&=\partial\mathcal{F}_m^i\backslash\Big(\!\mathcal{E}_\geq(c_i,\delta_i E_i)\!\cup\!\mathcal{E}_\leq(\bar c_i,\mu_i\bar E_i)\!\cup\!\mathcal{C}_\leq(c_i,c_i\!-\!p_m^i,\psi_i,E_i)\!\Big)\\
&\overset{\eqref{eq:boundary:Fmi}}{\subset}\Big(\cup_{k\in\{1,2,3,4\}}\mathcal{P}_{m,k}^i\Big)\backslash\Big(\mathcal{E}_\geq(c_i,\delta_i E_i)\cup\mathcal{E}_\leq(\bar c_i,\mu_i\bar E_i)\\
&\qquad\cup\mathcal{C}_\leq(c_i,c_i\!-\!p_m^i,\psi_i,E_i)\Big)\\
&\overset{\eqref{eq:set identities}}{=}\cup_{k\in\{1,2,3,4\}}\Big(\mathcal{P}_{m,k}^i\backslash\Big(\mathcal{E}_\geq(c_i,\delta_i E_i)\cup\mathcal{E}_\leq(\bar c_i,\mu_i\bar E_i)\\
&\qquad\cup\mathcal{C}_\leq(c_i,c_i\!-\!p_m^i,\psi_i,E_i)\Big)\Big)\\
&\overset{\eqref{eq:set identities}}{\subset}\cup_{k\in\{1,2,3,4\}}\Big(\big(\mathcal{P}_{m,k}^i\backslash\mathcal{E}_\geq(c_i,\delta_i E_i)\big)\\
&\qquad\cap\big(\mathcal{P}_{m,k}^i\backslash\mathcal{E}_\leq(\bar c_i,\mu_i\bar E_i)\big)\!\cap\!\big(\mathcal{P}_{m,k}^i\backslash\mathcal{C}_\leq(c_i,c_i\!-\!p_m^i,\psi_i,E_i)\big)\Big).
\end{align*}
Finally, we simplify this expression through the following facts
\begin{equation}
\label{eq:empty set differences}
\begin{aligned}
    &\mathcal{P}_{m,1}^i\backslash\mathcal{E}_\geq(c_i,\delta_i E_i)=\emptyset, \quad\mathcal{P}_{m,3}^i\backslash\mathcal{E}_\leq(\bar c_i,\mu_i\bar E_i)=\emptyset\\
    &\mathcal{P}_{m,4}^i\backslash\mathcal{C}_\leq(c_i,c_i\!-\!p_m^i,\psi_i,E_i)=\emptyset,
\end{aligned}
\end{equation}
which are an immediate consequence of~\eqref{eq:definitions of auxiliary sets} and yield
\begin{equation*}
\begin{aligned}
&\partial\mathcal{F}_m^i\backslash\mathcal{J}_m^i\overset{\eqref{eq:empty set differences}}{\subset}\mathcal{P}_{m,2}^i\backslash\Big(\mathcal{E}_\geq(c_i,\delta_i E_i)\cup\mathcal{E}_\leq(\bar c_i,\mu_i\bar E_i)\\
&\quad\cup\mathcal{C}_\leq(c_i,c_i\!-\!p_m^i,\psi_i,E_i)\Big)\\
&\overset{\eqref{eq:set identities}}{=}\!\mathcal{P}_{m,2}^i\!\cap\!\Big(\!\mathcal{E}_\geq(c_i,\!\delta_i E_i)\!\cup\!\mathcal{E}_\leq(\bar c_i,\!\mu_i\bar E_i)\!\cup\!\mathcal{C}_\leq(c_i,\!c_i\!-\!p_m^i,\!\psi_i,\!E_i)\!\Big)^\C\\
&\overset{\eqref{eq:set identities}}{=}\!\mathcal{P}_{m,2}^i\cap\mathcal{E}_<(c_i,\!\delta_i E_i)\cap\mathcal{E}_>(\bar c_i,\!\mu_i\bar E_i)\cap\mathcal{C}_>(c_i,\!c_i\!-\!p_m^i,\!\psi_i,\!E_i)\\
&=\!\mathcal{E}_<(c_i,\!\delta_i E_i)\!\cap\!\mathcal{E}(c_i,\!E_i)\!\cap\!\mathcal{E}_>(\bar c_i,\!\mu_i\bar E_i)\!\cap\!\mathcal{C}_>(c_i,c_i-p_m^i,\!\psi_i,\!E_i)\\
&=\!\mathcal{E}(c_i,E_i)\cap\mathcal{E}_>(\bar c_i,\mu_i\bar E_i)\cap\mathcal{C}_>(c_i,c_i-p_m^i,\psi_i,E_i)\\
&\overset{\eqref{eq:set identities}}{=}\!\mathcal{E}(c_i,E_i)\backslash(\mathcal{E}_\leq(\bar c_i,\mu_i\bar E_i)\cup\mathcal{C}_\leq(c_i,c_i-p_m^i,\psi_i,E_i)).
\end{aligned}
\end{equation*}

\begin{table*}
\begin{equation}\label{eq:boundary1}
\begin{aligned}
& \partial\mathcal{F}_0^*\backslash\mathcal{J}_0^* \subset\Big(\bigcup_{i\in\mathbb{I}}\mathcal{P}_0^i\Big)\backslash\Big(\bigcup_{i\in\mathbb{I}}\mathcal{H}_i(\epsilon_i,\nu_i)\Big)
=\bigcup_{i\in\mathbb{I}}\!\Big(\!\mathcal{P}_0^i\backslash\!\bigcup_{i^\prime\in\mathbb{I}}\!\big(\mathcal{H}_{i^\prime}(\epsilon_{i^\prime},\nu_{i^\prime})\big)\!\Big)
= \bigcup_{i\in\mathbb{I}}\!\Big(\!\bigcap_{i^\prime\in \mathbb{I}}\! \big(\mathcal{P}_0^i\backslash\mathcal{H}_{i^\prime}(\epsilon_{i^\prime},\nu_{i^\prime})\big)\!\Big)
= \bigcup_{i\in\mathbb{I}}\!\Big(\!\mathcal{P}_0^i\backslash\mathcal{H}_{i}(\epsilon_{i},\nu_i)\cap \!\bigcap_{i^\prime\in \mathbb{I},i^\prime \neq i} \! \mathcal{P}_0^i\!\Big)\\
&\!\!\!=\!\bigcup_{i\in\mathbb{I}}\big(\mathcal{P}_0^i\backslash\mathcal{H}_i(\epsilon_i,\nu_i)\big)=\bigcup_{i\in\mathbb{I}}\Big(\!\big(\mathcal{E}(c_i,\epsilon_iE_i)\cap\mathcal{E}_\geq(\bar c_i,\nu_i\bar E_i)\big)\!\cup\!\big((\mathcal{E}(\bar c_i,\nu_i\bar E_i)\backslash\mathcal{E}_<(c_i,E_i))\cap\mathcal{E}_\leq(c_i,\epsilon_i E_i)\big)\!\cup\!(\mathcal{E}(c_i,E_i)\backslash\mathcal{E}_>(\bar c_i,\nu_i\bar E_i))\!\Big)\backslash\mathcal{H}_i(\epsilon_i,\nu_i)\\ 
&\!\!\!=\!\bigcup_{i\in\mathbb{I}}\Big(\big(\mathcal{E}(c_i,\epsilon_iE_i)\cap\mathcal{E}_\geq(\bar c_i,\nu_i\bar E_i)\big)\cup\big((\mathcal{E}(\bar c_i,\nu_i\bar E_i)\backslash\mathcal{E}_<(c_i,E_i))\cap\mathcal{E}_\leq(c_i,\epsilon_i E_i)\big)\cup(\mathcal{E}(c_i,E_i)\backslash\mathcal{E}_>(\bar c_i,\nu_i\bar E_i))\Big)\cap\mathcal{H}_i(\epsilon_i,\nu_i)^\C\\
&\!\!\!=\!\bigcup_{i\in\mathbb{I}}\!\Big(\!\big(\mathcal{E}(c_i,\!\epsilon_iE_i\!)\!\cap\!\mathcal{E}_\geq(\bar c_i,\!\nu_i\bar E_i\!)\big)\!\cup\!\big((\mathcal{E}(\bar c_i,\!\nu_i\bar E_i\!)\backslash\mathcal{E}_<(c_i,\!E_i\!))\!\cap\!\mathcal{E}_\leq(c_i,\!\epsilon_i E_i\!)\big)\!\cup\!(\mathcal{E}(c_i,\!E_i\!)\backslash\mathcal{E}_>(\bar c_i,\!\nu_i\bar E_i\!))\!\Big)\!\cap\!\big(\mathcal{E}_>(c_i,\!\epsilon_iE_i\!)\!\cup\!\mathcal{E}_<(\bar c_i,\!\nu_i\bar E_i\!)\!\cup\!\mathcal{E}_<(c_i,\!E_i\!)\big)\\
&\!\!\!=\!\bigcup_{i\in\mathbb{I}}\!\Big(\!\mathcal{E}(c_i,\epsilon_iE_i)\cap\mathcal{E}_\geq(\bar c_i,\nu_i\bar E_i)\cap\mathcal{E}_>(c_i,\epsilon_iE_i)\Big)\!\cup\!\Big(\!\mathcal{E}(c_i,\epsilon_iE_i)\cap\mathcal{E}_\geq(\bar c_i,\nu_i\bar E_i)\cap\mathcal{E}_<(\bar c_i,\nu_i\bar E_i)\Big)\!\cup\!\Big(\mathcal{E}(c_i,\epsilon_iE_i)\cap\mathcal{E}_\geq(\bar c_i,\nu_i\bar E_i)\cap\mathcal{E}_<(c_i,E_i)\Big)\\
&\qquad\cup\Big(\big(\mathcal{E}(\bar c_i,\nu_i\bar E_i)\backslash\mathcal{E}_<(c_i,E_i)\big)\cap\mathcal{E}_\leq(c_i,\epsilon_i E_i)\cap\mathcal{E}_>(c_i,\epsilon_iE_i)\Big)\cup\Big(\big(\mathcal{E}(\bar c_i,\nu_i\bar E_i)\backslash\mathcal{E}_<(c_i,E_i)\big)\cap\mathcal{E}_\leq(c_i,\epsilon_i E_i)\cap\mathcal{E}_<(\bar c_i,\nu_i\bar E_i)\Big)\\
&\qquad\cup\Big(\big(\mathcal{E}(\bar c_i,\nu_i\bar E_i)\backslash\mathcal{E}_<(c_i,E_i)\big)\cap\mathcal{E}_\leq(c_i,\epsilon_i E_i)\cap\mathcal{E}_<(c_i,E_i)\Big)\cup\Big(\big(\mathcal{E}(c_i,E_i)\backslash\mathcal{E}_>(\bar c_i,\nu_i\bar E_i)\big)\cap\mathcal{E}_>(c_i,\epsilon_iE_i)\Big)\\
&\qquad\cup\Big(\big(\mathcal{E}(c_i,E_i)\backslash\mathcal{E}_>(\bar c_i,\nu_i\bar E_i)\big)\cap\mathcal{E}_<(\bar c_i,\nu_i\bar E_i)\Big)\cup\Big(\big(\mathcal{E}(c_i,E_i)\backslash\mathcal{E}_>(\bar c_i,\nu_i\bar E_i)\big)\cap\mathcal{E}_<(c_i,E_i)\Big) = \bigcup_{i\in\mathbb{I}}\mathcal{E}(c_i,E_i)\backslash\mathcal{E}_\geq(\bar c_i,\nu_i\bar E_i).
\end{aligned}
\end{equation}
\end{table*}

\subsubsection{Proof of Lemma \ref{lemma:decrease}}
We divide the proof into steps.
\textit{Step 1: For each $i \in \mathbb{I}$ and $m \in \{-1,1\}$, $\delta \in [\delta_i,1]$ and $\psi \in [\psi_i, \bar \psi_i)$, consider the sets
\begin{subequations}
\begin{align}
& \label{eq:tilde F^i_m3}
\begin{aligned}
& \tilde{\mathcal{F}}^i_m(\delta,\psi)\!:=\mathcal{E}_\leq(c_i,\delta E_i)\cap\mathcal{E}_\geq(c_i,E_i)\cap\mathcal{E}_\geq(\bar c_i,\mu_i\bar E_i)\\
& \qquad\cap\mathcal{C}_\geq(c_i,c_i-p_m^i,\psi,E_i).
\end{aligned}\\
& \label{eq:tilde P^i_m3}
\begin{aligned}
& \tilde{\mathcal{P}}^i_{m,3}(\delta,\psi) := \mathcal{E}_\leq (c_i,\delta E_i)\cap\mathcal{E}_\geq(c_i,E_i)\\
& \qquad\cap\mathcal{E}(\bar c_i,\mu_i\bar E_i)\cap\mathcal{C}_\geq(c_i,c_i-p_m^i,\psi,E_i)
\end{aligned}
\end{align}
\end{subequations}
Each maximal solution $x$ to the flow-only hybrid system
\begin{equation}
\label{eq:flow-only avoidance hybrid system}
\dot x = \kappa(x,i,m) =: \mathfrak{u}(x), \quad x \in \tilde{\mathcal{F}}^i_m(\delta,\psi) 
\end{equation}
has 
$T={\sup}_t \dom x < +\infty$ and $x(T) \in \tilde{\mathcal{P}}^i_{m,3}(\delta,\psi)$.}
We first prove that $T$ is finite. Consider the following nonnegative function $\mathbf{V}(x):=\tfrac{1}{2}\|E_i(x-p_m^i)\|^2$. Simple computations, \eqref{eq:u}, and \eqref{eq:propLine1} yield that for all $x \in \tilde{\mathcal{F}}_m^i(\delta,\psi)$
\begin{equation*}
\begin{aligned}
    \langle & \nabla \mathbf{V} (x), \kappa(x,i,m) \rangle \\
    & \! =\!-k_m (x-p_m^i)^\top E_i \pi^\perp(E_i(x-c_i)) E_i(x-p_m^i) \\
    & \! =\!-k_m \|\pi^\perp(E_i(x-c_i)) E_i(x-p_m^i)\| \!<\!-\mathfrak{e} \!<\! 0,
\end{aligned} 
\end{equation*}
where $\mathfrak{e}>0$ follows from $\pi^\perp(E_i(x-c_i)) E_i(x-p_m^i)$ vanishing only for $x \in \mathcal{L}(c_i,p_m^i-c_i)$ (Lemma~\ref{lemma:equilibria}), and $\mathcal{L}(c_i,p_m^i-c_i)$ is separated by a positive distance from $\tilde{\mathcal{F}}_m^i(\delta,\psi)$ (Lemma~\ref{lemma:empty1}).
Then, $T$ is finite, otherwise $\mathbf{V}$ evaluated along solutions would become negative. In order to show $x(T) \in \tilde{\mathcal{P}}^i_{m,3}(\delta,\psi)$, we resort to a viability argument based on tangent cones.
To this end, we need the next lemma.
\begin{lemma}
\label{lemma:tangent cones single subsets}
For all $x \in \tilde{\mathcal{F}}_m^i(\delta,\psi)\backslash \tilde{\mathcal{P}}^i_{m,3}(\delta,\psi)$,
\begin{subequations}
\begin{align}
\mathfrak{u}(x) & \in  \mathbf{T}_{\mathcal{E}_\leq(c_i,\delta E_i)}(x) \label{eq:tangent cone: ell delta}\\
\mathfrak{u}(x) & \in  \mathbf{T}_{\mathcal{E}_\geq(c_i,E_i)}(x)\label{eq:tangent cone: ell}\\
\mathfrak{u}(x) & \in  \mathbf{T}_{\mathcal{E}_\geq(\bar c_i,\mu_i\bar E_i)}(x)\label{eq:tangent cone: bar ell}\\
\mathfrak{u}(x) & \in  \mathbf{T}_{\mathcal{C}_\geq(c_i,c_i-p_m^i,\psi,E_i)}(x).\label{eq:tangent cone: cone}
\end{align}
\end{subequations}
\end{lemma}
$\tilde{\mathcal{F}}^i_m(\delta,\psi)$ in~\eqref{eq:tilde F^i_m3} is the intersection of four closed sets: if 
\begin{equation}
\label{eq:vector field in tangent cone of each subset}
\begin{aligned}
\mathfrak{u}(x) \in  & \mathbf{T}_{\mathcal{E}_\leq(c_i,\delta E_i)}(x) \cap \mathbf{T}_{\mathcal{E}_\geq(c_i,E_i)}(x) \\ 
& \cap  \mathbf{T}_{\mathcal{E}_\geq(\bar c_i,\mu_i\bar E_i)}(x) \cap \mathbf{T}_{\mathcal{C}_\geq(c_i,c_i-p_m^i,\psi,E_i)}(x),
\end{aligned}
\end{equation}
then $\mathfrak{u}(x) \in \mathbf{T}_{\tilde{\mathcal{F}}^i_m(\delta,\psi)}(x)$ by the next fact, which is an immediate corollary of \cite[Thm.~5]{rockafellar1979clarke}.
\begin{fact}[{\cite[Thm.~5]{rockafellar1979clarke}}]
{\itshape
Let $v \in \mathcal{A} \cap \mathcal{B}$ with $\mathcal{A}$, $\mathcal{B}$ closed subsets of $\real^n$. Suppose $\mathbf{T}_\mathcal{A}(v) \cap ( \mathbf{T}_\mathcal{B}(v)^\circ) \neq \emptyset$. Then, $\mathbf{T}_{\mathcal{A} \cap \mathcal{B}}(v) \supset \mathbf{T}_\mathcal{A}(v) \cap \mathbf{T}_\mathcal{B}(v)$.}
\end{fact}
The condition~\eqref{eq:vector field in tangent cone of each subset} has been checked in Lemma~\ref{lemma:tangent cones single subsets} for each $x \in \tilde{\mathcal{F}}_m^i(\delta,\psi)\backslash \tilde{\mathcal{P}}^i_{m,3}(\delta,\psi)$, hence 
\begin{equation}
\mathfrak{u}(x) \in \mathbf{T}_{\tilde{\mathcal{F}}^i_m(\delta,\psi)}(x) \quad \forall x \in \tilde{\mathcal{F}}_m^i(\delta,\psi)\backslash \tilde{\mathcal{P}}^i_{m,3}(\delta,\psi).
\end{equation}
Then, it can only be $x(T) \in \tilde{\mathcal{P}}^i_{m,3}(\delta,\psi)$, otherwise the solution could be further extended by viability results such as~\cite[Lemma~5.26(b)]{goebel2012hybrid}.

\textit{Step 2: For each $i \in \mathbb{I}$, both the sets $\mathcal{E}_{\leq}(\bar c_i,\bar E_i)$ and $\mathcal{C}(0,c_i,\vartheta_i(\delta,\mu_i),E_i)$ ($\delta \in [\delta_i,1]$) are forward invariant under the vector field $-k_0x$.}
For $x\in\mathcal{E}(\bar c_i,\bar E_i)$, we have 
\begin{equation}
\begin{aligned}
&    -k_0x^\top\bar E_i^2(x-\bar c_i)\!\overset{\eqref{eq:ci,Ei}}{=}\!-\tfrac{k_0}{2}\|\bar E_i x\|^2\!-\!\tfrac{k_0}{2}x^\top \bar E_i^2(x-c_i)\\
&\overset{\eqref{eq:dist_to_c_decreases}}{=}-\tfrac{k_0}{2}\|\bar E_i x\|^2 +\tfrac{k_0}{2} \big(1- \|\bar E_i(x-\bar c_i)\|^2\big)\\
&=-\tfrac{k_0}{2}\|\bar E_i x\|^2\leq 0,
\end{aligned}
\end{equation}
where the last equality follows from $x\in\mathcal{E}(\bar c_i,\bar E_i)$.
Therefore, $\{-k_0x\}\in\mathcal{P}_\leq(0,\bar E_i^2(x-\bar c_i))=\mathbf{T}_{\mathcal{E}(\bar c_i,\bar E_i)}(x)$. 
For $x\in\mathcal{C}(0,c_i,\vartheta_i(\delta,\mu_i),E_i)$, we have
\begin{equation}
    \begin{aligned}
      -k_0x^\top E_i\pi^{\vartheta_i(\delta,\mu_i)}(E_ic_i)E_ix=0.
    \end{aligned}
\end{equation}
Therefore, $\{-k_0x\}\in \mathcal{P}(0,E_i \pi^{\vartheta_i(\delta,\mu_i)}(E_i c_i)E_ix)=\mathbf{T}_{\mathcal{C}(0,c_i,\vartheta_i(\delta,\mu_i),E_i)}(x)$ (cf.~\eqref{eq:normal vector to cone surface}). Forward invariance follows then from the classical Nagumo's theorem.

\textit{Step 3: Proof of~\eqref{eq:strict decrease}.}
Let $\xi(t,j)=(x(t,j),i(t,j),m(t,j))\in\mathcal{F}_l\times\{l\}$ with $l\in\{-1,1\}$ and $(t,j)\in\dom\xi$. Hence, by~\eqref{eq:ri},
\begin{equation}
\label{eq:norm of x(t,j)}
\| x(t,j) \| \ge r_{i(t,j)}>0
\end{equation}
thanks to the discussion below~\eqref{eq:ri}. We further divide into mutually exclusive subcases.

\textit{Step 3a: $\xi(t,j)\in(\mathcal{F}_l\cap\mathcal{J}_l)\times\{l\}$ and the solution jumps.}
By~\eqref{eq:jump} and \eqref{eq:L_m1,L_m-1}, $\xi(t,j+1)=(x(t,j),i(t,j),0)$. Depending on $x(t,j)$, either the solution never jumps again or  reaches the set $\mathcal{J}_0^*$ in~\eqref{eq:F0J0star} at some $(t^\prime,j+1)$. 
Consider the former case. The flow map for $x$ in~\eqref{eq:flow} ensures $x(\tau,j+1)=\exp(-k_0(\tau-t)) x(t,j) $ for all $\tau \ge t$. By~\eqref{eq:norm of x(t,j)}, $s \ge t$ exists such that $\| x(s,j+1) \| = \min_{i^\prime \in \mathbb{I}} r_{i^\prime} /2 \!=:\! \sigma_0 \!>\!0$. So, 
\begin{equation}
\label{eq:decrease sigma 0}
\| x(s,j+1) \| + \sigma_0 = \min_{i^\prime \in \mathbb{I}} r_{i^\prime} \le r_{i(t,j)} \le \| x(t,j) \|
\end{equation}
and the claim of the proposition is proven.
Consider the latter case, i.e., there exist $t^\prime\geq t$ and $i^\prime \in \mathbb{I}$ such that $x(t^\prime,j+1)\in\mathcal{J}_0^{i^\prime} \subset \mathcal{J}_0^*$. $\mathcal{J}_0^{i^\prime}$ is a shrinking of the set $\mathcal{F}^{i^\prime}_1 \cup \mathcal{F}^{i^\prime}_{-1}$ by construction of the flow and jump sets (cf.~\eqref{eq:J0i} and \eqref{eq:F-11i}, where the union of the cones from~\eqref{eq:F-11i} gives $\real^n$ from the same arguments yielding~\eqref{eq:real^n union c}). Then, by Lemma~\ref{lemma:jump_map},  there exist $l^\prime$ such that $\xi(t^\prime, j+2) \in (\mathcal{F}_{l^\prime} \backslash \mathcal{J}_{l^\prime}) \times \{ l^\prime \}$ and the solution is forced to flow as in Step~3b below.

\textit{Step 3b: $\xi(t,j)\in(\mathcal{F}_l \backslash \mathcal{J}_l)\times\{l\}$ and the solution flows.}
If $x(t,j) \in \mathcal{F}^i_l \backslash \mathcal{J}^i_l$, there exists $\delta \in (\delta_i,1]$ and $\psi \in (\psi_i, \bar \psi_i)$ such that $x(t,j) \in \tilde{\mathcal{F}}^i_l(\delta,\psi)$ as defined in~\eqref{eq:tilde F^i_m3}. From the facts
\begin{subequations}
\label{eq:auxiliary for intersection of tilde F^i_l and J^i_l}
\begin{align}
&    \mathcal{H}_i(\delta,\mu_i)\cap\mathcal{E}_\geq(c_i,\delta_iE_i)=\emptyset\\
&    \begin{aligned}
\mathcal{C}_\geq(c_i,c_i-p_l^i,\psi,E_i)\cap\mathcal{C}_\leq(c_i,&c_i-p_l^i,\psi_i,E_i)\\
& =\{ c_i \} \notin \mathcal{E}_\ge(c_i,E_i),
\end{aligned}
\end{align}
\end{subequations}
$\tilde{\mathcal{F}}_l^i(\delta,\psi)\cap\mathcal{J}_l^i$ is a subset of $\tilde{\mathcal{P}}_{l,3}^i(\delta,\psi)$ because
\begin{align}
    &\tilde{\mathcal{F}}_l^i(\delta,\psi)\cap\mathcal{J}_l^i =\mathcal{E}_\leq(c_i,\delta E_i) \cap \mathcal{E}_\geq(c_i,E_i) \nonumber\\
    & \hspace*{.4cm} \cap \mathcal{E}_\geq(\bar c_i,\mu_i\bar E_i) \cap \mathcal{C}_\geq(c_i,c_i-p_l^i,\psi,E_i) \cap \mathcal{J}_l^i \nonumber\\
    &\subset\mathcal{H}_i(\delta,\mu_i)\cap\mathcal{C}_\geq(c_i,c_i-p_l^i,\psi,E_i) \nonumber\\
    &\hspace*{.4cm}\cap\big(\mathcal{E}_\geq(c_i,\delta_iE_i)\cup\mathcal{E}_\leq(\bar c_i,\mu_i\bar E_i)\cup\mathcal{C}_\leq(c_i,c_i-p_l^i,\psi_i,E_i)\big) \nonumber\\
    &=\big(\mathcal{H}_i(\delta,\mu_i)\cap\mathcal{C}_\geq(c_i,c_i-p_l^i,\psi,E_i)\cap\mathcal{E}_\geq(c_i,\delta_iE_i)\big) \nonumber\\
    &\hspace*{.4cm}\cup\big(\mathcal{H}_i(\delta,\mu_i)\cap\mathcal{C}_\geq(c_i,c_i-p_l^i,\psi,E_i)\cap\mathcal{E}_\leq(\bar c_i,\mu_i\bar E_i)\big) \nonumber\\
    & \hspace*{.4cm}\cup\!\big(\mathcal{H}_i(\delta,\mu_i)\!
    \cap\!\mathcal{C}_\geq(c_i,c_i\!-\!p_l^i,\psi,E_i)\!\cap\!\mathcal{C}_\leq(c_i,c_i\!-\!p_l^i,\psi_i,E_i))\big) 
    \nonumber\\
    &\overset{\eqref{eq:auxiliary for intersection of tilde F^i_l and J^i_l}}{=}\mathcal{H}_i(\delta,\mu_i)\cap\mathcal{C}_\geq(c_i,c_i-p_l^i,\psi,E_i)\cap\mathcal{E}_\leq(\bar c_i,\mu_i\bar E_i) \nonumber\\
    &=\mathcal{E}_\leq(c_i,\delta E_i)\!\cap\!\mathcal{E}_\geq(c_i,E_i)\!\cap\!\mathcal{E}(\bar c_i,\mu_i\bar E_i)\!\cap\!\mathcal{C}_\geq(c_i,c_i\!-\!p^i_l,\psi,E_i)\nonumber\\
	&\overset{\eqref{eq:tilde P^i_m3}}{=}\tilde{\mathcal{P}}_{l,3}^i(\delta,\psi). \label{eq:intersection of tilde F^i_l and J^i_l}
\end{align}
This fact combined with Step~1, shows that the solution leaves the set $\mathcal{F}_l^i$ in finite time through the set $\tilde{\mathcal{P}}_{l,3}^i(\delta,\psi)$, where it jumps.
Then, we have $x(t^\prime,j+1)=x(t^\prime,j) \in \tilde{\mathcal{P}}_{m(t^\prime,j),3}^{i(t^\prime,j)}(\delta,\psi)$, $i(t^\prime,j+1)= i(t^\prime,j)=:\iota$ and $m(t^\prime,j)=0$.
Then, by~\eqref{eq:intersection of tilde F^i_l and J^i_l}, $x(t^\prime,j+1)\in\mathcal{E}(c_\iota,\delta E_\iota)\cap\mathcal{E}(\bar c_\iota,\mu_\iota\bar E_\iota)$, and, by Lemma~\ref{lemma:intersection}, $x(t^\prime, j+1)\in\mathcal{C}(0,c_\iota,\vartheta_\iota(\delta,\mu_\iota),E_\iota)$.
We have shown in Step~2 that both the sets $\mathcal{E}_{\leq}(\bar c_\iota,\bar E_\iota)$ and $\mathcal{C}(0,c_\iota,\vartheta_\iota(\delta,\mu_\iota),E_\iota)$ are forward invariant under the stabilization flow map for $x$, i.e., $-k_0x$. Since the obstacles are weakly disjoint, the solution then flows in $\mathcal{E}_\le(\bar c_\iota,\bar E_\iota) \cap \mathcal{C}(0,c_\iota,\vartheta_\iota(\delta,\mu_\iota),E_\iota)$ until it reaches the set $\mathcal{E}(c_\iota,\delta_\iota E_\iota)$ at $(t'',j+1)$. We either have $\|x(t'',j+1)\|<r_\iota$ or $\|x(t'',j+1)\|\geq r_\iota$. Consider the former case. Define 
\begin{equation}\label{eq:sigma1}
    \sigma_1\!:=\!\!\!\min_{{i,i^\prime\in\mathbb{I},\,i \neq i^\prime}} \!\!\! \mathbf{dist}(\mathcal{E}_\leq(c_i,\delta_i E_i),\mathcal{E}_\leq(c_{i^\prime},\delta_{i^\prime} E_{i^\prime}))>0,
\end{equation}
which is positive because obstacle are compact, pairwise disjoint sets. Since $x(t'',j+1)\in\mathcal{E}(c_\iota,\delta_\iota E_\iota)$ and the obstacles are weakly pairwise disjoint, the solution can only flow up to the time $(s,j+1)$ such that $\sigma_1$ is traversed, i.e.,
\begin{equation}
\label{eq:decrease sigma 1}
\| x(s,j+1) \| + \sigma_1 \le \| x(t'',j+1) \| < r_\iota \le \| x(t,j) \|,
\end{equation}
and the claim of the proposition is proven. Consider the latter case, i.e., $\|x(t'',j+1)\|\geq r_\iota$. Then, the definition of the set $\mathcal{R}_\iota(\delta_\iota,\mu_\iota)$ in~\eqref{eq:Ri_deltai_mui} implies that $x(t'',j+1)\in\mathcal{R}_\iota(\delta_\iota,\mu_\iota)$. However, thanks to~\eqref{eq:Ri:disjoint}, the solution can only flow with stabilization mode while in $\mathcal{R}_\iota(\delta_\iota,\mu_\iota)$, so that $(t''',j+1)$ exists such that $\| x(t''',j+1)\| = r_\iota$. Define 
\begin{equation}
    \sigma_2:=\min_{i,i^\prime\in\mathbb{I},i \neq i^\prime} \mathbf{dist}(\mathcal{R}_i(\delta_i,\mu_i),\mathcal{E}_\leq(c_{i^\prime},\delta_{i^\prime} E_{i^\prime}))>0,
\end{equation}
which is positive because the considered sets are compact and pairwise sufficiently disjoint. Before a jump to avoidance mode at $(s,j+1)$ can occur, we have
\begin{equation}
\label{eq:decrease sigma 2}
\| x(s,j+1) \| + \sigma_2 \le \| x(t''',j+1) \| = r_\iota \le \| x(t,j) \|.
\end{equation}

\textit{Step 3c: $\xi(t,j)\in(\mathcal{F}_l \cap \mathcal{J}_l)\times\{l\}$ and the solution flows.}
The solution cannot flow forever, as established in Step~1. If it flows until its component $x$ reaches the set $\mathcal{P}_{l,3}^{i(t,j)}(\delta_{i(t,j)},\psi_{i(t,j)})$, the second part of the argument of Step~3b still applies, in particular \eqref{eq:decrease sigma 1} or \eqref{eq:decrease sigma 2}. If it jumps beforehand, Step~3a applies. Then we do not have circularity. By~combining \eqref{eq:decrease sigma 0}, \eqref{eq:decrease sigma 1} and \eqref{eq:decrease sigma 2}, \eqref{eq:strict decrease} is proven with $\sigma:=\min\{\sigma_0,\sigma_1,\sigma_2\}>0$.

\subsubsection{Proof of Lemma~\ref{lemma:tangent cones single subsets}}

As for~\eqref{eq:tangent cone: ell delta}, we have that $\mathbf{T}_{\mathcal{E}_\leq(c_i,\delta E_i)}(x)$ is either $\real^n$ for $x \in \mathcal{E}_<(c_i,\delta E_i)$ or $\mathcal{P}_\leq(0,E_i^2(x-c_i))$ for $x \in \mathcal{E}(c_i,\delta E_i)$. For all $x\in\mathcal{E}_\le(c_i,\delta E_i) \supset \tilde{\mathcal{F}}_m^i(\delta,\psi)$, $\mathfrak{u}(x)\in\mathcal{P}(0,E_i^2(x-c_i))$ (see \eqref{eq:vector field tangent cone |m|=1}), so \eqref{eq:tangent cone: ell delta} is proven. A similar argument yields \eqref{eq:tangent cone: ell}. As for~\eqref{eq:tangent cone: bar ell}, we note that 
\begin{equation}
\label{eq:tilde F minus tilde P}
\begin{aligned}
    &\tilde{\mathcal{F}}_m^i(\delta,\psi)\backslash\tilde{\mathcal{P}}^i_{m,3}(\delta,\psi) \!\!=\!\mathcal{E}_\leq(c_i,\!\delta E_i)\!\cap\!\mathcal{E}_\geq(c_i,\!E_i)\\
    &\hspace*{15pt} \! \cap\!\mathcal{E}_\geq(\bar c_i,\!\mu_i\bar E_i)\! \cap\! \mathcal{C}_\geq\!(c_i,c_i\!-\!p_m^i,\psi,E_i)\backslash\Big(\mathcal{E}_\leq(c_i,\delta E_i)\\
    &\hspace*{15pt} \! \cap\!\mathcal{E}_\geq\!(c_i,E_i) \!\cap\! \mathcal{E}(\bar c_i,\mu_i\bar E_i)\! \cap\!\mathcal{C}_\geq(c_i,c_i-p_m^i,\psi,E_i)\Big)\\
    &\overset{\eqref{eq:set identities:c}}{=}\mathcal{E}_\leq(c_i,\delta E_i)\cap\mathcal{E}_\geq(c_i,E_i)\cap\mathcal{E}_>(\bar c_i,\mu_i\bar E_i)\\
    &\hspace*{15pt} \cap \mathcal{C}_\geq(c_i,c_i-p_m^i,\psi,E_i).
    \end{aligned}
\end{equation}
Then, for all $x \in \tilde{\mathcal{F}}_m^i(\delta,\psi)\backslash \tilde{\mathcal{P}}^i_{m,3}(\delta,\psi)$, $\mathbf{T}_{\mathcal{E}_\geq(\bar c_i,\mu_i\bar E_i)}(x) = \real^n$ thanks to~\eqref{eq:tilde F minus tilde P}, and this proves \eqref{eq:tangent cone: bar ell}. 
As for \eqref{eq:tangent cone: cone}, we have that $\mathbf{T}_{\mathcal{C}_\geq(c_i,c_i-p_m^i,\psi,E_i)}(x)$ is either $\real^n$ for $x \in \mathcal{C}_>(c_i,c_i-p_m^i,\psi,E_i)$ or $\mathcal{P}_\geq(0,n^i_m(x))$ for $x \in \mathcal{C}(c_i,c_i-p_m^i,\psi,E_i)$ with $n_m^i(x):=E_i\pi^{\psi}(E_i(c_i-p_m^i))E_i(x-c_i)$ by~\eqref{eq:normal vector to cone surface}. If we prove that $\mathfrak{u}(x)^\top n_m^i(x) \ge 0$ for all $x \in \mathcal{C}(c_i,c_i-p_m^i,\psi,E_i) \supset \tilde{\mathcal{F}}_m^i(\delta,\psi)$, then \eqref{eq:tangent cone: cone} holds. Indeed, this last step follows from
\begin{equation}
    \begin{aligned}
       &\!\mathfrak{u}(x)^\top\!\cdot n_m^i(x)\!=\!-k_m(x-p_m^i)^\top\! E_i\pi^\perp(E_i(x-c_i))E_i^{-1}\\
       &\qquad\quad \cdot E_i\pi^{\psi}(E_i(c_i-p_m^i))E_i(x-c_i)\\
       &=k_m(p_m^i-c_i)^\top E_i\pi^\perp(E_i(x-c_i))\\
       &\qquad\quad \cdot \pi^{\psi}(E_i(c_i-p_m^i))E_i(x-c_i)\\
       &\overset{\eqref{eq:def:piTheta},\eqref{eq:proj-refl-maps}}{=}k_m(p_m^i-c_i)^\top E_i\pi^\perp(E_i(x-c_i))\\
       &\qquad\quad \cdot \big(\cos^2(\psi)I_n-\pi^\parallel(E_i(c_i-p_m^i))\big)E_i(x-c_i)\\
       &=-k_m(p_m^i-c_i)^\top E_i\pi^\perp(E_i(x-c_i))\\
       &\qquad\quad \cdot \pi^\parallel(E_i(c_i-p_m^i))E_i(x-c_i) \\
       &\overset{\eqref{eq:proj-refl-maps}}{=}k_m(c_i-p_m^i)^\top E_i\pi^\perp(E_i(x-c_i))E_i(c_i-p_m^i)\\
       &\qquad\quad \cdot (c_i - p_m^i)^\top\! E_i^2(x-c_i)\|E_i(c_i-p_m^i)\|^{-2}\! \ge 0
      \end{aligned}
\end{equation}
because for all $x\in\mathcal{C}(c_i,c_i-p_m^i,\psi,E_i)$, one has $x\in\mathcal{P}_\geq(c_i,E_i^2(c_i-p_m^i))$.

\bibliographystyle{IEEEtran}
\bibliography{IEEEabrv,references}
\begin{IEEEbiographynophoto}{Soulaimane Berkane}
received his Engineering and M.Sc. degrees in Automatic Control from Ecole Nationale Polytechnique, Algeria, in 2013, and his PhD in Electrical Engineering from the University of Western Ontario, Canada, in 2017. He held postdoctoral positions at the University of Western Ontario, Canada, and at KTH Royal Institute of Technology, Sweden, between 2018 and 2019. 
He is currently an assistant professor at the Department of Computer Science and Engineering, University of Quebec in Outaouais, Canada. His research interests are in the area of nonlinear control theory with applications to robotic and autonomous systems.\vspace{-0.8cm}
\end{IEEEbiographynophoto}
\begin{IEEEbiographynophoto}{Andrea Bisoffi} received the M.Sc. degree in Automatic Control Engineering from Politecnico di Milano, Italy, in 2013 
and the Ph.D. degree in Mechatronics from the University of Trento, Italy, in 2017. In 2015–-2016 he was a visiting scholar in 
the Control Group at the University of Cambridge, UK. 
He was a postdoctoral researcher at KTH Royal Institute of Technology, Sweden in 2017-2019 and is currently one at the University of Groningen, The Netherlands.
His current research interests include hybrid and nonlinear control systems, with applications to mechanical, and robotic systems.
\end{IEEEbiographynophoto}

\begin{IEEEbiographynophoto}{Dimos V. Dimarogonas}
received the Diploma in Electrical and Computer Engineering in 2001 and the Ph.D. in
Mechanical Engineering in 2007, both from the National
Technical University of Athens (NTUA), Greece. From 2007 to 2010, he held postdoctoral positions at KTH Royal Institute of Technology and Massachusetts Institute of Technology. He
is currently a Professor in Automatic Control,
School of Electrical Engineering and Computer Science, KTH Royal Institute of Technology. His current
research interests include multi-agent and hybrid systems with applications to autonomous systems. He serves on the Editorial Board of Automatica and the IEEE Transactions on Control of Network Systems. Dr. Dimarogonas is a recipient of an ERC Starting Grant in 2014, an ERC Consolidator Grant in 2019, and a Wallenberg Academy Fellowship in 2015. He is a member of the Technical Chamber of Greece and a Senior Member of the IEEE.
\end{IEEEbiographynophoto}

\end{document}